\documentclass[10pt]{amsart}

\pdfoutput=1

\usepackage[text={420pt,660pt},centering]{geometry}
\makeatletter
\@mparswitchfalse \makeatother

\usepackage{cancel}
\usepackage{graphicx}
\IfFileExists{MnSymbol.sty}{\usepackage{MnSymbol}}{\usepackage{amssymb}}
\usepackage{mathtools}
\usepackage[colorlinks=true, pdfstartview=Fit, linkcolor=blue, citecolor=blue, urlcolor=blue,pagebackref=false]{hyperref}
\usepackage{orcidlink}
\usepackage{microtype}

\usepackage{bm}
\IfFileExists{dsfont.sty}{\usepackage{dsfont}}{}
\providecommand{\mathds}[1]{\mathbf{#1}}
\usepackage{mathrsfs}
\usepackage{xcolor}
\IfFileExists{accents.sty}{\usepackage{accents}}{}

\IfFileExists{ulem.sty}{\usepackage[normalem]{ulem}}{}

\definecolor{darkgreen}{rgb}{0,0.5,0}
\definecolor{darkblue}{rgb}{0,0,0.7}
\definecolor{darkred}{rgb}{0.9,0.1,0.1}

\newtheorem{proposition}{Proposition}
\newtheorem{theorem}[proposition]{Theorem}
\newtheorem{lemma}[proposition]{Lemma}
\newtheorem{corollary}[proposition]{Corollary}

\theoremstyle{remark}

\theoremstyle{definition}

\numberwithin{equation}{section}
\numberwithin{proposition}{section}
\numberwithin{figure}{section}
\numberwithin{table}{section}

\newcommand{\N}{\mathbb{N}}
\newcommand{\Q}{\mathbb{Q}}
\newcommand{\R}{\mathbb{R}}

\newcommand{\E}{\mathbb{E}}
\renewcommand{\P}{\mathbb{P}}

\renewcommand{\S}{\mathbf{S}}

\newcommand{\ep}{\varepsilon}

\renewcommand{\le}{\leqslant}
\renewcommand{\ge}{\geqslant}

\renewcommand{\geq}{\geqslant}

\renewcommand{\subset}{\subseteq}
\renewcommand{\bar}{\overline}

\newcommand{\Ll}{\left}
\newcommand{\Rr}{\right}
\renewcommand{\d}{\mathrm{d}}

\newcommand{\mcl}{\mathcal}

\DeclareMathOperator{\tr}{tr}
\DeclareMathOperator{\supp}{supp}

\newcommand{\sP}{\mathscr{P}}

\newcommand{\la}{\left\langle}
\newcommand{\ra}{\right\rangle}

\newcommand{\one}{\mathds{1}}

\title{Free energy of spherical vector spin glasses}
\date{}
\hypersetup{pdftitle={Free energy of spherical vector spin glasses}}
\begin{document}
\author[H.-B.\ Chen]{Hong-Bin Chen\,\orcidlink{0000-0001-6412-0800}}
\address[Hong-Bin Chen]{NYU-ECNU Institute of Mathematical Sciences, NYU Shanghai, China}
\email{\href{mailto:hongbin.chen@nyu.edu}{hongbin.chen@nyu.edu}}
\raggedbottom \begin{abstract}
We identify the limit free energy of spherical vector spin glasses by a variational formula. The covariance function need not be convex. Using the global symmetry of the spin configurations, we introduce a family of constraint measures that gives a disintegration of the original spin distribution. We first identify the limit free energy of each constrained model and then apply the Laplace method to recover that of the original model.
\end{abstract}
\maketitle
\setcounter{tocdepth}{2}

\section{Introduction}
\label{s:introduction}

We study the limit free energy of spherical vector spin glasses, continuing the investigation in \cite{panchenkoTalagrand2007overlap,ko2020multiple,ko2019crisanti}. We allow non-convex interactions and identify the limit by a variational formula.

The following notation is used throughout. Write $\N=\{1,2,\dots\}$ and fix $D\in\N$. We write $\S^D$ for the real symmetric $D\times D$ matrices, and $\S_+^D$ and $\S_{++}^D$ for their positive-semidefinite and positive-definite cones. For $a,b\in\S^D$, we write $a\preccurlyeq b$ if $b-a\in\S_+^D$. For vectors, matrices, or tensors $a,b$ of the same dimensions, we write $a\cdot b:=\sum_{\bm i}a_{\bm i}b_{\bm i}$ and $|a|:=(a\cdot a)^{1/2}$. The matrix transpose is denoted by $\intercal$.

For $k\in\N$, let $\Delta_k\in\S_+^{D^k}$ act on $(\R^D)^{\otimes k}$ in its standard tensor-product basis, and assume that
\begin{equation}
\label{eq:mixture-summability}
\sum_{k\ge1}\rho^k|\Delta_k|<\infty,\qquad\forall\rho>0.
\end{equation}
Define $\xi:\R^{D\times D}\to\R$ by
\begin{equation}
\label{eq:xi-mixture}
\xi(r):=\sum_{k\ge1}\langle\Delta_k,r^{\otimes k}\rangle.
\end{equation}
Assumption \eqref{eq:mixture-summability} is equivalent to absolute convergence of \eqref{eq:xi-mixture} everywhere. The assumption $\Delta_k\in\S_+^{D^k}$ allows the construction of $H_N$ below but does not imply convexity of $\xi$.

For $N\in\N$, let $S_N:=\{x\in\R^N:|x|^2=N\}$, let $\mu_N$ be normalized surface measure on $S_N$, and put $P_N:=\mu_N^{\otimes D}$ on $S_N^D$. We view configurations $\sigma\in S_N^D$ as $D\times N$ matrices and write $\sigma_i\in\R^D$ for the $i$th column. Thus, under $P_N$, the rows of $\sigma$ are independent with common law $\mu_N$. For $\sigma,\sigma'\in\R^{D\times N}$, their \textbf{overlap} is $\frac1N\sigma\sigma'^\intercal\in\R^{D\times D}$. 
Let $H_N$ be the centered Gaussian field on $\R^{D\times N}$ with covariance (see~\eqref{eq:mixed-Hamiltonian} for its construction)
\begin{equation}\label{e.EHH=xi}
\E H_N(\sigma)H_N(\sigma')=N\xi\Big(\frac1N\sigma\sigma'^\intercal\Big).
\end{equation}
Fix $h\in\R^D$ throughout the paper. We consider the free energy
\begin{equation}
\label{eq:physical-free-energy}
\mathsf F_N:=-\frac1N\E\log\int_{S_N^D}\exp\Big(H_N(\sigma)+h\cdot\sum_{i=1}^N\sigma_i\Big)P_N(\mathrm d\sigma).
\end{equation}

We enrich this free energy by Gaussian fields whose covariance is parametrized by a matrix path. Let $\mcl Q$ be the set of right-continuous increasing paths $q:[0,1)\to\S_+^D$ with respect to $\preccurlyeq$. For $\rho\in[1,\infty]$, put $\mcl Q_\rho:=\mcl Q\cap L^\rho([0,1),\S^D)$. For $r\in\S_+^D$, define
\begin{equation}\label{e.Q_infty,<=}
\mcl Q_{\infty,\preccurlyeq r}:=\{p\in\mcl Q_\infty:p(s)\preccurlyeq r\text{ for every }s\in[0,1)\}.
\end{equation}
For $q\in\mcl Q_\infty$, write $q(1):=\lim_{s\uparrow1}q(s)$. The path pairing is $\langle q,p\rangle_{L^2}:=\int_0^1q(s)\cdot p(s)\,\mathrm d s$. Let $\mathfrak R$ be the continuous Ruelle probability cascade with uniform overlap distribution. For $q\in\mcl Q_\infty$, conditionally on $\mathfrak R$, let $z_i^q:\supp\mathfrak R\to\R^D$, $i\in\{1,\ldots,N\}$, be independent centered Gaussian fields, independent of $H_N$, with covariance $\E\big[z_i^q(\alpha)z_j^q(\alpha')^{\intercal}\mid\mathfrak R\big]=\one_{\{i=j\}}q(\alpha\wedge\alpha')$. We refer to Section~\ref{s:enriched-properties} for more details. For $t\ge0$ and $q\in\mcl Q_\infty$, define the enriched free energy by
\begin{align}
\bar F_N(t,q):={}&-\frac1N\E\log\iint_{S_N^D\times\supp\mathfrak R}\exp\Big(\sqrt{2t}H_N(\sigma)-Nt\xi\Ll(N^{-1}\sigma\sigma^{\intercal}\Rr)\notag\\
&\qquad+\sqrt2\sum_{i=1}^Nz_i^q(\alpha)\cdot\sigma_i-q(1)\cdot\sigma\sigma^{\intercal}+h\cdot\sum_{i=1}^N\sigma_i\Big)P_N(\mathrm d\sigma)\mathfrak R(\mathrm d\alpha).
\label{eq:global-enriched-intro}
\end{align}
The expectation averages the Gaussian fields and the cascade. The free energy $\bar F_N(t,\cdot)$ admits a unique $L^1$-continuous extension to $\mcl Q_1$ by \eqref{eq:global-enriched-q-Lipschitz}. 

Our variational formula uses a constrained model with prescribed \textbf{self-overlap} $N^{-1}\sigma\sigma^\intercal$. In \cite{ko2020multiple,ko2019crisanti}, the product surface measure is restricted to an $\varepsilon$-neighborhood of the prescribed self-overlap. Here we use a probability measure on configurations with exactly that self-overlap, as described in \eqref{eq:stiefel-intro}--\eqref{eq:nu-pushforward}. The constraint is parametrized by the sets
\begin{equation}
\label{eq:elliptope}
\mcl A:=\{a\in\S_+^D:(a^2)_{dd}=1\text{ for }d\in\{1,\dots,D\}\}\quad\text{and}\quad \mcl A^\circ:=\mcl A\cap\S_{++}^D.
\end{equation}
For $q\in\mcl Q_1$ and $\ell\in\S_{++}^D$, define $L_{\ell,q}:[0,1)\to\S_{++}^D$ by
\begin{equation}
\label{eq:initial-path-multiplier}
L_{\ell,q}(s):=\ell+2\Big(sq(s)-\int_0^s q(r)\,\mathrm d r\Big),\qquad s\in[0,1).
\end{equation}
Since $q$ is increasing, we have $L_{\ell,q}(s)\succcurlyeq\ell\succ0$. For each $v\in\R^D$, define $\psi_v:\mcl Q_1\to\R$ by 
\begin{align}
\label{eq:psi-general-path}
\psi_v(q):=\int_0^1\tr q(s)\,\mathrm d s-\frac12\min_{\ell\in\S_{++}^D}\bigg\{&\tr\ell-D-\log\det\ell\notag\\
&\qquad+2\tr(\ell^{-1}q(0))+v^{\intercal}\ell^{-1}v+\int_0^1\frac1{s^2}\log\frac{\det L_{\ell,q}(s)}{\det\ell}\,\mathrm d s\bigg\}.
\end{align}
Lemma~\ref{cor:initial-general-path} proves that the integral is finite and the minimum is attained. For $a\in\mcl A^\circ$, define the induced covariance function $\xi_a:\R^{D\times D}\to\R$ by
\begin{equation}
\label{eq:fixed-a-xi-intro}
\xi_a(r):=\xi(ara).
\end{equation}
For $a\in\mcl A^\circ$, $t\ge0$, and $q\in\mcl Q_2$, define $\mcl J_{a,t,q}:\mcl Q_2\times\mcl Q_{\infty,\preccurlyeq I_D}\to\R$ by
\begin{equation}
\label{eq:J-fixed-a-intro}
\mcl J_{a,t,q}(q',p):=\psi_{ah}(q')+\langle q-q',p\rangle_{L^2}+t\int_0^1\xi_a(p(s))\,\mathrm d s.
\end{equation}
The corresponding Hopf formula $f_a:[0,\infty)\times\mcl Q_2\to\R$ is defined by
\begin{equation}
\label{eq:Hopf-fixed-a-intro}
f_a(t,q):=\sup_{p\in\mcl Q_{\infty,\preccurlyeq I_D}}\inf_{q'\in\mcl Q_2}\mcl J_{a,t,q}(q',p).
\end{equation}
For $a\in\mcl A^\circ$ and $q\in\mcl Q$, define the transformed path $q_a:[0,1)\to\S_+^D$ by
\begin{equation}
\label{eq:qa-intro}
q_a(s):=aq(s)a,\qquad s\in[0,1).
\end{equation}

\begin{theorem}
\label{thm:main}
We have
\begin{equation*}
\lim_{N\to\infty}\mathsf F_N=\inf_{a\in\mcl A^\circ}\Big\{f_a\Big(\tfrac12,0\Big)-\tfrac12\xi(a^2)-\log\det a\Big\}.
\end{equation*}
Moreover, for every $t\ge0$ and $q\in\mcl Q_2$, we have
\begin{equation*}
\lim_{N\to\infty}\bar F_N(t,q)=\inf_{a\in\mcl A^\circ}\Big\{f_a(t,q_a)-\log\det a\Big\}.
\end{equation*}
Both infima above are attained. 
\end{theorem}

The proof is given at the end of Section~\ref{s:positive-outer}.

\smallskip
\noindent\textbf{Outline of proof.}
We first identify the limit at each fixed positive-definite self-overlap and then optimize over this constraint. Our choice of a normalized measure on the constraint has two useful properties. Its limit initial condition is convex in the matrix path, which permits the Hamilton--Jacobi approach used in \cite{chenIssaMourrat2026Ising,chenMourrat2026spherical}. It also gives an exact disintegration of $P_N$, from which the determinant term in Theorem~\ref{thm:main} arises. 

For $a\in\mcl A^\circ$, the configurations with self-overlap $a^2$ are parametrized by $\sigma=au$, where
\begin{equation}
\label{eq:stiefel-intro}
\mathbb U_N:=\{u\in\R^{D\times N}:uu^{\intercal}=NI_D\},\qquad N\ge D.
\end{equation}
Write $\operatorname{Orth}(N):=\{o\in\R^{N\times N}:oo^{\intercal}=I_N\}$ for the orthogonal group, and let $O_N$ be Haar distributed on $\operatorname{Orth}(N)$. Define the probability measure $\nu_N$ on $\mathbb U_N$ by
\begin{equation}
\label{eq:nu-pushforward}
\nu_N:=\operatorname{Law}\big(\sqrt N\,(I_D\ 0)O_N\big),
\end{equation}
where $(I_D\ 0)$ is the $D\times N$ matrix whose first $D$ columns form $I_D$ and whose remaining entries are zero.
For $u,u'\in\mathbb U_N$, the covariance identity \eqref{e.EHH=xi} gives
\begin{equation}
\label{eq:fixed-a-covariance-intro}
\E H_N(au)H_N(au')=N\xi_a(N^{-1}uu'^{\intercal}).
\end{equation}
For $t\ge0$ and $q\in\mcl Q_\infty$, the constrained free energy is defined by
\begin{align}
    \bar F_{N,a}(t,q):=-\frac1N\E\log \iint_{\mathbb U_N\times\supp\mathfrak R}&\exp\Big(\sqrt{2t}H_N(au)-Nt\xi_a(I_D)\notag
    \\
    &+\sqrt2\sum_{i=1}^Nz_i^q(\alpha)\cdot u_i-Nq(1)\cdot I_D+ah\cdot\sum_{i=1}^Nu_i\Big)\nu_N(\mathrm d u)\mathfrak R(\mathrm d\alpha). \label{e.F_N,a_intro}
\end{align}
Lemma~\ref{cor:initial-general-path} and \eqref{eq:enriched-fixed-a-initial} identify the function in \eqref{eq:psi-general-path} as
\begin{equation*}
\lim_{N\to\infty}\bar F_{N,a}(0,q)=\psi_{ah}(q),\qquad a\in\mcl A^\circ,\quad q\in\mcl Q_1.
\end{equation*}
Proposition~\ref{prop:stiefel-convexity} proves that $\psi_{ah}$ is convex on $\mcl Q_1$. This is the convexity needed in the argument below, and it does not require convexity of $\xi$. 

Fix $a\in\mcl A^\circ$, $t\ge0$, and, for this outline, $q\in\mcl Q_\infty$. The Hamilton--Jacobi lower bound of \cite[Theorem~3.4]{mourrat2023upper}, together with the Hopf representation in \cite[Theorem~4.7(3)]{chenXia2025HJ}, gives the following inequality, proved in Proposition~\ref{prop:HJ-lower},
\begin{equation*}
f_a(t,q)\le\liminf_{N\to\infty}\bar F_{N,a}(t,q).
\end{equation*}
For the reverse inequality, Proposition~\ref{prop:block-critical} gives $p^+\in\mcl Q_{\infty,\preccurlyeq I_D}$ such that, with $q^+:=q+t\nabla\xi_a(p^+)$, we have 
\begin{equation*}
p^+\in\partial\psi_{ah}(q^+)\quad\text{and}\quad\lim_{N\to\infty}\bar F_{N,a}(t,q)=\mcl J_{a,t,q}(q^+,p^+).
\end{equation*}
Here $\partial\psi_{ah}$ is the subdifferential relative to $\mcl Q_2$, defined in \eqref{eq:initial-subdifferential}. We explain first how the constraint is handled in this representation and then how convexity gives the upper bound.

The reference measure $\nu_N$ does not factor over the columns of $u$, so the cavity results (e.g.\ those in~\cite{chen2025free}) for vector spins with product reference measures do not apply directly. Section~\ref{s:block-approximation} instead compares $\bar F_{N,a}$ with $\bar F_{N,a}^{[M]}$. For $N=nM$, the latter is obtained by replacing $\nu_N$ with $\nu_M^{\otimes n}$, as defined at the beginning of Section~\ref{s:block-restriction}. Thus $M$ is the number of columns in each block, and a block can be viewed as one vector spin in $\R^{DM}$. The same section treats the remaining block when $M$ does not divide $N$.

The column-rotation estimate in Lemma~\ref{lem:block-moment} and the polar normalization in \eqref{eq:block-normalization-map} lead to Lemma~\ref{lem:block-comparison}, which gives
\begin{equation*}
\limsup_{N\to\infty}\left|\bar F_{N,a}(t,q)-\bar F_{N,a}^{[M]}(t,q)\right|\le CM^{-1/4}.
\end{equation*}
For fixed $M$, the convergence results for compact product measures \cite[Theorem~1.1]{ho2026convergence} and \cite[Theorem~1.1]{chenHoIssa2026formula} apply to the block model. Letting $N\to\infty$ first and then $M\to\infty$ proves convergence of $\bar F_{N,a}(t,q)$ in Proposition~\ref{prop:block-limit}. Lemma~\ref{lem:block-product-critical} applies the critical-point formula of \cite[Proposition~7.3]{chen2025free} to each block model. To pass its critical relation to the limit $M\to\infty$, Section~\ref{s:block-critical} uses the uniform semiconcavity estimate in Lemma~\ref{lem:block-semiconcavity}, in addition to convergence of the function values. This proves the representation of $\lim_{N\to\infty}\bar F_{N,a}(t,q)$ above.

The final convexity argument is given in Corollary~\ref{cor:fixed-a-limit}. The inclusion $p^+\in\partial\psi_{ah}(q^+)$ implies, for every $q'\in\mcl Q_2$, that 
\begin{equation*}
\mcl J_{a,t,q}(q',p^+)-\mcl J_{a,t,q}(q^+,p^+)=\psi_{ah}(q')-\psi_{ah}(q^+)-\langle p^+,q'-q^+\rangle_{L^2}\ge0.
\end{equation*}
Thus $q^+$ attains $\inf_{q'\in\mcl Q_2}\mcl J_{a,t,q}(q',p^+)$, and we have
\begin{equation*}
\lim_{N\to\infty}\bar F_{N,a}(t,q)=\mcl J_{a,t,q}(q^+,p^+)=\inf_{q'\in\mcl Q_2}\mcl J_{a,t,q}(q',p^+)\stackrel{\eqref{eq:Hopf-fixed-a-intro}}{\le} f_a(t,q).
\end{equation*}
Combining this with the lower bound proves $\lim_{N\to\infty}\bar F_{N,a}(t,q)=f_a(t,q)$. 

We now pass from the constrained model to Theorem~\ref{thm:main}. The disintegration in Section~\ref{s:outer-a} uses the self-overlap $b=a^2$ as the integration variable. Define
\begin{equation*}
\mcl B^\circ:=\{b\in\S_{++}^D:b_{dd}=1\text{ for }d\in\{1,\ldots,D\}\}=\{a^2:a\in\mcl A^\circ\},
\end{equation*}
and let $\mathrm d b_{\mathrm{off}}:=\prod_{1\le d<d'\le D}\mathrm d b_{dd'}$ denote Lebesgue measure in the $D(D-1)/2$ independent off-diagonal entries $(b_{dd'})_{d<d'}$. Lemma~\ref{lem:disintegration} gives, for $N>D$ and every nonnegative measurable $G:S_N^D\to[0,\infty]$, the identity
\begin{equation*}
\int_{S_N^D}G(\sigma)P_N(\mathrm d\sigma)=c_{N,D}\iint_{\mathbb U_N\times\mcl B^\circ}G(b^{1/2}u)\nu_N(\mathrm d u)(\det b)^{(N-D-1)/2}\,\mathrm d b_{\mathrm{off}}.
\end{equation*}
The normalizing constant $c_{N,D}$ is given in \eqref{eq:disintegration-constant}. Lemma~\ref{lem:matrix-density} proves that 
\begin{equation*}
\lim_{N\to\infty}\frac1N\log c_{N,D}=0.
\end{equation*}

For $q\in\mcl Q_\infty$, let $Z_N(t,q)$ and $Z_{N,a}(t,q)$ be the associated partition functions such that $\bar F_N(t,q)=-\frac{1}{N}\E\log Z_N(t,q)$ and $\bar F_{N,a}(t,q)=-\frac{1}{N}\E\log Z_{N,a}(t,q)$. 
In the following identity, the constrained field with path $q_a$ is realized as $az_i^q$. Lemma~\ref{lem:global-enriched-disintegration} gives, almost surely, 
\begin{equation*}
Z_N(t,q)=c_{N,D}\int_{\mcl B^\circ}Z_{N,b^{1/2}}(t,q_{b^{1/2}})(\det b)^{(N-D-1)/2}\,\mathrm d b_{\mathrm{off}},
\end{equation*}
where $q_{b^{1/2}}$ is given as in~\eqref{eq:qa-intro} with $a=b^{1/2}$.
Since $\det b=(\det a)^2$, the determinant density contributes the term $-\log\det a$ to the limit free energy. For large $N$, the Laplace principle therefore suggests the heuristic
\begin{equation*}
\bar F_N(t,q)\approx\inf_{a\in\mcl A^\circ}\{\bar F_{N,a}(t,q_a)-\log\det a\}.
\end{equation*}
To justify this approximation, we need uniform control of the random constrained free energies before taking the logarithm and expectation. For the random free energy $F_{N,a}(t,q):=-N^{-1}\log Z_{N,a}(t,q)$, Lemma~\ref{lem:enriched-local-uniform} proves, for every $\kappa\in(0,1)$, that
\begin{equation*}
\lim_{N\to\infty}\E\sup_{a\in\mcl A^\circ:\,\det(a^2)\ge\kappa}|F_{N,a}(t,q_a)-f_a(t,q_a)|=0.
\end{equation*}
This estimate justifies replacing $F_{N,a}(t,q_a)$ by $f_a(t,q_a)$ uniformly on the compact set $\{a\in\mcl A^\circ:\det(a^2)\ge\kappa\}$ at a vanishing cost in expectation. Lemma~\ref{lem:enriched-boundary} allows $\kappa>0$ to be chosen so that matrices with $\det(a^2)<\kappa$ cannot minimize the objective. These estimates, together with the Laplace calculation in Proposition~\ref{prop:global-enriched}, yield
\begin{equation*}
\lim_{N\to\infty}\bar F_N(t,q)=\inf_{a\in\mcl A^\circ}\{f_a(t,q_a)-\log\det a\}.
\end{equation*}

Finally, for $N>D$, \eqref{eq:measure-disintegration} gives the following identity for the partition function of $\mathsf F_N$,
\begin{align*}
    \int_{S_N^D}\exp\Big(H_N(\sigma)+h\cdot\sum_{i=1}^N\sigma_i\Big)P_N(\mathrm d\sigma)=c_{N,D}\int_{\mcl B^\circ}e^{N\xi(b)/2}Z_{N,b^{1/2}}\Big(\tfrac12,0\Big)(\det b)^{(N-D-1)/2}\,\mathrm d b_{\mathrm{off}}.
\end{align*}
At $t=1/2$ and $q=0$, the correction in \eqref{eq:enriched-fixed-a-H} is $-N\xi_{b^{1/2}}(I_D)/2=-N\xi(b)/2$. The factor $e^{N\xi(b)/2}$ removes this correction. Therefore, the corresponding Laplace heuristic is
\begin{equation*}
\mathsf F_N\approx\inf_{a\in\mcl A^\circ}\Big\{\bar F_{N,a}\Big(\tfrac12,0\Big)-\tfrac12\xi(a^2)-\log\det a\Big\}.
\end{equation*}
Section~\ref{s:positive-outer} applies the same disintegration and Laplace argument to obtain the formula for $\mathsf F_N$ in Theorem~\ref{thm:main}, including attainment.

\smallskip\noindent\textbf{Related works.} 
The Parisi formula for the Sherrington--Kirkpatrick model was established through the interpolation bound in \cite{guerra2003broken} and the matching bound in \cite{talagrand2006parisi}. Variational formulas were subsequently obtained for multi-species Sherrington--Kirkpatrick models in \cite{panchenko2015multispecies}, the Potts spin glass in \cite{panchenko2018potts}, and models with vector spins in \cite{panchenko2018vector}.

For scalar spherical spins, the free-energy formula was proved for even $p$-spin interactions in \cite{talagrand2006spherical} and extended to include odd $p$ in \cite{chen2013aizenman}. The Parisi formula for multi-species spherical models was obtained in \cite{batesSohn2022free} under a convexity assumption on the covariance function. In these multi-species models, the covariance depends on one scalar overlap for each species. The covariance in \eqref{e.EHH=xi} may instead depend on every entry of $N^{-1}\sigma\sigma'^{\intercal}$.

The closest earlier results involving full matrix overlaps concern several spherical systems coupled by a constraint on their mutual overlaps. A Guerra interpolation bound for these constrained systems was derived in \cite{panchenkoTalagrand2007overlap}. The matching bound was proved in \cite{ko2020multiple} for copies of spherical models with even $p$-spin interactions, and the corresponding Crisanti--Sommers formula was obtained in \cite{ko2019crisanti}. For these copies, the covariance is a sum of even powers of individual overlap entries. The tensor expansion \eqref{eq:xi-mixture} also allows products of different entries and need not be convex. Moreover, the free energies in those works retain the mass of a neighborhood of the prescribed self-overlap. Our normalized constraint measure separates this contribution from the constrained free energy, and the term $-\log\det a$ restores it in Theorem~\ref{thm:main}.

The Hamilton--Jacobi approach provides a way to study such non-convex covariance functions. A bound by the solution of an infinite-dimensional Hamilton--Jacobi equation was proved in \cite{mourrat2023upper}. A critical-point representation of the limit was obtained in \cite{chen2025free} under the assumption that the enriched free energy converges. Convexity of the initial condition turns this representation into a Hopf formula in the centered Ising setting of \cite{chenIssaMourrat2026Ising} and in the multi-species spherical setting of \cite{chenMourrat2026spherical}. In the latter case, this convexity holds for every deterministic external field. Our proof follows this strategy at each fixed self-overlap, with convexity on the cone of matrix-valued paths supplied by Proposition~\ref{prop:stiefel-convexity}. The variable self-overlap then requires the additional optimization over $a\in\mcl A^\circ$.

For compactly supported vector spins with a product reference measure, convergence without assuming convexity of the covariance was proved in \cite{ho2026convergence}. The limit was subsequently identified in \cite{chenHoIssa2026formula} by minimizing, over the empirical spin distribution, the sum of a relative entropy and a constrained Hamilton--Jacobi solution. These results apply directly to the product block models used in Section~\ref{s:block-approximation}. Passing from these blocks to $\nu_N$ requires the comparison in Lemma~\ref{lem:block-comparison}, since $\nu_N$ is not a product over the columns. 

Non-convex multi-species Ising models with arbitrary deterministic external fields were treated in \cite{boursierMourrat2026external} by first fixing the magnetization of each species. The resulting initial condition is convex in the path, and optimizing over the magnetizations gives a variational formula over Brownian martingales. This shares with \cite{chenHoIssa2026formula} and the present paper the strategy of identifying constrained free energies before optimizing over the constraint. The constraints are different: magnetization in \cite{boursierMourrat2026external}, empirical spin distribution in \cite{chenHoIssa2026formula}, and the full self-overlap $a^2$ here. In our setting, disintegration of the product surface measure gives the explicit cost $-\log\det a$, while the constrained limit is the Hopf formula \eqref{eq:Hopf-fixed-a-intro} with the explicit initial condition \eqref{eq:psi-general-path}.

The constrained free-energy formulas also have applications to ground states and the structure of their variational minimizers. For the coupled spherical models considered above, \cite{auffingerZhou2022properties} derives a variational formula for the ground-state energy by taking the zero-temperature limit of the Crisanti--Sommers formula. The same work studies the minimizers at positive temperature and gives a criterion for replica symmetry at zero temperature.

\smallskip\noindent\textbf{Organization of the paper.}
Section~\ref{s:enriched-initial} constructs the enriched model and proves the spherical initial-condition formula, its convexity, and the concentration estimates. Section~\ref{s:block-approximation} compares the model with reference measure $\nu_N$ on $\mathbb U_N$ against products of the measures $\nu_M$ on $\mathbb U_M$. Section~\ref{s:block-critical} proves the representation $\lim_{N\to\infty}\bar F_{N,a}(t,q)=\mcl J_{a,t,q}(q^+,p^+)$. Section~\ref{s:fixed-a-free-energy} identifies $\lim_{N\to\infty}\bar F_{N,a}=f_a$ by combining the Hamilton--Jacobi lower bound with the previous representation. Section~\ref{s:main-proof} applies the Laplace method to optimize over $a$.

\smallskip\noindent\textbf{AI use.}
ChatGPT 5.6 was used to identify a new constraint argument that makes the approach of~\cite{chenIssaMourrat2026Ising,chenMourrat2026spherical} applicable. ChatGPT 6 assisted with the writing of the paper. The author has revised the manuscript to improve its clarity and readability throughout.

\smallskip\noindent\textbf{Acknowledgements.} HBC acknowledges funding from the NYU Shanghai Start-Up Fund and support from the NYU--ECNU Institute of Mathematical Sciences at NYU Shanghai.

\section{Enriched model and initial condition}
\label{s:enriched-initial}

Proposition~\ref{prop:stiefel-initial} evaluates the initial condition for finite-step paths, and Lemma~\ref{cor:initial-general-path} identifies $\lim_{N\to\infty}\bar F_{N,a}(0,q)=\psi_{ah}(q)$ for every $q\in\mcl Q_1$. Proposition~\ref{prop:stiefel-convexity} proves the convexity required in Proposition~\ref{prop:HJ-lower}. The proof of Proposition~\ref{prop:stiefel-initial} evaluates an integral over a Gaussian matrix $\omega\in\R^{D\times N}$ with a matrix multiplier, then restricts $N^{-1}\omega\omega^{\intercal}$ to a neighborhood of $I_D$. The sampled-branch law and variance estimates of \cite[Corollary~A.2 and Lemma~A.4]{chenMourrat2026spherical} are used in the initial-condition computation. Lemma~\ref{lem:total-concentration-continuous} controls fluctuations over the Gaussian fields and cascade weights.

\subsection{Construction of \texorpdfstring{$H_N$}{H\_N} and Haar invariance}
\label{s:mixed-model}
\label{s:fixed-a}

For $k\in\N$, $(i_1,\ldots,i_k)\in\{1,\ldots,N\}^k$, and $\bm d\in\{1,\dots,D\}^k$, let $g^{(k)}_{i_1,\ldots,i_k,\bm d}$ be centered jointly Gaussian random variables. The vectors indexed by $\bm d$ are independent for distinct $(k,i_1,\ldots,i_k)$, and their covariance is
\begin{equation}
\label{eq:g-covariance}
\E g^{(k)}_{i_1,\ldots,i_k,\bm d}g^{(k)}_{j_1,\ldots,j_k,\bm d'} =\one_{\{(i_1,\ldots,i_k)=(j_1,\ldots,j_k)\}}(\Delta_k)_{\bm d,\bm d'}.
\end{equation}
For $\sigma\in\R^{D\times N}$, define
\begin{equation}
\label{eq:mixed-Hamiltonian}
H_N(\sigma):=\sum_{k\ge1}\frac1{N^{(k-1)/2}} \sum_{i_1,\ldots,i_k\le N}\sum_{\bm d\in\{1,\dots,D\}^k} g^{(k)}_{i_1,\ldots,i_k,\bm d}\prod_{j=1}^k\sigma_{d_ji_j}.
\end{equation}

The degree-$k$ tensor construction is described in \cite[Section~6.3]{mourrat2023upper}. For each fixed $\sigma\in\R^{D\times N}$, the variance of the sum of terms in \eqref{eq:mixed-Hamiltonian} with $K<k\le L$ is at most $N\sum_{k>K}|\frac1N\sigma\sigma^{\intercal}|^k|\Delta_k|$. Assumption \eqref{eq:mixture-summability} makes this bound tend to zero as $K\to\infty$, so the partial sums converge in $L^2$. Contracting the tensor indices using \eqref{eq:g-covariance} and \eqref{eq:xi-mixture} gives \eqref{e.EHH=xi} for every $\sigma,\sigma'\in\R^{D\times N}$. 

\begin{lemma}[Properties of $\xi$]
\label{lem:mixture-monotonicity}
For every $r\in\R^{D\times D}$, we have $\xi(r^{\intercal})=\xi(r)$. Consequently, if $r\in\S^D$, then $\nabla\xi(r)\in\S^D$. If $r,s\in\S_+^D$ satisfy $r\succcurlyeq s$, then we have $\nabla\xi(r)\succcurlyeq\nabla\xi(s)\succcurlyeq0$.
\end{lemma}

\begin{proof}
Since $\Delta_k$ is symmetric as an operator on $(\R^D)^{\otimes k}$, we have $\xi(r^{\intercal})=\xi(r)$ from~\eqref{eq:xi-mixture}. Differentiating $\xi(r^{\intercal})=\xi(r)$ at a symmetric matrix gives $\nabla\xi(r)=\nabla\xi(r)^{\intercal}$.

For $r\succcurlyeq s\succcurlyeq0$, we have
\begin{equation}
\label{eq:tensor-telescoping}
r^{\otimes k}-s^{\otimes k} =\sum_{j=1}^ks^{\otimes(j-1)}\otimes(r-s)\otimes r^{\otimes(k-j)}.
\end{equation}
For $c\succcurlyeq0$, differentiating the $k$th term of \eqref{eq:xi-mixture} in direction $c$ gives a sum of $k$ terms containing one factor $c$ and $k-1$ factors $r$. Each term is nonnegative, since it pairs two positive-semidefinite operators. Applying \eqref{eq:tensor-telescoping} to the $k-1$ factors shows that this directional derivative at $r$ is at least its value at $s$. The derivative series converge uniformly on bounded sets by \eqref{eq:mixture-summability}, so summing over $k$ gives $c\cdot\nabla\xi(r)\ge c\cdot\nabla\xi(s)\ge0$. Since $\S_+^D$ is self-dual and the gradients are symmetric, we obtain $\nabla\xi(r)\succcurlyeq\nabla\xi(s)\succcurlyeq0$.
\end{proof}

We next record the invariance of the reference measure and identify the constrained model in the coordinates of $\mathbb U_N$.

Recall $\mcl A^\circ$ from \eqref{eq:elliptope} and fix $a\in\mcl A^\circ$. For $N\ge D$, the set $\mathbb U_N=\{u\in\R^{D\times N}:uu^{\intercal}=NI_D\}$ carries the probability measure $\nu_N=\operatorname{Law}(\sqrt N\,(I_D\ 0)O_N)$, where $O_N$ is Haar distributed on $\operatorname{Orth}(N)$, as in \eqref{eq:nu-pushforward}.

\begin{lemma}[Invariance and uniqueness of $\nu_N$]
\label{lem:stiefel-invariance}
For every bounded measurable $\varphi:\mathbb U_N\to\R$ and every $o\in\operatorname{Orth}(N)$, we have
\begin{equation}
\label{eq:stiefel-right-invariance}
\int_{\mathbb U_N}\varphi(uo)\,\nu_N(\d u)=\int_{\mathbb U_N}\varphi(u)\,\nu_N(\d u).
\end{equation}
Moreover, $\nu_N$ is the unique probability measure on $\mathbb U_N$ invariant under $u\mapsto uo$ for every $o\in\operatorname{Orth}(N)$.
\end{lemma}

\begin{proof}
By \eqref{eq:nu-pushforward}, $u:=\sqrt N\,(I_D\ 0)O_N$ has law $\nu_N$. The right invariance of Haar probability on $\operatorname{Orth}(N)$ gives $O_N o\stackrel{\mathrm D}=O_N$, which proves \eqref{eq:stiefel-right-invariance}.

Let $\mu$ be any probability measure on $\mathbb U_N$ invariant under right multiplication by $\operatorname{Orth}(N)$, and let $\mathfrak h_N$ be Haar probability on $\operatorname{Orth}(N)$. For every $u\in\mathbb U_N$, there is an $o_u\in\operatorname{Orth}(N)$ such that $u=\sqrt N\,(I_D\ 0)o_u$. For any bounded measurable $\varphi$, right invariance of $\mu$ and left invariance of $\mathfrak h_N$ give
\begin{align*}
\int_{\mathbb U_N}\varphi(u)\,\mu(\d u)&=\iint_{\operatorname{Orth}(N)\times\mathbb U_N}\varphi(uo)\,\mathfrak h_N(\d o)\,\mu(\d u)\\
&=\int_{\operatorname{Orth}(N)}\varphi\big(\sqrt N\,(I_D\ 0)o\big)\,\mathfrak h_N(\d o)=\int_{\mathbb U_N}\varphi(u)\,\nu_N(\d u).
\end{align*}
Thus $\mu=\nu_N$.
\end{proof}

The parametrization of fixed self-overlap in \eqref{eq:stiefel-intro} gives $\sigma=au$ and $u=a^{-1}\sigma$ whenever $\sigma\sigma^{\intercal}=Na^2$. The overlap identity $N^{-1}(au)(au')^{\intercal}=a(N^{-1}uu'^{\intercal})a$ yields \eqref{eq:fixed-a-covariance-intro}, and $h\cdot\sum_{i=1}^N(au)_i=ah\cdot\sum_{i=1}^Nu_i$ because $a$ is symmetric. 

\subsection{Cascades and enrichment}
\label{s:enriched-properties}

We use a single continuous Ruelle probability cascade $\mathfrak R$ supported on the unit sphere of a separable Hilbert space, as constructed in \cite[Section~4]{chen2025free}. The scalar product of $\alpha,\alpha'\in\supp\mathfrak R$ is denoted by $\alpha\wedge\alpha'$, and the overlap of two conditionally independent samples from $\mathfrak R$ has the uniform distribution on $[0,1]$ after averaging over $\mathfrak R$.
For $q\in\mcl Q_\infty$, we put $q(1):=\lim_{s\uparrow1}q(s)$. Conditionally on $\mathfrak R$, let $z_i^q:\supp\mathfrak R\to\R^D$, $i\in\{1,\ldots,N\}$, be independent centered Gaussian fields, independent of $H_N$, with
\begin{equation}
\label{eq:enriched-cascade-covariance}
\E\big[z_i^q(\alpha)z_j^q(\alpha')^{\intercal}\mid\mathfrak R\big]=\one_{\{i=j\}}q(\alpha\wedge\alpha').
\end{equation}
The jointly measurable construction is given in \cite[Proposition~4.1]{chen2025free}. Every expectation involving the cascade and these fields integrates the Gaussian randomness first. We now express the constrained free energy \eqref{e.F_N,a_intro} in terms of its Hamiltonian and partition function, and define the corresponding random free energy. For $a\in\mcl A^\circ$ and $t\ge0$, define $H_{N,a}^{t,q}:\mathbb U_N\times\supp\mathfrak R\to\R$, its partition function, its random free energy, and its averaged free energy by 
\begin{align}
H_{N,a}^{t,q}(u,\alpha)&:=\sqrt{2t}H_N(au)-Nt\xi_a(I_D)+\sqrt2\sum_{i=1}^Nz_i^q(\alpha)\cdot u_i-Nq(1)\cdot I_D+ah\cdot\sum_{i=1}^Nu_i,\label{eq:enriched-fixed-a-H}\\
Z_{N,a}(t,q)&:=\iint_{\mathbb U_N\times\supp\mathfrak R}\exp\big(H_{N,a}^{t,q}(u,\alpha)\big)\nu_N(\mathrm d u)\mathfrak R(\mathrm d\alpha),\label{eq:fixed-a-Z}\\
F_{N,a}(t,q)&:=-\frac1N\log Z_{N,a}(t,q),\qquad\text{and}\qquad
\bar F_{N,a}(t,q)\stackrel{\eqref{e.F_N,a_intro}}{=}\E F_{N,a}(t,q).\label{eq:enriched-fixed-a}
\end{align}
The bracket $\langle\cdot\rangle_{N,a,t,q}$ is expectation under the probability measure obtained by normalizing the integrand in \eqref{eq:fixed-a-Z}, with independent replicas when the observable has several arguments. Thus $\langle\cdot\rangle_{N,a,t,q}$ is conditional on all the disorder, while $\E\langle\cdot\rangle_{N,a,t,q}$ also averages over that disorder.

For a finite-step path, we write
\begin{equation}
\label{eq:finite-step-q}
q=\sum_{j=0}^kq_j\one_{[m_j,m_{j+1})},\qquad\text{with}\quad 0=m_0<m_1<\cdots<m_k<m_{k+1}=1\quad\text{and}\quad 0\preccurlyeq q_0\preccurlyeq\cdots\preccurlyeq q_k.
\end{equation}
For the path in \eqref{eq:finite-step-q}, we also use the $k$-level discrete Ruelle probability cascade with parameters $(m_j)_{j=1}^k$ and normalized leaf weights $(v_\alpha)_{\alpha\in\N^k}$, as in \cite[Section~4]{chen2025free}. With $\N^0:=\{\varnothing\}$, its vertices form the rooted tree
\begin{equation*}
\N^0\cup\N^1\cup\cdots\cup\N^k.
\end{equation*}
Each vertex in $\N^j$ with $j\in\{0,\ldots,k-1\}$ is joined to its extensions in $\N^{j+1}$. For $\alpha=(\alpha_1,\ldots,\alpha_k)\in\N^k$, set $\alpha|_0:=\varnothing$ and $\alpha|_j:=(\alpha_1,\ldots,\alpha_j)\in\N^j$ for $j\in\{1,\ldots,k\}$. We write $\alpha\wedge\alpha':=\max\{j\in\{0,\ldots,k\}:\alpha|_j=\alpha'|_j\}$. For independent standard Gaussian vectors $\eta_{i,\gamma}\in\R^D$, indexed by $i\in\{1,\ldots,N\}$ and tree vertices $\gamma$, independently of the weights, set
\begin{equation*}
z_i^q(\alpha):=q_0^{1/2}\eta_{i,\varnothing}+\sum_{j=1}^k(q_j-q_{j-1})^{1/2}\eta_{i,\alpha|_j}\quad\text{with covariance}\quad \E z_i^q(\alpha)z_i^q(\alpha')^{\intercal}=q_{\alpha\wedge\alpha'}.
\end{equation*}
By \cite[Proposition~4.6]{chen2025free}, the continuous and discrete constructions have the same averaged logarithmic partition function for the finite-step models considered below. We use the discrete representation to compute that function. For the quadratic Gaussian integrals below, the recursive identity extends by truncation and monotone convergence when the Gaussian integrals are finite, as verified in \eqref{eq:matrix-gaussian-recursion}.

\begin{lemma}[Lipschitz continuity of $\bar F_{N,a}(t,\cdot)$]
\label{lem:fixed-path-Lipschitz}
For every $N\in\N$ with $N\ge D$, $a\in\mcl A^\circ$, $t\ge0$, and bounded $q,q'\in\mcl Q$, we have
\begin{equation}
\label{eq:fixed-enriched-q-Lipschitz}
|\bar F_{N,a}(t,q)-\bar F_{N,a}(t,q')|\le\sqrt D\int_0^1|q(s)-q'(s)|\,\mathrm d s.
\end{equation}
Consequently, $\bar F_{N,a}(t,\cdot)$ has a unique $L^1$-continuous extension to $\mcl Q_1$.
\end{lemma}

\begin{proof}
The argument of \cite[Proposition~5.1, (5.7)]{chen2025free} applies to the Hamiltonian \eqref{eq:enriched-fixed-a-H} and covariance \eqref{eq:enriched-cascade-covariance}, with reference measure $\nu_N$. The covariance terms evaluated at a replica paired with itself cancel against the derivative of $-Nq(1)\cdot I_D$, since $uu^{\intercal}/N=I_D$. Moreover, \cite[Proposition~4.8]{chen2025free} identifies the averaged two-replica cascade overlap. Together with $|u(u')^{\intercal}/N|\le\sqrt D$, these facts give \eqref{eq:fixed-enriched-q-Lipschitz} for finite-step paths. Approximation proves the extension.
\end{proof}

\subsection{Concentration over Gaussian fields and cascade}
\label{s:can-concentration}

Let $\mu$ be a deterministic probability measure on a measurable space $\Omega$, let $G_0$ be a centered Gaussian field on $\Omega$ independent of $\mathfrak R$, and let $A:\Omega\to\R$ be a bounded measurable function. Let $J$ be finite. For every $j\in J$, let $q_j\in\mcl Q_\infty$ and let $z_j(\alpha)\in\R^D$ be conditionally independent centered Gaussian cascade fields, given $\mathfrak R$, with covariance $q_j(\alpha\wedge\alpha')$, independent of $G_0$. For bounded measurable maps $v_j:\Omega\to\R^D$, define
\begin{equation*}
G(x,\alpha):=G_0(x)+A(x)+\sum_{j\in J}z_j(\alpha)\cdot v_j(x)\quad\text{and}\quad Z:=\iint_{\Omega\times\supp\mathfrak R}e^{G(x,\alpha)}\mu(\mathrm d x)\mathfrak R(\mathrm d\alpha).
\end{equation*}
Assume joint measurability of these fields and set
\begin{equation}
\label{eq:can-concentration-variance}
V:=\sup_{x\in\Omega}\left\{\E G_0(x)^2+\sum_{j\in J}v_j(x)^{\intercal}q_j(1)v_j(x)\right\}<\infty.
\end{equation}

\begin{lemma}[Concentration of free energy]
\label{lem:total-concentration-continuous}
There exists a universal constant $C$ such that
\begin{equation}
\label{eq:can-total-concentration}
\E|\log Z-\E\log Z|\le C\sqrt{1+V}.
\end{equation}
\end{lemma}

\begin{proof}
We follow \cite[proof of Lemma~B.1]{chenMourrat2026spherical}, specifying the changes for matrix-valued paths. Put $M:=\E[\log Z\mid\mathfrak R]$. Conditional on $\mathfrak R$, the centered Gaussian field $G-A$ has variance at most $V$ by \eqref{eq:can-concentration-variance}. Thus \cite[proof of Lemma~B.1, (B.6)]{chenMourrat2026spherical} gives
\begin{equation}
\label{eq:can-conditional-variance}
\E\big[(\log Z-M)^2\mid\mathfrak R\big]\le CV.
\end{equation}
For $\delta\in(0,1/2]$, set $\widetilde q_j(s):=q_j(s\vee\delta)$. Let $\widetilde Z$ be obtained from $Z$ by replacing $z_j$ by Gaussian fields with covariance paths $\widetilde q_j$, independent of the original fields conditional on $\mathfrak R$, and retaining $G_0$ and $A$. Put $\widetilde M:=\E[\log\widetilde Z\mid\mathfrak R]$. Monotonicity gives $0\preccurlyeq(\widetilde q_j-q_j)(s)\preccurlyeq q_j(1)\one_{[0,\delta)}(s)$. Matrix Cauchy--Schwarz and \eqref{eq:can-concentration-variance} therefore give
\begin{equation}
\label{eq:can-matrix-covariance-bound}
\left|\sum_{j\in J}v_j(x)^{\intercal}(\widetilde q_j-q_j)(\alpha\wedge\alpha')v_j(x')\right|\le V\one_{\{\alpha\wedge\alpha'<\delta\}}.
\end{equation}
For every $(x,\alpha)\in\Omega\times\supp\mathfrak R$, the original and modified scalar Gaussian fields have the same variance, namely $\E G_0(x)^2+\sum_{j\in J}v_j(x)^{\intercal}q_j(1)v_j(x)$, because $\alpha\wedge\alpha=1$ and $\widetilde q_j(1)=q_j(1)$. Thus their equal-index covariance terms cancel in Gaussian interpolation. Interpolate between $q_j$ and $\widetilde q_j$, conditional on $\mathfrak R$. At each interpolation parameter, \cite[Proposition~4.8]{chen2025free} gives $\E\langle\one_{\{\alpha\wedge\alpha'<\delta\}}\rangle=\delta$. Applying \cite[proof of Lemma~B.1, (B.7)--(B.8)]{chenMourrat2026spherical} with \eqref{eq:can-matrix-covariance-bound} yields the comparison
\begin{equation}
\label{eq:can-conditional-mean-comparison}
\E|M-\widetilde M|\le\frac12V\delta.
\end{equation}

To bound $\operatorname{Var}(\log\widetilde Z)$, approximate the paths $\widetilde q_j$ from below by finite-step paths $\widehat q_j$ on a common partition $0=m_0<m_1=\delta<\cdots<m_{K+1}=1$. Write $\widehat q_{j,\ell}$ for their values on $[m_\ell,m_{\ell+1})$ and set $\widehat q_{j,-1}:=0$. The corresponding fields $z_j^{\widehat q_j}:\N^K\to\R^D$ are given by
\begin{equation*}
z_j^{\widehat q_j}(\alpha):=\sum_{\ell=0}^K(\widehat q_{j,\ell}-\widehat q_{j,\ell-1})^{1/2}\eta_{j,\alpha|_\ell},
\end{equation*}
where the $\eta_{j,\gamma}$ are independent standard Gaussian vectors in $\R^D$. The matrices under the square roots are positive semidefinite, and for every $x\in\Omega$ their coefficients satisfy
\begin{equation*}
\sum_{j\in J}\sum_{\ell=0}^K\big|(\widehat q_{j,\ell}-\widehat q_{j,\ell-1})^{1/2}v_j(x)\big|^2=\sum_{j\in J}v_j(x)^{\intercal}\widehat q_j(1)v_j(x)\le V.
\end{equation*}
This is the squared-norm bound used in \cite[proof of Lemma~B.1, Step~4]{chenMourrat2026spherical}. The contribution of $G_0$ is included through \eqref{eq:can-concentration-variance}. For finite-step paths $\widetilde q_j$, \cite[Corollary~A.2, (A.9), and Lemma~A.3]{chenMourrat2026spherical} therefore give
\begin{equation}
\label{eq:can-modified-variance}
\operatorname{Var}(\log\widetilde Z)\le C(V+\delta^{-2}).
\end{equation}
The constant in \eqref{eq:can-modified-variance} is independent of the number of steps. The approximation argument in \cite[proof of Lemma~B.1, Step~3, (B.10)]{chenMourrat2026spherical} therefore proves \eqref{eq:can-modified-variance} for arbitrary $\widetilde q_j\in\mcl Q_\infty$.

Finally, the triangle inequality gives
\begin{equation*}
\E|\log Z-\E\log Z|\le\E|\log Z-M|+2\E|M-\widetilde M|+\E|\widetilde M-\E\widetilde M|.
\end{equation*}
By \eqref{eq:can-conditional-variance}, \eqref{eq:can-conditional-mean-comparison}, \eqref{eq:can-modified-variance}, and conditional Jensen's inequality, we obtain
\begin{equation*}
\E|\log Z-\E\log Z|\le C(\sqrt V+V\delta+\delta^{-1}).
\end{equation*}
Choosing $\delta:=\min\{1/2,(1+V)^{-1/2}\}$ proves \eqref{eq:can-total-concentration}.
\end{proof}

\subsection{The spherical initial condition on \texorpdfstring{$\mathbb U_N$}{\unichar{"1D54C}\_N}}
\label{s:spherical-initial}
\label{s:initial-condition}

Fix $v\in\R^D$ and a finite-step path $q$ as in \eqref{eq:finite-step-q}. The initial free energy $\psi_{N,v}(q)$ is defined by
\begin{equation}
\label{eq:psi-N-stiefel}
\psi_{N,v}(q):=-\frac1N\E\log\sum_\alpha v_\alpha\int_{\mathbb U_N} \exp\Ll(\sqrt2\sum_{i\le N}z_i^q(\alpha)\cdot u_i+\sum_{i\le N}v\cdot u_i-Nq_k\cdot I_D\Rr)\nu_N(\d u).
\end{equation}
With $v=ah$, \eqref{eq:enriched-fixed-a} at $t=0$ satisfies
\begin{equation}
\label{eq:enriched-fixed-a-initial}
\bar F_{N,a}(0,q)=\psi_{N,ah}(q).
\end{equation}
For $b\in\S^D$, define $\ell_1,\ldots,\ell_{k+1}\in\S^D$ by
\begin{equation}
\label{eq:L-recursion}
\ell_{k+1}:=b+2q_k\quad\text{and}\quad \ell_i:=\ell_{i+1}-2m_i(q_i-q_{i-1}), \quad i\in\{1,\ldots,k\}.
\end{equation}
To display the dependence on $b$, we write $\ell_i=\ell_i(b)$ for $i\in\{1,\ldots,k+1\}$. When $k=0$, we use only $\ell_1(b):=b+2q_0$. For fixed $q$, define the open convex set
\begin{equation}
\label{eq:initial-B-domain}
\mathcal B_q:=\{b\in\S^D:\ell_1(b)\succ0\}.
\end{equation}
All $\ell_i(b)$ are positive definite for $b\in\mathcal B_q$. Define $G_v:\mathcal B_q\to\R$ by
\begin{equation}
\label{eq:G-explicit}
G_v(b):=\frac12\left[\tr b-D-\log\det \ell_{k+1}(b)+2\tr(\ell_1(b)^{-1}q_0)+v^{\intercal}\ell_1(b)^{-1}v+\sum_{i=1}^k\frac1{m_i}\log\frac{\det \ell_{i+1}(b)}{\det \ell_i(b)}\right].
\end{equation}
The dependence of $G_v$ on $q$ is suppressed. Proposition~\ref{prop:stiefel-initial} identifies $\lim_{N\to\infty}\psi_{N,v}(q)$ with $-\inf_{b\in\mathcal B_q}G_v(b)$. 

\begin{lemma}[Coercivity of $G_v$]
\label{lem:initial-coercivity}
For every finite-step $q$ and $v\in\R^D$, the function $G_v:\mathcal B_q\to\R$ satisfies
\begin{equation}
\label{eq:G-coercive}
\lim_{n\to\infty}G_v(b_n)=+\infty
\end{equation}
for every sequence $(b_n)$ leaving every compact subset of $\mathcal B_q$. Consequently, $G_v$ attains its minimum on $\mathcal B_q$. Moreover, for fixed step times $0=m_0<\cdots<m_k<m_{k+1}=1$, every compact set $K$ of admissible data $(q_0,\ldots,q_k,v)$, and every $M\in\R$, there exists $R\in[1,\infty)$ such that
\begin{equation*}
\{b\in\mathcal B_q:G_v(b)\le M\}\subseteq\{b\in\mathcal B_q:|b|\le R,\ \ell_1(b)\succcurlyeq R^{-1}I_D\},\qquad (q_0,\ldots,q_k,v)\in K.
\end{equation*}
\end{lemma}
\begin{proof}
If \eqref{eq:G-coercive} failed, we could pass to a subsequence with $G_v(b_n)$ bounded above. Since $\ell_1(b)$ is an affine translate of $b$, a further subsequence would satisfy either $\lim_{n\to\infty}|b_n|=+\infty$, or $(b_n)$ would be bounded with $\lim_{n\to\infty}\lambda_{\min}(\ell_1(b_n))=0$. We rule out both cases.

First suppose that $\lim_{n\to\infty}|b_n|=+\infty$. The condition $\ell_1(b_n)\succ0$ bounds the negative eigenvalues of $b_n$ from below by a constant depending only on $q$. Thus $\tr b_n\ge |b_n|-C_{D,q}$. Since $\ell_{i+1}(b_n)\succcurlyeq\ell_i(b_n)$, every ratio $\det\ell_{i+1}(b_n)/\det\ell_i(b_n)$ is at least one. Consequently, the logarithmic part of \eqref{eq:G-explicit} satisfies
\begin{equation*}
-\log\det\ell_{k+1}(b_n) +\sum_{i=1}^k\frac1{m_i}\log\frac{\det\ell_{i+1}(b_n)}{\det\ell_i(b_n)} \ge-C_{D,q}\log\Ll(2+|b_n|\Rr).
\end{equation*}
The two inverse-matrix terms in \eqref{eq:G-explicit} are nonnegative. Hence \eqref{eq:G-coercive} follows from the linear lower bound on $\tr b_n$.

Next, suppose that $(b_n)$ is bounded and $\lim_{n\to\infty}\lambda_{\min}(\ell_1(b_n))=0$. Every determinant ratio $\det\ell_{i+1}/\det\ell_i$ is at least one, and $m_i<1$, so we have
\begin{equation*}
-\log\det\ell_{k+1}+\sum_{i=1}^k\frac1{m_i}\log\frac{\det\ell_{i+1}}{\det\ell_i}\ge-\log\det\ell_{k+1}+\sum_{i=1}^k\log\frac{\det\ell_{i+1}}{\det\ell_i}=-\log\det\ell_1.
\end{equation*}
This bound also holds with an empty sum when $k=0$. Since the other eigenvalues of $\ell_1(b_n)$ are bounded above, we have $\lim_{n\to\infty}\det\ell_1(b_n)=0$. The trace term in \eqref{eq:G-explicit} remains bounded, and its inverse-matrix terms are nonnegative, so \eqref{eq:G-coercive} follows. The estimates above depend on $(q,v)$ only through bounds for $|q_j|$ and $|v|$, with fixed step times. They therefore apply uniformly to any compact family $K$ of admissible $(q_0,\ldots,q_k,v)$. Every sublevel set is bounded and has $\lambda_{\min}(\ell_1(b))$ bounded away from zero, uniformly on $K$. This proves both assertions of the lemma.
\end{proof}

For $N\in\N$ with $N\ge D$, let $\gamma_N$ be the standard Gaussian probability measure on $\R^{D\times N}$ and define $\widehat a_N:\R^{D\times N}\to\S_+^D$ by $\widehat a_N(\omega):=N^{-1}\omega\omega^{\intercal}$. Write $\omega_i\in\R^D$ for the $i$th column of $\omega$. For $b\in\mathcal B_q$, define
\begin{align}
\label{eq:canonical-Z}
    Z_{N,b,v}(q):=\sum_\alpha v_\alpha\int_{\R^{D\times N}}\exp\bigg(\sqrt2\sum_{i\le N}z_i^q(\alpha)\cdot\omega_i&+\sum_{i\le N}v\cdot\omega_i-Nq_k\cdot\widehat a_N(\omega)\notag\\
    &-\frac N2(b-I_D)\cdot\Ll(\widehat a_N(\omega)-I_D\Rr)\bigg)\gamma_N(\d\omega).
\end{align}
We write $\langle\cdot\rangle_b$ for the Gibbs expectation obtained by normalizing the integrand in \eqref{eq:canonical-Z}, including the cascade weights $v_\alpha$, and $\E\langle\cdot\rangle_b$ for the corresponding annealed expectation.

\begin{lemma}
\label{lem:canonical-evaluation}
For every $b\in\mathcal B_q$ and $N\ge D$, we have
\begin{equation}
\label{eq:canonical-G-identity}
\frac1N\E\log Z_{N,b,v}(q)=G_v(b).
\end{equation}
\end{lemma}
\begin{proof}
Fix $b\in\mathcal B_q$ and write $\ell_i=\ell_i(b)$. In this proof, put $c_0:=q_0$ and $c_j:=q_j-q_{j-1}$ for $j\in\{1,\ldots,k\}$, and let $g_0,\ldots,g_k$ be independent centered Gaussian vectors with covariances $2c_0,\ldots,2c_k$, respectively. We first integrate the columns in \eqref{eq:canonical-Z}. The standard Gaussian density and the correction terms give
\begin{align}
\label{eq:terminal-precision}
    -\tfrac12\sum_{i\le N}|\omega_i|^2-Nq_k\cdot\widehat a_N(\omega)-\tfrac N2(b-I_D)\cdot\Ll(\widehat a_N(\omega)-I_D\Rr)=-\tfrac12\sum_{i\le N}\omega_i^{\intercal}\ell_{k+1}\omega_i+\tfrac N2(\tr b-D).
\end{align}
For $y\in\R^D$, Gaussian integration gives
\begin{equation*}
\frac1{(2\pi)^{D/2}}\int_{\R^D}e^{y\cdot\omega-\omega^{\intercal}\ell_{k+1}\omega/2}\,\mathrm d\omega=\frac{e^{y^{\intercal}\ell_{k+1}^{-1}y/2}}{(\det\ell_{k+1})^{1/2}}.
\end{equation*}
By \eqref{eq:terminal-precision}, the Gaussian columns can be integrated independently for fixed cascade weights and increments, giving
\begin{equation}
\label{eq:canonical-product-Z}
Z_{N,b,v}(q)=e^{N(\tr b-D-\log\det\ell_{k+1})/2}\sum_\alpha v_\alpha\exp\Ll(\frac12\sum_{i=1}^N(v+\sqrt2z_i^q(\alpha))^{\intercal}\ell_{k+1}^{-1}(v+\sqrt2z_i^q(\alpha))\Rr).
\end{equation}

We next evaluate the expectation of the logarithm of the cascade sum in \eqref{eq:canonical-product-Z}. 
Define $X_i:\R^D\to\R$, $i\in\{1,\ldots,k+1\}$, by
\begin{equation}
\label{eq:cascade-terminal}
X_{k+1}(y):=\frac12y^{\intercal}\ell_{k+1}^{-1}y
\end{equation}
and, for $i\in\{k,\ldots,1\}$, by the backward recursion
\begin{equation}
\label{eq:cascade-recursion}
X_i(y):=\frac1{m_i}\log\E_{g_i}\exp\Ll(m_iX_{i+1}(y+g_i)\Rr).
\end{equation}
The increments of $\sqrt2z_i^q(\alpha)$ along the ancestral line have the laws of $g_0,\ldots,g_k$, independently across the column index. Applying \cite[Corollary~A.2, (A.7)]{chenMourrat2026spherical} to \eqref{eq:cascade-terminal}--\eqref{eq:cascade-recursion} gives
\begin{equation}
\label{eq:canonical-cascade-expectation}
\frac1N\E\log\sum_\alpha v_\alpha\exp\Ll(\sum_{i=1}^NX_{k+1}(v+\sqrt2z_i^q(\alpha))\Rr)=\E_{g_0}X_1(v+g_0).
\end{equation}

To calculate the recursion, let $c\succcurlyeq0$, $\ell\succ0$, and $m>0$ satisfy $\ell-2mc\succ0$, and let $g$ be centered Gaussian with covariance $2c$. Completing the square gives
\begin{equation}
\label{eq:matrix-gaussian-recursion}
\frac1m\log\E_g\exp\Ll(\frac m2(y+g)^{\intercal}\ell^{-1}(y+g)\Rr) =\frac1{2m}\log\frac{\det\ell}{\det(\ell-2mc)} +\frac12y^{\intercal}(\ell-2mc)^{-1}y.
\end{equation}
Indeed, with $g=\sqrt2c^{1/2}\eta$ for a standard Gaussian vector $\eta\in\R^D$, the determinant and quadratic coefficient are identified by
\begin{equation*}
\det(I_D-2mc^{1/2}\ell^{-1}c^{1/2})^{-1/2} =\left(\frac{\det\ell}{\det(\ell-2mc)}\right)^{1/2},
\end{equation*}
\begin{equation*}
\ell^{-1}+2m\ell^{-1}c^{1/2}(I_D-2mc^{1/2}\ell^{-1}c^{1/2})^{-1}c^{1/2}\ell^{-1}=(\ell-2mc)^{-1}.
\end{equation*}
The second identity follows by multiplication by $\ell-2mc$, without assuming that $c$ and $\ell$ commute. For the root increment, we also have
\begin{equation}
\label{eq:matrix-gaussian-m0}
\E_g\frac12(y+g)^{\intercal}\ell^{-1}(y+g) =\frac12y^{\intercal}\ell^{-1}y+\tr(\ell^{-1}c).
\end{equation}

Applying \eqref{eq:matrix-gaussian-recursion} to \eqref{eq:cascade-recursion} and using \eqref{eq:L-recursion}, we obtain
\begin{equation}
\label{eq:cascade-quadratic-form}
X_i(y)=A_i+\frac12y^{\intercal}\ell_i^{-1}y, \qquad i\in\{1,\ldots,k+1\},
\end{equation}
where
\begin{equation}
\label{eq:cascade-constants}
A_{k+1}:=0\quad\text{and}\quad A_i:=A_{i+1}+\frac1{2m_i}\log\frac{\det\ell_{i+1}}{\det\ell_i}.
\end{equation}
For the root increment, equations \eqref{eq:cascade-quadratic-form} and \eqref{eq:matrix-gaussian-m0} give
\begin{equation}
\label{eq:cascade-root}
\E_{g_0}X_1(v+g_0) =A_1+\tr(\ell_1^{-1}q_0)+\frac12v^{\intercal}\ell_1^{-1}v.
\end{equation}
By \eqref{eq:canonical-product-Z}, \eqref{eq:canonical-cascade-expectation}, \eqref{eq:cascade-constants}, and \eqref{eq:cascade-root}, we obtain \eqref{eq:canonical-G-identity}.
\end{proof}

For the next lemma and its proof, put $c_0:=q_0$ and $c_j:=q_j-q_{j-1}$ for $j\in\{1,\ldots,k\}$. In this lemma, $\mcl N(m,\Sigma)$ denotes the Gaussian law on $\R^D$ with mean $m\in\R^D$ and covariance $\Sigma\in\S_+^D$, including singular $\Sigma$.

\begin{lemma}
\label{lem:canonical-column-law}
For $b\in\mathcal B_q$ and $N\ge D$, the columns $(\omega_i)_{i=1}^N$ under $\E\langle\cdot\rangle_b$ are independent and identically distributed, with a common law independent of $N$. More precisely, for each $i$, let $g_{j,i}:=\sqrt2c_j^{1/2}\eta_{i,\alpha|_j}$ for $j\in\{0,\ldots,k\}$, where $\alpha$ is sampled from the Gibbs measure in \eqref{eq:canonical-Z}, and set $y_{j,i}:=v+\sum_{r=0}^jg_{r,i}$. Then $g_{0,i}\sim\mcl N(0,2q_0)$, and, for $j\in\{1,\ldots,k\}$, we have
\begin{equation}
\label{eq:canonical-increment-Gaussian}
\mathrm{Law}(g_{j,i}\mid g_{0,i},\ldots,g_{j-1,i})=\mcl N\big(2m_jc_j\ell_j^{-1}y_{j-1,i},\,2c_j+4m_jc_j\ell_j^{-1}c_j\big).
\end{equation}
Conditional on $\alpha$ and $(g_{j,i})_{j=0}^k$, the column $\omega_i$ has distribution
\begin{equation}
\label{eq:canonical-column-conditional}
\mathrm{Law}(\omega_i\mid\alpha,(g_{j,i})_{j=0}^k)=\mcl N(\ell_{k+1}^{-1}y_{k,i},\ell_{k+1}^{-1}).
\end{equation}
For every compact set $K\subset\mathcal B_q$, there is a constant $C_{q,v,K}<\infty$ such that
\begin{equation}
\label{eq:canonical-column-fourth-moment}
\sup_{\substack{N\ge D,\ b\in K\\i\in\{1,\ldots,N\}}}\E\langle|\omega_i|^4\rangle_b\le C_{q,v,K}.
\end{equation}
\end{lemma}
\begin{proof}
For fixed $\alpha$ and $(\eta_{i,\alpha|_j})_{i\le N,\,0\le j\le k}$, completing the square in \eqref{eq:canonical-Z} gives \eqref{eq:canonical-column-conditional}. The sampled-path tilt in \cite[Corollary~A.2, (A.8), and Lemma~A.4(2)]{chenMourrat2026spherical} applies to the sum of the $N$ column functionals, so the collections $(g_{j,i})_{j=0}^k$ are independent across $i$. Integrating the increments at levels $j+1,\ldots,k$ backwards and using $\ell_{r+1}-\ell_r=2m_r c_r$ for $r\in\{1,\ldots,k\}$ gives the following conditional law for $j\in\{1,\ldots,k\}$:
\begin{equation}
\label{eq:canonical-one-column-law}
\mathrm{Law}(g_{j,i}\mid g_{0,i},\ldots,g_{j-1,i})(\mathrm d g)\propto\exp\Big(\frac{m_j}2(y_{j-1,i}+g)^{\intercal}\ell_{j+1}^{-1}(y_{j-1,i}+g)\Big)\,\mcl N(0,2c_j)(\mathrm d g).
\end{equation}
Writing $g=\sqrt2c_j^{1/2}z$ with $z\sim\mcl N(0,I_D)$, its quadratic precision is
\begin{equation}
\label{eq:canonical-tilted-precision}
I_D-2m_jc_j^{1/2}\ell_{j+1}^{-1}c_j^{1/2}=\big(I_D+2m_jc_j^{1/2}\ell_j^{-1}c_j^{1/2}\big)^{-1}\succ0.
\end{equation}
Completing the square in \eqref{eq:canonical-one-column-law} using \eqref{eq:canonical-tilted-precision} gives \eqref{eq:canonical-increment-Gaussian}, including singular $c_j$, since no inverse of $c_j$ is required. For $j=0$, no cascade tilt acts on the root increment $g_{0,i}$, which has law $\mcl N(0,2q_0)$. By \eqref{eq:canonical-increment-Gaussian}, $y_{j,i}$ is obtained from $y_{j-1,i}$ by an affine Gaussian transition with deterministic covariance. Equation \eqref{eq:canonical-column-conditional} therefore gives a Gaussian marginal law for $\omega_i$, independent of $N$, whose fourth moment is locally bounded for $b\in\mathcal B_q$.
\end{proof}

By Lemma~\ref{lem:initial-coercivity}, choose a minimizer $b_*\in\mathcal B_q$ of $G_v$. 

\begin{lemma}
\label{lem:canonical-localization}
There exists $C_{q,v,b_*}<\infty$ such that, for every $N\ge D$, we have
\begin{equation}
\label{eq:canonical-localization}
\E\left\langle\left|\frac1N\sum_{i=1}^N\omega_i\omega_i^{\intercal}-I_D\right|^2\right\rangle_{b_*}\le\frac{C_{q,v,b_*}}N.
\end{equation}
\end{lemma}
\begin{proof}
Recall from \eqref{eq:canonical-Z} that $\widehat a_N=N^{-1}\sum_{i=1}^N\omega_i\omega_i^{\intercal}$. The independence and moment bound needed below are given by Lemma~\ref{lem:canonical-column-law}, in particular \eqref{eq:canonical-column-fourth-moment}. Independence here holds under $\E\langle\cdot\rangle_b$, not under the Gibbs measure conditional on the disorder.

For $b\in\mathcal B_q$ and $c\in\S^D$, differentiating \eqref{eq:canonical-G-identity} in direction $c$ gives
\begin{equation}
c\cdot\nabla G_v(b)=\frac12c\cdot\big(I_D-\E\langle\widehat a_N\rangle_b\big).
\label{eq:canonical-multiplier-derivative}
\end{equation}
Since $b_*$ minimizes $G_v$ in the open set $\mathcal B_q$, we have $c\cdot\nabla G_v(b_*)=0$ for every $c\in\S^D$. Equation \eqref{eq:canonical-multiplier-derivative} gives $\E\langle\widehat a_N\rangle_{b_*}=I_D$. Since the vectors $\omega_1,\ldots,\omega_N\in\R^D$ are identically distributed under $\E\langle\cdot\rangle_{b_*}$, we obtain
\begin{equation}
\label{eq:canonical-moment-identity}
\E\langle\omega_1\omega_1^{\intercal}\rangle_{b_*}=\E\langle\widehat a_N\rangle_{b_*}=I_D.
\end{equation}
For distinct $i,j\in\{1,\ldots,N\}$, Lemma~\ref{lem:canonical-column-law} gives the independence of $\omega_i$ and $\omega_j$ under $\E\langle\cdot\rangle_{b_*}$. The matrices $\omega_i\omega_i^{\intercal}-I_D$ and $\omega_j\omega_j^{\intercal}-I_D$ are therefore independent and centered by \eqref{eq:canonical-moment-identity}, which gives
\begin{equation*}
\E\left\langle(\omega_i\omega_i^{\intercal}-I_D)\cdot(\omega_j\omega_j^{\intercal}-I_D)\right\rangle_{b_*}=0.
\end{equation*}
Expanding the square of $N^{-1}\sum_i(\omega_i\omega_i^{\intercal}-I_D)$ and using \eqref{eq:canonical-column-fourth-moment}, we obtain
\begin{equation*}
\E\left\langle|\widehat a_N-I_D|^2\right\rangle_{b_*}=\frac1N\E\left\langle|\omega_1\omega_1^{\intercal}-I_D|^2\right\rangle_{b_*}\le\frac{C_{q,v,b_*}}N,
\end{equation*}
which proves \eqref{eq:canonical-localization}.
\end{proof}

The following estimate controls the normalized logarithm of \eqref{eq:canonical-Z}. 

\begin{lemma}
\label{lem:canonical-log-moment}
For every $b\in\mathcal B_q$, there exists $C_{q,v,b}<\infty$ such that, for all $N\ge D$, we have
\begin{equation}
\label{eq:canonical-log-second-moment}
\E\left[\left(\frac1N\log Z_{N,b,v}(q)\right)^2\right]=G_v(b)^2+\frac1{N^2}\operatorname{Var}\big(\log Z_{N,b,v}(q)\big)\le C_{q,v,b}.
\end{equation}
\end{lemma}
\begin{proof}
Fix $b\in\mathcal B_q$, and put $g_{0,i}:=\sqrt2q_0^{1/2}\eta_{i,\varnothing}$. Define the quadratic forms $Q_i(y):=y^{\intercal}\ell_i^{-1}y/2$ for $i\in\{1,k+1\}$ and
\begin{equation*}
Y_N:=\sum_\alpha v_\alpha\exp\Ll(\sum_{i=1}^NQ_{k+1}(v+\sqrt2z_i^q(\alpha))\Rr).
\end{equation*}
Equation \eqref{eq:canonical-product-Z} implies $\operatorname{Var}(\log Z_{N,b,v}(q))=\operatorname{Var}(\log Y_N)$. When $k=0$, the cascade consists of its root, so
\begin{equation*}
\operatorname{Var}(\log Y_N)=N\operatorname{Var}(Q_1(v+g_{0,1})).
\end{equation*}
Suppose $k\ge1$. Let $(w_\alpha)$ be the unnormalized cascade weights, and put $T:=\sum_\alpha w_\alpha$ and
\begin{equation*}
U_N:=\sum_\alpha w_\alpha\exp\Ll(\sum_{i=1}^NQ_{k+1}(v+\sqrt2z_i^q(\alpha))\Rr),
\end{equation*}
so that $Y_N=U_N/T$. Iterating the tilted-cascade invariance in \cite[Lemma~A.1, (A.6)]{chenMourrat2026spherical} and the Gaussian recursion \eqref{eq:matrix-gaussian-recursion} shows that, conditional on $(g_{0,i})_{i=1}^N$, $U_N$ has the same law as
\begin{equation*}
\exp\Ll(NC_{q,b}+\sum_{i=1}^NQ_1(v+g_{0,i})\Rr)T',
\end{equation*}
where $C_{q,b}\in\R$ is deterministic and $T'$ has the law of $T$ and is independent of the root increments. Since $T$ is also independent of these increments, the variance inequality gives
\begin{equation*}
\operatorname{Var}(\log Y_N)\le2N\operatorname{Var}(Q_1(v+g_{0,1}))+4\operatorname{Var}(\log T).
\end{equation*}
The first variance is finite because $Q_1$ is quadratic and $g_{0,1}$ is Gaussian. The second is finite by \cite[proof of Corollary~A.2, (A.9)]{chenMourrat2026spherical}. Thus
\begin{equation}
\label{eq:G-random-concentration}
\operatorname{Var}\big(\log Z_{N,b,v}(q)\big)\le C_{q,v,b}N.
\end{equation}
Combining \eqref{eq:G-random-concentration} with \eqref{eq:canonical-G-identity} proves \eqref{eq:canonical-log-second-moment}.
\end{proof}

To compare the Gaussian integral over $\R^{D\times N}$ in \eqref{eq:canonical-Z} with the integral over $\mathbb U_N$ in \eqref{eq:psi-N-stiefel}, we use the following decomposition. 

\begin{lemma}[Gaussian polar decomposition]
\label{lem:gaussian-polar}
For every nonnegative measurable function $\Phi:\R^{D\times N}\to[0,+\infty]$, we have
\begin{equation}
\label{eq:gaussian-polar-integral}
\int_{\R^{D\times N}}\Phi(\omega)\gamma_N(\mathrm d\omega)=\iint_{\mathbb U_N\times\R^{D\times N}}\Phi\big((N^{-1}\omega\omega^{\intercal})^{1/2}u\big)\nu_N(\mathrm d u)\gamma_N(\mathrm d\omega).
\end{equation}
Moreover, for every $\delta\in(0,1)$, we have
\begin{equation}
\label{eq:Qhat-shell-mass}
\lim_{N\to\infty}\frac1N\log\gamma_N\Ll(\Big\{\omega\in\R^{D\times N}:\big|N^{-1}\omega\omega^{\intercal}-I_D\big|\le\delta\Big\}\Rr)=0.
\end{equation}
\end{lemma}

\begin{proof}
A standard Gaussian matrix $\omega\in\R^{D\times N}$ has full row rank almost surely. Put $b:=N^{-1}\omega\omega^{\intercal}$ and $u:=b^{-1/2}\omega$. Then $uu^{\intercal}=NI_D$ and $\omega=b^{1/2}u$. Right-orthogonal invariance of $\gamma_N$ and Lemma~\ref{lem:stiefel-invariance} imply that the conditional law of $u$ given $b$ is $\nu_N$ and does not depend on $b$. This proves \eqref{eq:gaussian-polar-integral}. Under a coupling by independent standard Gaussian columns $(\omega_i)_{i\in\N}$, the strong law of large numbers gives $\lim_{N\to\infty}\frac1N\sum_{i=1}^N\omega_i\omega_i^{\intercal}=I_D$ almost surely, which implies \eqref{eq:Qhat-shell-mass}.
\end{proof}

For $N\ge D$, let $\lambda_N$ be the law of $\widehat a_N(\omega)$ under $\gamma_N$. For $\delta\in(0,1/2)$, set
\begin{equation*}
A_\delta:=\{c\in\S_{++}^D:|c-I_D|\le\delta\}.
\end{equation*}
For $b\in\mathcal B_q$, define $Z_{N,b,v}^\delta(q)$ by restricting the $\gamma_N$ integral in \eqref{eq:canonical-Z} to $\{\omega:\widehat a_N(\omega)\in A_\delta\}$. 

For $c\in A_\delta$, define $W_c:\mathbb U_N\times\N^k\to\R$ by
\begin{align}
W_c(u,\alpha)&:=\sqrt2\sum_{i\le N}z_i^q(\alpha)\cdot(c^{1/2}u_i),\notag\\
\frac1N\E W_c(u,\alpha)W_{c'}(u',\alpha')&=2q_{\alpha\wedge\alpha'}\cdot\left(c^{1/2}\frac{u(u')^{\intercal}}N (c')^{1/2}\right).\label{eq:canonical-shell-cross-covariance}
\end{align}
Applying \eqref{eq:gaussian-polar-integral} to the restricted integral in \eqref{eq:canonical-Z}, conditional on the cascade and $z^q$, gives
\begin{align}
\label{eq:canonical-shell-polar}
    Z_{N,b,v}^\delta(q)=\sum_\alpha v_\alpha\iint_{\mathbb U_N\times A_\delta}\exp\bigg(W_c(u,\alpha)&+\sum_{i\le N}v\cdot(c^{1/2}u_i)\notag\\
    &-Nq_k\cdot c-\frac N2(b-I_D)\cdot(c-I_D)\bigg)\nu_N(\d u)\lambda_N(\d c).
\end{align}
Let $Z_{N,b,v}^{\delta,1}(q)$ denote the right-hand side of \eqref{eq:canonical-shell-polar} with $W_c(u,\alpha)$ replaced by $W_{I_D}(u,\alpha)$, retaining $\sum_{i\le N}v\cdot(c^{1/2}u_i)$ and the two correction terms.
Gaussian interpolation gives the following bound.

\begin{lemma}
\label{lem:canonical-shell-gaussian}
For all $b\in\mathcal B_q$, $\delta\in(0,1/2)$, and $N\ge D$, we have
\begin{equation}
\label{eq:shell-gaussian-log-comparison}
\left|\frac1N\E\log Z_{N,b,v}^\delta(q) -\frac1N\E\log Z_{N,b,v}^{\delta,1}(q)\right| \le C_{D,q}\delta.
\end{equation}
\end{lemma}
\begin{proof}

For $c\in A_\delta$, diagonalize $c$ orthogonally. The scalar bound $|\sqrt\lambda-1|\le|\lambda-1|$ for $\lambda>0$, applied to its eigenvalues, gives
\begin{equation}
\label{eq:Qhat-sqrt-shell}
|c^{1/2}-I_D|\le|c-I_D|\le\delta.
\end{equation}
Since $uu^{\intercal}=u'(u')^{\intercal}=NI_D$, we have $|u/\sqrt N|_{\mathrm{op}}=|u'/\sqrt N|_{\mathrm{op}}=1$, and thus $|u(u')^{\intercal}/N|_{\mathrm{op}}\le1$. Equations \eqref{eq:Qhat-sqrt-shell} and~\eqref{eq:canonical-shell-cross-covariance} imply
\begin{equation}
\label{eq:shell-covariance-comparison}
    \sup_{\substack{c,c'\in A_\delta,\,u,u'\in\mathbb U_N\\ \alpha,\alpha'}}
    \frac1N\left|\E W_c(u,\alpha)W_{c'}(u',\alpha')-\E W_{I_D}(u,\alpha)W_{I_D}(u',\alpha')\right|
    \le C_{D,q}\delta.
\end{equation}
Interpolate $\sqrt{1-\rho}\,W_c+\sqrt\rho\,W'_{I_D}$ for $\rho\in[0,1]$, where $W'_{I_D}$ is an independent copy of $W_{I_D}$. Gaussian integration by parts, applied to both the variances at each $(u,\alpha)$ and the covariances between two replicas, together with \eqref{eq:shell-covariance-comparison}, proves \eqref{eq:shell-gaussian-log-comparison}.
\end{proof}

Define $\widetilde Z_{N,v}^\delta(q)$ by
\begin{equation}
\label{eq:frozen-shell-partition}
\widetilde Z_{N,v}^\delta(q):=\lambda_N(A_\delta)\sum_\alpha v_\alpha\int_{\mathbb U_N}\exp\Ll(W_{I_D}(u,\alpha)+\sum_{i\le N}v\cdot u_i-Nq_k\cdot I_D\Rr)\nu_N(\d u).
\end{equation}
The next lemma removes the $c$-dependent terms from $Z_{N,b,v}^{\delta,1}(q)$.

\begin{lemma}
\label{lem:canonical-shell-deterministic}
For all $b\in\mathcal B_q$, $\delta\in(0,1/2)$, and $N\ge D$, we have
\begin{equation}
\label{eq:shell-deterministic-log-comparison}
\left|\frac1N\E\log Z_{N,b,v}^{\delta,1}(q) -\frac1N\E\log\widetilde Z_{N,v}^\delta(q)\right| \le C_{D,q,b,v}\delta.
\end{equation}
\end{lemma}
\begin{proof}
For $u\in\mathbb U_N$ and $c\in A_\delta$, the identities $|u|_{\mathrm{op}}=\sqrt N$ and $|\mathbf 1_N|=\sqrt N$ give
\begin{equation}
\label{eq:shell-source-comparison}
\frac1N\left|\sum_{i\le N}v\cdot(c^{1/2}-I_D)u_i\right| =\frac1N\left|v^{\intercal}(c^{1/2}-I_D)u\mathbf 1_N\right| \le C_D|v|\delta.
\end{equation}
Moreover, we have
\begin{equation}
\label{eq:shell-correction-comparison}
\left|q_k\cdot(c-I_D)+\frac12(b-I_D)\cdot(c-I_D)\right| \le\left(|q_k|+\frac12|b-I_D|\right)\delta.
\end{equation}

The partition function in \eqref{eq:frozen-shell-partition} is obtained from $Z_{N,b,v}^{\delta,1}(q)$ by replacing $c^{1/2}$ by $I_D$ in $\sum_{i\le N}v\cdot(c^{1/2}u_i)$ and replacing $c$ by $I_D$ in both correction terms. Equations \eqref{eq:shell-source-comparison} and \eqref{eq:shell-correction-comparison} bound the difference of the two exponents by $NC_{D,q,b,v}\delta$, uniformly in $(u,c,\alpha)$. Integrating the resulting pointwise bounds gives
\begin{equation*}
e^{-NC_{D,q,b,v}\delta}\widetilde Z_{N,v}^\delta(q)\le Z_{N,b,v}^{\delta,1}(q)\le e^{NC_{D,q,b,v}\delta}\widetilde Z_{N,v}^\delta(q).
\end{equation*}
Taking logarithms and expectations, then dividing by $N$, proves \eqref{eq:shell-deterministic-log-comparison}.
\end{proof}

From \eqref{eq:psi-N-stiefel} and \eqref{eq:frozen-shell-partition}, we have
\begin{equation*}
\frac1N\E\log\widetilde Z_{N,v}^\delta(q)=\frac1N\log\lambda_N(A_\delta)-\psi_{N,v}(q).
\end{equation*}
Lemmas~\ref{lem:canonical-shell-gaussian} and~\ref{lem:canonical-shell-deterministic} now give
\begin{equation}
\label{eq:shell-log-comparison}
\left|\frac1N\E\log Z_{N,b,v}^\delta(q)-\frac1N\log\lambda_N(A_\delta)+\psi_{N,v}(q)\right| \le C_{D,q,b,v}\delta.
\end{equation}

The preceding lemmas give the limit of \eqref{eq:psi-N-stiefel}.

\begin{proposition}
\label{prop:stiefel-initial}
For every path in \eqref{eq:finite-step-q}, we have
\begin{equation}
\label{eq:stiefel-determinant-formula}
\lim_{N\to\infty}\psi_{N,v}(q)=-\inf_{b\in\mathcal B_q}G_v(b).
\end{equation}
\end{proposition}
\begin{proof}
We prove the two inequalities in \eqref{eq:stiefel-determinant-formula}.

\smallskip
\noindent\emph{Step 1: The lower bound for $\psi_{N,v}(q)$.} Fix $b\in\mathcal B_q$ and $\delta\in(0,1/2)$. Since $Z_{N,b,v}(q)\ge Z_{N,b,v}^\delta(q)$, Lemma~\ref{lem:canonical-evaluation} and equations \eqref{eq:shell-log-comparison} and \eqref{eq:Qhat-shell-mass} give
\begin{equation}
\label{eq:initial-upper}
\limsup_{N\to\infty}(-\psi_{N,v}(q))\le G_v(b)+C_{D,q,b,v}\delta.
\end{equation}
Letting $\delta\downarrow0$ and taking the infimum over $b\in\mathcal B_q$ gives $\liminf_{N\to\infty}\psi_{N,v}(q)\ge-\inf_{b\in\mathcal B_q}G_v(b)$.

\smallskip
\noindent\emph{Step 2: The upper bound for $\psi_{N,v}(q)$.} Fix $\delta\in(0,1/2)$ and put $P_{N,\delta}:=\langle\one_{\{|\widehat a_N-I_D|\le\delta\}}\rangle_{b_*}$. Equation \eqref{eq:canonical-localization} and Markov's inequality give
\begin{equation}
\label{eq:PNdelta}
\E[1-P_{N,\delta}]\le\frac{C_{q,v,b_*}}{N\delta^2}\quad\text{and}\quad\P(P_{N,\delta}<1/2)\le\frac{2C_{q,v,b_*}}{N\delta^2}.
\end{equation}
Recall that $Z_{N,b_*,v}(q)$ is defined in \eqref{eq:canonical-Z} and that $Z_{N,b_*,v}^\delta(q)$ is its restriction to $\widehat a_N\in A_\delta$, as expressed in \eqref{eq:canonical-shell-polar}. Normalizing the integrand in \eqref{eq:canonical-Z} therefore gives
\begin{equation}
\label{eq:canonical-shell-ratio}
Z_{N,b_*,v}^\delta(q)=Z_{N,b_*,v}(q)P_{N,\delta}\quad\text{and}\quad 0<Z_{N,b_*,v}^\delta(q)\le Z_{N,b_*,v}(q).
\end{equation}
Taking logarithms, we obtain
\begin{equation}
\label{eq:log-PNdelta}
\frac1N\E\log Z_{N,b_*,v}(q)-\frac1N\E\log Z_{N,b_*,v}^\delta(q) =-\frac1N\E\log P_{N,\delta}.
\end{equation}

It remains to prove
\begin{align}\label{e.limEP=0}
    \lim_{N\to\infty}-\frac{1}{N}\E\log P_{N,\delta}=0.
\end{align}
On $\{P_{N,\delta}\ge1/2\}$, we have $0\le-\log P_{N,\delta}\le\log2$. On $\{P_{N,\delta}<1/2\}$, we use the following lower bound on the restricted partition function.
To bound $\log Z_{N,b_*,v}^\delta(q)$, let $\gamma_N^\delta$ be the probability measure on $\R^{D\times N}$ obtained by conditioning $\gamma_N$ on $\widehat a_N\in A_\delta$, namely
\begin{equation*}
\gamma_N^\delta(\mathrm d\omega):=\frac{\one_{\{\widehat a_N(\omega)\in A_\delta\}}}{\lambda_N(A_\delta)}\gamma_N(\mathrm d\omega).
\end{equation*}
The event $\{\widehat a_N\in A_\delta\}$ is invariant under $\omega\mapsto-\omega$. Hence we have $\int\omega_i\,\gamma_N^\delta(\mathrm d\omega)=0$ for every $i\in\{1,\ldots,N\}$. Consequently, for each $\alpha$ and each realization of $z^q$, we obtain
\begin{equation}
\label{eq:canonical-shell-linear-mean}
\int\sum_{i=1}^N\big(\sqrt2z_i^q(\alpha)+v\big)\cdot\omega_i\,\gamma_N^\delta(\mathrm d\omega)=0.
\end{equation}
On $\{\widehat a_N\in A_\delta\}$, we have $|\widehat a_N|\le\sqrt D+\delta$ and therefore obtain
\begin{equation*}
\left|q_k\cdot\widehat a_N+\frac12(b_*-I_D)\cdot(\widehat a_N-I_D)\right|\le |q_k|(\sqrt D+\delta)+\frac12|b_*-I_D|\delta.
\end{equation*}
Applying Jensen's inequality to each integral against $\gamma_N^\delta$, we obtain
\begin{align*}
Z_{N,b_*,v}^\delta(q)\ge\lambda_N(A_\delta)e^{-C_{q,b_*,\delta}N}\sum_\alpha v_\alpha\exp\Ll(\int\sum_{i=1}^N\big(\sqrt2z_i^q(\alpha)+v\big)\cdot\omega_i\,\gamma_N^\delta(\mathrm d\omega)\Rr)=\lambda_N(A_\delta)e^{-C_{q,b_*,\delta}N},
\end{align*}
where the equality follows from \eqref{eq:canonical-shell-linear-mean} and $\sum_\alpha v_\alpha=1$. Taking logarithms, we obtain
\begin{equation}
\label{eq:canonical-shell-log-lower}
\log Z_{N,b_*,v}^\delta(q)\ge\log\lambda_N(A_\delta)-C_{q,b_*,\delta}N.
\end{equation}
Together, \eqref{eq:canonical-shell-ratio} and \eqref{eq:canonical-shell-log-lower} imply
\begin{equation*}
\big|\log Z_{N,b_*,v}^\delta(q)\big|\le\big|\log Z_{N,b_*,v}(q)\big|+C_{q,b_*,\delta}N-\log\lambda_N(A_\delta).
\end{equation*}
Since $\lambda_N(A_\delta)>0$ for every $N\ge D$ and $\lim_{N\to\infty}N^{-1}\log\lambda_N(A_\delta)=0$ by \eqref{eq:Qhat-shell-mass}, equation \eqref{eq:canonical-log-second-moment} gives
\begin{equation}
\label{eq:canonical-shell-second-moment}
\sup_{N\ge D}\E\left[\left(\frac1N\log Z_{N,b_*,v}^\delta(q)\right)^2\right]\le C_{q,v,b_*,\delta}.
\end{equation}
On $\{P_{N,\delta}<1/2\}$, \eqref{eq:canonical-shell-ratio} gives $-\log P_{N,\delta}=\log Z_{N,b_*,v}(q)-\log Z_{N,b_*,v}^\delta(q)$. Applying Cauchy--Schwarz on this event, with \eqref{eq:canonical-log-second-moment} and \eqref{eq:canonical-shell-second-moment}, we obtain
\begin{equation*}
0\le-\frac1N\E\log P_{N,\delta}\stackrel{\text{Cauchy--Schwarz}}{\le}\frac{\log2}{N}+C_{q,v,b_*,\delta}\P(P_{N,\delta}<1/2)^{1/2}\stackrel{\eqref{eq:PNdelta}}{\le}\frac{\log2}{N}+\frac{C_{q,v,b_*,\delta}}{\sqrt N},
\end{equation*}
verifying~\eqref{e.limEP=0}.

Consequently, \eqref{eq:log-PNdelta}, \eqref{e.limEP=0}, and \eqref{eq:canonical-G-identity} give
\begin{equation}
\label{eq:canonical-shell-limit}
\lim_{N\to\infty}\frac1N\E\log Z_{N,b_*,v}^\delta(q)=G_v(b_*).
\end{equation}
Combining \eqref{eq:canonical-shell-limit}, \eqref{eq:shell-log-comparison}, and \eqref{eq:Qhat-shell-mass}, we obtain
\begin{equation*}
\limsup_{N\to\infty}\psi_{N,v}(q)\le-G_v(b_*)+C_{D,q,b_*,v}\delta.
\end{equation*}
Letting $\delta\downarrow0$ and using \eqref{eq:initial-upper} proves \eqref{eq:stiefel-determinant-formula}.
\end{proof}

Proposition~\ref{prop:stiefel-initial} evaluates $\lim_{N\to\infty}\psi_{N,v}(q)$ for every finite-step path $q$. Lemma~\ref{lem:fixed-path-Lipschitz}, with $a=I_D$ and external field $v$, gives a unique $L^1$-continuous extension $\psi_{N,v}:\mcl Q_1\to\R$. The following lemma identifies its pointwise limit with $\psi_v$ from \eqref{eq:psi-general-path}. 

\begin{lemma}[Initial condition for $q\in\mcl Q_1$]
\label{cor:initial-general-path}
For every $v\in\R^D$, the functions $\psi_{N,v}$ converge pointwise on $\mcl Q_1$ to $\psi_v$ given by \eqref{eq:psi-general-path}, that is,
\begin{equation}
\label{eq:initial-general-limit}
\lim_{N\to\infty}\psi_{N,v}(q)=\psi_v(q),\qquad q\in\mcl Q_1.
\end{equation}
Both $\psi_{N,v}$ and $\psi_v$ are $\sqrt D$-Lipschitz with respect to $|\cdot|_{L^1}$.
\end{lemma}
\begin{proof}
Fix $v\in\R^D$. We first extend the finite-step limit continuously and then identify this extension with \eqref{eq:psi-general-path}.

\smallskip
\noindent\emph{Step 1: Existence of the pointwise limit on $\mcl Q_1$.} Lemma~\ref{lem:fixed-path-Lipschitz}, with $a=I_D$ and external field $v$, gives
\begin{equation}
\label{eq:psi-Lipschitz}
|\psi_{N,v}(q)-\psi_{N,v}(q')|\le\sqrt D\int_0^1|q(s)-q'(s)|\,\mathrm d s
\end{equation}
for finite-step $q,q'$. Taking $N\to\infty$ gives the same estimate for the limits in Proposition~\ref{prop:stiefel-initial}. For $q\in\mcl Q_1$ and $n\in\N$, define the finite-step path $q_n:[0,1)\to\S_+^D$ by
\begin{equation}
\label{eq:initial-path-approximation}
q_n(s):=q(n^{-1}\lfloor ns\rfloor),\qquad s\in[0,1).
\end{equation}
Then $q_n(0)=q(0)$ and $0\preccurlyeq q_n(s)\preccurlyeq q(s)$. Since an increasing matrix path is continuous almost everywhere, dominated convergence gives $\lim_{n\to\infty}|q_n-q|_{L^1}=0$. Thus finite-step paths are dense in $\mcl Q_1$, and \eqref{eq:psi-Lipschitz} defines unique continuous extensions of $\psi_{N,v}$ and the finite-step limit. Denote the latter extension temporarily by $\widetilde\psi_v:\mcl Q_1\to\R$. For every $n\in\N$, we have
\begin{equation*}
\limsup_{N\to\infty}|\psi_{N,v}(q)-\widetilde\psi_v(q)|\le 2\sqrt D\,|q_n-q|_{L^1}+\limsup_{N\to\infty}|\psi_{N,v}(q_n)-\widetilde\psi_v(q_n)|=2\sqrt D\,|q_n-q|_{L^1}.
\end{equation*}
Letting $n\to\infty$ proves \eqref{eq:initial-general-limit} with $\widetilde\psi_v$ in place of $\psi_v$.

\smallskip
\noindent\emph{Step 2: Finiteness and attainment in \eqref{eq:psi-general-path}.} Fix $q\in\mcl Q_1$ and define $A_q:[0,1)\to\S_+^D$ by $A_q(s):=L_{\ell,q}(s)-\ell$. Definition \eqref{eq:initial-path-multiplier} gives the following expression, which is independent of $\ell$,
\begin{equation}
A_q(s)=2\Big(sq(s)-\int_0^s q(r)\,\mathrm d r\Big)=2\int_0^s(q(s)-q(r))\,\mathrm d r\succcurlyeq0,\qquad s\in[0,1).
\label{eq:initial-path-increment}
\end{equation}
For each $\varepsilon\in(0,1)$, the function $s\mapsto s^{-1}\int_0^s\tr q(r)\,\mathrm d r$ is absolutely continuous on $[\varepsilon,1]$, since $q\in L^1$. The quotient rule gives
\begin{equation*}
\frac{\mathrm d}{\mathrm d s}\left(\frac1s\int_0^s\tr q(r)\,\mathrm d r\right)\stackrel{\eqref{eq:initial-path-increment}}{=}\frac{\tr A_q(s)}{2s^2}\qquad\text{for almost every }s\in(\varepsilon,1).
\end{equation*}
Integrating over $[\varepsilon,1]$, we obtain
\begin{equation*}
\int_\varepsilon^1\frac{\tr A_q(s)}{s^2}\,\mathrm d s=2\left(\int_0^1\tr q(s)\,\mathrm d s-\frac1\varepsilon\int_0^\varepsilon\tr q(r)\,\mathrm d r\right).
\end{equation*}
Right-continuity at $0$ gives $\lim_{\varepsilon\downarrow0}\varepsilon^{-1}\int_0^\varepsilon\tr q(r)\,\mathrm d r=\tr q(0)$. Since $A_q\succcurlyeq0$, monotone convergence as $\varepsilon\downarrow0$ yields
\begin{equation}
\int_0^1\frac{\tr A_q(s)}{s^2}\,\mathrm d s=2\Big(\int_0^1\tr q(s)\,\mathrm d s-\tr q(0)\Big)<\infty.
\label{eq:initial-path-integrability}
\end{equation}
For $\ell\succcurlyeq\delta I_D$ with $\delta>0$, applying $\log(1+x)\le x$ to the eigenvalues of $\ell^{-1/2}A_q(s)\ell^{-1/2}\succcurlyeq0$ gives 
\begin{equation*}
0\le\log\frac{\det L_{\ell,q}(s)}{\det\ell}=\log\det\big(I_D+\ell^{-1/2}A_q(s)\ell^{-1/2}\big)\le\tr(\ell^{-1}A_q(s)).
\end{equation*}
Moreover, $\ell^{-1}\preccurlyeq\delta^{-1}I_D$ implies $A_q(s)^{1/2}\ell^{-1}A_q(s)^{1/2}\preccurlyeq\delta^{-1}A_q(s)$. Taking traces and using their invariance under cyclic permutations, we obtain
\begin{equation}
0\le\log\frac{\det L_{\ell,q}(s)}{\det\ell}\le\tr\big(A_q(s)^{1/2}\ell^{-1}A_q(s)^{1/2}\big)\le\delta^{-1}\tr A_q(s).
\label{eq:initial-path-log-bound}
\end{equation}
Equations \eqref{eq:initial-path-integrability}--\eqref{eq:initial-path-log-bound} prove finiteness of the integral in \eqref{eq:psi-general-path}. Let $K_{v,q}:\S_{++}^D\to\R$ denote the expression inside braces in that formula. For every $\delta>0$, the same bounds give the integrable majorant $\delta^{-1}s^{-2}\tr A_q(s)$ for the logarithmic integrand, uniformly over $\ell\succcurlyeq\delta I_D$. Dominated convergence therefore proves continuity of $K_{v,q}$ on $\S_{++}^D$. Its terms involving $q$ and $v$ are nonnegative. At $\ell=I_D$, its definition reduces to
\begin{equation*}
K_{v,q}(I_D)=2\tr q(0)+|v|^2+\int_0^1\frac1{s^2}\log\det L_{I_D,q}(s)\,\mathrm d s.
\end{equation*}
Consequently, we obtain
\begin{align}
K_{v,q}(\ell)&\ge\tr\ell-D-\log\det\ell,\notag\\
K_{v,q}(I_D)&\stackrel{\eqref{eq:initial-path-log-bound}}{\le}2\tr q(0)+|v|^2+\int_0^1\frac{\tr A_q(s)}{s^2}\,\mathrm d s\stackrel{\eqref{eq:initial-path-integrability}}{=}2\int_0^1\tr q(s)\,\mathrm d s+|v|^2.
\label{eq:initial-general-coercivity}
\end{align}
If $\lambda_1,\ldots,\lambda_D$ are the eigenvalues of $\ell$, the first lower bound equals $\sum_{d=1}^D(\lambda_d-1-\log\lambda_d)$. Each summand is nonnegative and tends to $+\infty$ as $\lambda_d\downarrow0$ or $\lambda_d\to\infty$. Hence the sublevel sets of $K_{v,q}$ are compact subsets of $\S_{++}^D$, and $K_{v,q}$ attains its minimum.

\smallskip
\noindent\emph{Step 3: Identification of \eqref{eq:psi-general-path} for finite-step paths.} First let $q$ be given by \eqref{eq:finite-step-q}. The substitution $\ell:=b+2\int_0^1q(s)\,\mathrm d s$ is a bijection from $\mathcal B_q$ onto $\S_{++}^D$. Indeed, \eqref{eq:L-recursion} gives
\begin{equation*}
\ell=\ell_1(b)\quad\text{and}\quad L_{\ell,q}(s)=\ell+2\sum_{i=1}^j m_i(q_i-q_{i-1})=\ell_{j+1}(b),\qquad s\in[m_j,m_{j+1}),\quad j\in\{0,\ldots,k\}.
\end{equation*}
Integrating on these intervals and collecting successive logarithms gives
\begin{align*}
\int_0^1\frac1{s^2}\log\frac{\det L_{\ell,q}(s)}{\det\ell}\,\mathrm d s
&=\sum_{j=1}^k\Big(\frac1{m_j}-\frac1{m_{j+1}}\Big)\log\frac{\det\ell_{j+1}(b)}{\det\ell_1(b)}\\
&=\log\det\ell_1(b)-\log\det\ell_{k+1}(b)+\sum_{i=1}^k\frac1{m_i}\log\frac{\det\ell_{i+1}(b)}{\det\ell_i(b)}.
\end{align*}
Substituting this identity and $\tr b=\tr\ell-2\int_0^1\tr q(s)\,\mathrm d s$ into \eqref{eq:G-explicit}, we obtain
\begin{equation*}
G_v(b)=\frac12K_{v,q}(\ell)-\int_0^1\tr q(s)\,\mathrm d s.
\end{equation*}
Proposition~\ref{prop:stiefel-initial} therefore proves \eqref{eq:psi-general-path} for finite-step $q$.

\smallskip
\noindent\emph{Step 4: Passage to $q\in\mcl Q_1$.} Fix $q\in\mcl Q_1$ and use $q_n$ from \eqref{eq:initial-path-approximation}. By \eqref{eq:initial-path-increment}, $A_{q_n}(s)$ converges to $A_q(s)$ for almost every $s$. Also, \eqref{eq:initial-path-integrability} and $q_n(0)=q(0)$ imply
\begin{equation*}
\lim_{n\to\infty}\int_0^1\frac{\tr A_{q_n}(s)}{s^2}\,\mathrm d s=\int_0^1\frac{\tr A_q(s)}{s^2}\,\mathrm d s.
\end{equation*}
Since $A_{q_n},A_q\succcurlyeq0$, we have $|A_{q_n}-A_q|\le\tr A_{q_n}+\tr A_q$. Apply Fatou's lemma to the nonnegative functions $s^{-2}(\tr A_{q_n}+\tr A_q-|A_{q_n}-A_q|)$. Their limit is $2s^{-2}\tr A_q$ for almost every $s\in(0,1)$, so the preceding convergence yields
\begin{equation}
\limsup_{n\to\infty}\int_0^1\frac{|A_{q_n}(s)-A_q(s)|}{s^2}\,\mathrm d s\le 2\int_0^1\frac{\tr A_q(s)}{s^2}\,\mathrm d s-2\int_0^1\frac{\tr A_q(s)}{s^2}\,\mathrm d s=0.
\label{eq:initial-increment-convergence}
\end{equation}
For $X,Y\succcurlyeq\delta I_D$, differentiation along the segment from $Y$ to $X$ gives
\begin{equation*}
|\log\det X-\log\det Y|=\left|\int_0^1\tr\big((Y+r(X-Y))^{-1}(X-Y)\big)\,\mathrm d r\right|\le\sqrt D\,\delta^{-1}|X-Y|.
\end{equation*}
The terms in $K_{v,q_n}(\ell)-K_{v,q}(\ell)$ outside the integral vanish because $q_n(0)=q(0)$. Consequently, we have
\begin{equation}
\lim_{n\to\infty}\sup_{\ell\succcurlyeq\delta I_D}|K_{v,q_n}(\ell)-K_{v,q}(\ell)|\le\sqrt D\,\delta^{-1}\lim_{n\to\infty}\int_0^1\frac{|A_{q_n}(s)-A_q(s)|}{s^2}\,\mathrm d s\stackrel{\eqref{eq:initial-increment-convergence}}{=}0.
\label{eq:initial-objective-convergence}
\end{equation}
By \eqref{eq:initial-general-coercivity} and $q_n\preccurlyeq q$, the minimizers of all $K_{v,q_n}$ and of $K_{v,q}$ lie in a single compact subset of $\S_{++}^D$. Choose $\delta>0$ such that every matrix in this set is at least $\delta I_D$. For this $\delta$, the difference of the minimum values satisfies
\begin{equation*}
\limsup_{n\to\infty}\left|\min_{\ell\succ0}K_{v,q_n}(\ell)-\min_{\ell\succ0}K_{v,q}(\ell)\right|\le\lim_{n\to\infty}\sup_{\ell\succcurlyeq\delta I_D}|K_{v,q_n}(\ell)-K_{v,q}(\ell)|\stackrel{\eqref{eq:initial-objective-convergence}}{=}0.
\end{equation*}
Taking $n\to\infty$ in the finite-step formula, using $\lim_{n\to\infty}|q_n-q|_{L^1}=0$ and Step~1, proves $\widetilde\psi_v(q)=\psi_v(q)$ with $\psi_v$ given by \eqref{eq:psi-general-path}. This proves \eqref{eq:initial-general-limit} and the asserted Lipschitz bound.
\end{proof}

The formula \eqref{eq:psi-general-path} also gives the convexity needed in Section~\ref{s:fixed-a-free-energy}. The proof uses the affine dependence of $L_{\ell,q}$ on $q$ and the concavity of $\log\det$. 

\begin{proposition}[Convexity of $\psi_v$]
\label{prop:stiefel-convexity}
For every $v\in\R^D$, $q^0,q^1\in\mcl Q_1$, and $\rho\in[0,1]$, we have
\begin{equation}
\label{eq:stiefel-convexity}
\psi_v\Ll((1-\rho)q^0+\rho q^1\Rr) \le(1-\rho)\psi_v(q^0)+\rho\psi_v(q^1). 
\end{equation}
\end{proposition}

\begin{proof}
Fix $v\in\R^D$, $q^0,q^1\in\mcl Q_1$, and $\rho\in[0,1]$, and put $q^\rho:=(1-\rho)q^0+\rho q^1$. For every $\ell\succ0$, definition \eqref{eq:initial-path-multiplier} gives
\begin{equation*}
L_{\ell,q^\rho}(s)=(1-\rho)L_{\ell,q^0}(s)+\rho L_{\ell,q^1}(s),\qquad s\in[0,1).
\end{equation*}
Since $L_{\ell,q}(s)\succcurlyeq\ell\succ0$ for every $q\in\mcl Q_1$, concavity of $\log\det$ gives
\begin{equation*}
\log\frac{\det L_{\ell,q^\rho}(s)}{\det\ell}\ge(1-\rho)\log\frac{\det L_{\ell,q^0}(s)}{\det\ell}+\rho\log\frac{\det L_{\ell,q^1}(s)}{\det\ell}.
\end{equation*}
The proof of Lemma~\ref{cor:initial-general-path} shows that these functions are integrable against $s^{-2}\,\mathrm d s$. Multiplying by $-1/(2s^2)$ and integrating, we find that
\begin{equation*}
q\mapsto\int_0^1\tr q(s)\,\mathrm d s-\tr(\ell^{-1}q(0))-\frac12\int_0^1\frac1{s^2}\log\frac{\det L_{\ell,q}(s)}{\det\ell}\,\mathrm d s
\end{equation*}
is convex on $\mcl Q_1$. Adding $-\frac12(\tr\ell-D-\log\det\ell+v^\intercal\ell^{-1}v)$ preserves convexity, since this term is independent of $q$. Taking the supremum over $\ell\succ0$ in \eqref{eq:psi-general-path} proves \eqref{eq:stiefel-convexity}.
\end{proof}

\subsection{The optimizing multiplier and regularity}
\label{s:initial-multiplier}

For a finite-step path $q$ as in \eqref{eq:finite-step-q}, recall $G_v:\mathcal B_q\to\R$ from \eqref{eq:G-explicit} and $\ell_i(b)$ from \eqref{eq:L-recursion}. Lemma~\ref{lem:initial-coercivity} gives a minimizer of $G_v$. We now prove its uniqueness and use its defining equation to compute the path derivative of $\psi_v$.

\begin{proposition}
\label{prop:initial-multiplier}
Let $b_*=b_*(q,v)\in\mathcal B_q$ be a minimizer of $G_v$, whose existence is established in Lemma~\ref{lem:initial-coercivity}. This minimizer is unique and is the unique solution $b\in\mathcal B_q$ of the following equation with $\ell_i=\ell_i(b)$:
\begin{equation}
\label{eq:initial-multiplier-moment}
I_D=\ell_{k+1}^{-1}+\ell_1^{-1}(2q_0+vv^{\intercal})\ell_1^{-1}+\sum_{i=1}^k\frac{\ell_i^{-1}-\ell_{i+1}^{-1}}{m_i}.
\end{equation}
For fixed $0=m_0<m_1<\cdots<m_k<m_{k+1}=1$, the map $(q_0,\ldots,q_k,v)\mapsto b_*(q,v)$ is continuous on the set $0\preccurlyeq q_0\preccurlyeq\cdots\preccurlyeq q_k$. Evaluate $\ell_i$ at $b_*$ and define $p(s):=p_j$ for $s\in[m_j,m_{j+1})$, where
\begin{equation}
\label{eq:initial-finite-gradient}
p_j:=\ell_1^{-1}(2q_0+vv^{\intercal})\ell_1^{-1}+\sum_{i=1}^j\frac{\ell_i^{-1}-\ell_{i+1}^{-1}}{m_i},\qquad j\in\{0,\ldots,k\}.
\end{equation}
For $c_0,\ldots,c_k\in\S^D$ and $q^\varepsilon:=\sum_{j=0}^k(q_j+\varepsilon c_j)\one_{[m_j,m_{j+1})}\in\mcl Q$, $\varepsilon\in[0,\varepsilon_0)$, with $\varepsilon_0>0$, we have
\begin{equation}\label{e.path_der}
\left.\frac{\d}{\d\varepsilon}\psi_v(q^\varepsilon)\right|_{\varepsilon=0+}=\sum_{j=0}^k(m_{j+1}-m_j)p_j\cdot c_j=\int_0^1p(s)\cdot\sum_{j=0}^kc_j\one_{[m_j,m_{j+1})}(s)\,\d s.
\end{equation}
In particular, $0\preccurlyeq p_0\preccurlyeq\cdots\preccurlyeq p_k=I_D-\ell_{k+1}^{-1}\prec I_D$.
\end{proposition}
\begin{proof}
For $b\in\mathcal B_q$, set $S_i:=\ell_i(b)^{-1}$ and $C:=2q_0+vv^{\intercal}\succcurlyeq0$, and write $V_i(c):=\tr(S_i c S_i c)$ for $c\in\S^D$. Since $\ell_i(b+\varepsilon c)=\ell_i(b)+\varepsilon c$ for every $i$, we have $(c\cdot\nabla_b)S_i=-S_icS_i$, $(c\cdot\nabla_b)\log\det\ell_i=S_i\cdot c$ and $(c\cdot\nabla_b)^2\log\det\ell_i=-V_i(c)$. Moreover, $(c\cdot\nabla_b)^2\tr(S_1C)=2\tr(S_1cS_1cS_1C)$, while $(c\cdot\nabla_b)^2\tr b=0$. Substituting these derivatives into \eqref{eq:G-explicit} gives
\begin{equation}
\label{eq:multiplier-Hessian}
(c\cdot\nabla_b)^2G_v(b)=\frac12V_{k+1}(c)+\frac12\sum_{i=1}^k\frac{V_i(c)-V_{i+1}(c)}{m_i}+\tr(S_1cS_1cS_1C).
\end{equation}
By \eqref{eq:L-recursion}, $S_i-S_{i+1}\succcurlyeq0$. Since $cS_ic$ and $cS_{i+1}c$ are positive semidefinite, we obtain
\begin{equation*}
V_i(c)-V_{i+1}(c)=\tr((S_i-S_{i+1})cS_ic)+\tr(S_{i+1}c(S_i-S_{i+1})c)\ge0.
\end{equation*}
Also, $\tr(S_1cS_1cS_1C)\ge0$, because $S_1cS_1cS_1\succcurlyeq0$ and $C\succcurlyeq0$. Consequently, \eqref{eq:multiplier-Hessian} yields
\begin{equation*}
(c\cdot\nabla_b)^2G_v(b)\ge\frac12\tr(S_{k+1}cS_{k+1}c)>0\qquad\text{if }c\ne0.
\end{equation*}
Hence $G_v$ is strictly convex on $\mathcal B_q$. For the first derivative, the terms $\tr b$, $-\log\det\ell_{k+1}$ and $\tr(S_1C)$ in twice \eqref{eq:G-explicit} contribute $I_D$, $-S_{k+1}$ and $-S_1CS_1$, respectively. Each logarithmic ratio contributes $(S_{i+1}-S_i)/m_i$. Therefore, we obtain
\begin{equation}
\label{eq:multiplier-gradient-explicit}
2\nabla_b G_v(b)=I_D-S_{k+1}-S_1CS_1-\sum_{i=1}^k\frac{S_i-S_{i+1}}{m_i}.
\end{equation}
Since the minimizer lies in the open set $\mathcal B_q$, its gradient vanishes. Equation \eqref{eq:multiplier-gradient-explicit} is then exactly \eqref{eq:initial-multiplier-moment}. Conversely, any solution of \eqref{eq:initial-multiplier-moment} is a critical point of $G_v$, hence the unique minimizer by strict convexity.

The compact-uniform coercivity in Lemma~\ref{lem:initial-coercivity} gives a common compact set containing the minimizers for data $(q_0,\ldots,q_k,v)$ in a sufficiently small compact neighborhood. To see this, choose the continuous competitor $b$ satisfying $\ell_1(b)=I_D$, which bounds the minimum values uniformly in that neighborhood. Thus every convergent sequence of data has a subsequence of minimizers converging to a minimizer of the limiting $G_v$. By uniqueness, the whole sequence converges. Hence $(q_0,\ldots,q_k,v)\mapsto b_*(q,v)$ is continuous for fixed $(m_i)_{i=1}^k$. 

To identify the path derivative~\eqref{e.path_der}, vary $q_j$ in a symmetric direction $c$, keeping $b$ and all other $q_i$ fixed. With $m_{k+1}=1$, the recursion \eqref{eq:L-recursion} gives
\begin{equation*}
(c\cdot\nabla_{q_j})\log\det\ell_i=2(\ell_i^{-1}\cdot c)\begin{cases}m_{j+1}-m_j,&i\le j,\\m_{j+1},&i=j+1,\\0,&i\ge j+2.\end{cases}
\end{equation*}
Here the case $i\le j$ is empty for $j=0$. Applying this identity to the logarithmic terms of \eqref{eq:G-explicit} yields
\begin{align}
\label{eq:initial-logarithm-derivative}
&(c\cdot\nabla_{q_j})\left[-\log\det\ell_{k+1}+\sum_{i=1}^k\frac{\log\det\ell_{i+1}-\log\det\ell_i}{m_i}\right]\notag\\
&\qquad=-2\one_{\{j=0\}}\ell_1^{-1}\cdot c-2(m_{j+1}-m_j)\sum_{i=1}^j\frac{\ell_i^{-1}-\ell_{i+1}^{-1}}{m_i}\cdot c.
\end{align}
The inverse-matrix term satisfies
\begin{equation*}
(c\cdot\nabla_{q_j})\tr\big(\ell_1^{-1}(2q_0+vv^{\intercal})\big)=2\one_{\{j=0\}}\ell_1^{-1}\cdot c-2(m_{j+1}-m_j)\big(\ell_1^{-1}(2q_0+vv^{\intercal})\ell_1^{-1}\big)\cdot c.
\end{equation*}
The terms involving $\one_{\{j=0\}}$ in \eqref{eq:initial-logarithm-derivative} and the preceding identity cancel, and \eqref{eq:initial-finite-gradient} gives
\begin{equation}
\label{eq:initial-envelope-derivative}
-(c\cdot\nabla_{q_j})G_v(b_*)=(m_{j+1}-m_j)p_j\cdot c.
\end{equation}
The envelope theorem applies to $\psi_v(q)=-\min_{b\in\mathcal B_q}G_v(b)$ because the minimizer is unique and continuous, as just proved, and the admissible domains can be identified by the common coordinate $\ell_1(b)\in\S_{++}^D$. Therefore the directional derivative of $\psi_v$ equals the partial derivative of $-G_v$ at $b_*$, and \eqref{eq:initial-envelope-derivative} proves~\eqref{e.path_der}. 
Finally, \eqref{eq:L-recursion} gives $\ell_1^{-1}\succcurlyeq\cdots\succcurlyeq\ell_{k+1}^{-1}$, and \eqref{eq:initial-multiplier-moment} gives $p_k=I_D-\ell_{k+1}^{-1}$. Thus the claimed order of the $p_j$ follows.
\end{proof}

\begin{corollary}
\label{cor:initial-uniform}
For every $L,H\in[0,\infty)$, we have
\begin{equation}
\label{eq:initial-uniform}
\lim_{N\to\infty}\sup_{\substack{q\in\mcl Q_{\infty,\preccurlyeq LI_D}\\ |v|\le H}}|\psi_{N,v}(q)-\psi_v(q)|=0.
\end{equation}
\end{corollary}
\begin{proof}
By \cite[Lemma~3.4]{chen2025free}, every $L^\infty$-bounded subset of $\mcl Q_\infty$ is precompact in $L^r$ for $1\le r<\infty$. The set $\mcl Q_{\infty,\preccurlyeq LI_D}$ from \eqref{e.Q_infty,<=} is bounded by $\sqrt D L$ in $L^\infty$. It is also closed in $L^1$: the inequalities $0\preccurlyeq q\preccurlyeq LI_D$ pass to an almost-everywhere limit and then to its right-continuous representative. Thus $\mcl Q_{\infty,\preccurlyeq LI_D}$ is compact in $L^1$. Lemma~\ref{cor:initial-general-path} gives pointwise convergence on this compact set for every $v\in\R^D$. The bound in \eqref{eq:psi-Lipschitz} is uniform in $N$. For $u\in\mathbb U_N$, the identity $uu^{\intercal}=NI_D$ gives
\begin{equation*}
\frac1N\left|(v-v')\cdot\sum_{i=1}^Nu_i\right|\le\frac1N|v-v'|\,|u|_{\mathrm{op}}\,|\mathbf1_N|=|v-v'|.
\end{equation*}
Comparing the exponents in \eqref{eq:psi-N-stiefel} at $v$ and $v'$ therefore gives
\begin{equation}
\label{eq:initial-source-Lipschitz}
|\psi_{N,v}(q)-\psi_{N,v'}(q)|\le|v-v'|\quad\text{and}\quad|\psi_v(q)-\psi_{v'}(q)|\le|v-v'|.
\end{equation}
The first inequality extends from finite-step paths to $\mcl Q_1$ by \eqref{eq:psi-Lipschitz}. The second follows by taking $N\to\infty$ and using Lemma~\ref{cor:initial-general-path}. For $\varepsilon>0$, choose a finite $\varepsilon$-net $\{(q^j,v^j):j\le J\}$ in the metric $|q-q'|_{L^1}+|v-v'|$ on this compact product. Equations \eqref{eq:psi-Lipschitz} and~\eqref{eq:initial-source-Lipschitz} give
\begin{equation*}
\sup_{\substack{q\in\mcl Q_{\infty,\preccurlyeq LI_D}\\|v|\le H}}|\psi_{N,v}(q)-\psi_v(q)|\le\max_{j\le J}|\psi_{N,v^j}(q^j)-\psi_{v^j}(q^j)|+2(\sqrt D+1)\varepsilon.
\end{equation*}
Let $N\to\infty$ for the finite net, and then let $\varepsilon\downarrow0$. This proves \eqref{eq:initial-uniform}.
\end{proof}

For $q\in\mcl Q_2$, the \textbf{subdifferential} relative to the path cone is
\begin{equation}
\label{eq:initial-subdifferential}
\partial\psi_v(q):=\left\{p\in L^2([0,1),\S^D):\psi_v(q')\ge\psi_v(q)+\int_0^1(q'-q)\cdot p\,\mathrm d s\ \text{for every }q'\in\mcl Q_2\right\}.
\end{equation}
We write $\mcl Q_2^*:=\{r\in L^2([0,1),\S^D):\int_0^1r\cdot p\,\mathrm d s\ge0\text{ for every }p\in\mcl Q_2\}$ for the \textbf{dual cone}.

\begin{lemma}
\label{lem:initial-subgradient}
For every $v\in\R^D$ and $q\in\mcl Q_2$, the set $\partial\psi_v(q)\cap\mcl Q_{\infty,\preccurlyeq I_D}$ is nonempty. If $q,q'\in\mcl Q_2$ and $q'-q\in\mcl Q_2^*$, then $\psi_v(q')\ge\psi_v(q)$. For $c\in\S_+^D$, $L\ge0$, and $s_0\in[0,1)$, we have
\begin{equation}
\label{eq:initial-growth}
0\le\psi_v\big(q+Lc\one_{[s_0,1)}\big)-\psi_v(q)\le L\operatorname{tr}(c)(1-s_0).
\end{equation}
\end{lemma}
\begin{proof}
Let $q$ be finite-step as in \eqref{eq:finite-step-q}, and let $p(s):=p_j$ on $[m_j,m_{j+1})$, with $p_j$ given by \eqref{eq:initial-finite-gradient}. Proposition~\ref{prop:initial-multiplier} gives $p\in\mcl Q_{\infty,\preccurlyeq I_D}$. Fix another finite-step path $q'\in\mcl Q_2$ and refine the partitions of $q$ and $q'$ to a common partition. Inserting a step at which $q$ is constant leaves $p$ unchanged: the inserted difference $q_j-q_{j-1}$ vanishes in \eqref{eq:L-recursion}, so the corresponding term in \eqref{eq:initial-finite-gradient} vanishes. For $\varepsilon\in[0,1]$, put $q^\varepsilon:=q+\varepsilon(q'-q)\in\mcl Q_2$. Then \eqref{e.path_der}, applied on the common partition, gives
\begin{equation*}
\lim_{\varepsilon\downarrow0}\frac{\psi_v(q^\varepsilon)-\psi_v(q)}{\varepsilon}=\langle q'-q,p\rangle_{L^2}.
\end{equation*}
By Proposition~\ref{prop:stiefel-convexity}, we have
\begin{equation*}
\psi_v(q')-\psi_v(q)\ge\frac{\psi_v(q^\varepsilon)-\psi_v(q)}{\varepsilon},\qquad\varepsilon\in(0,1].
\end{equation*}
Letting $\varepsilon\downarrow0$ proves the inequality in \eqref{eq:initial-subdifferential} for finite-step $q'$. For arbitrary $q'\in\mcl Q_2$, choose finite-step $q'_n$ with $\lim_{n\to\infty}|q'_n-q'|_{L^2}=0$. The preceding inequality and \eqref{eq:psi-Lipschitz} give
\begin{equation*}
\psi_v(q')-\psi_v(q)\ge\langle q'-q,p\rangle_{L^2}.
\end{equation*}
Hence $p\in\partial\psi_v(q)$ by \eqref{eq:initial-subdifferential}. 

For general $q\in\mcl Q_2$, choose finite-step $q_n$ with $\lim_{n\to\infty}|q_n-q|_{L^2}=0$ and select $p_n\in\partial\psi_v(q_n)\cap\mcl Q_{\infty,\preccurlyeq I_D}$ using \eqref{eq:initial-finite-gradient}. The bound $0\preccurlyeq p_n\preccurlyeq I_D$ and \cite[Lemma~3.4]{chen2025free} give a subsequence converging in $L^2$ to $p\in\mcl Q_{\infty,\preccurlyeq I_D}$. For every $q'\in\mcl Q_2$, the defining inequality \eqref{eq:initial-subdifferential} for $p_n\in\partial\psi_v(q_n)$ reads
\begin{equation*}
\psi_v(q')\ge\psi_v(q_n)+\int_0^1(q'-q_n)\cdot p_n\,\mathrm d s.
\end{equation*}
To pass to the limit, we use
\begin{equation*}
\left|\int_0^1(q'-q_n)\cdot p_n\,\mathrm d s-\int_0^1(q'-q)\cdot p\,\mathrm d s\right|\le|q'-q|_{L^2}|p_n-p|_{L^2}+\sqrt D|q_n-q|_{L^2}.
\end{equation*}
Since $\lim_{n\to\infty}|p_n-p|_{L^2}=0$, $\lim_{n\to\infty}|q_n-q|_{L^2}=0$, and $\lim_{n\to\infty}\psi_v(q_n)=\psi_v(q)$ by \eqref{eq:psi-Lipschitz}, we obtain $p\in\partial\psi_v(q)$ from \eqref{eq:initial-subdifferential}. If $q'-q\in\mcl Q_2^*$, then $\langle q'-q,p\rangle_{L^2}\ge0$ because $p\in\mcl Q_2$, so \eqref{eq:initial-subdifferential} gives $\psi_v(q')\ge\psi_v(q)$. Finally, put $d:=Lc\one_{[s_0,1)}$, and choose $p\in\partial\psi_v(q)\cap\mcl Q_{\infty,\preccurlyeq I_D}$ and $p'\in\partial\psi_v(q+d)\cap\mcl Q_{\infty,\preccurlyeq I_D}$. By \eqref{eq:initial-subdifferential}, we have
\begin{equation*}
0\le\langle d,p\rangle_{L^2}\le\psi_v(q+d)-\psi_v(q)\le\langle d,p'\rangle_{L^2}\le L(1-s_0)\tr(c),
\end{equation*}
which proves \eqref{eq:initial-growth}.
\end{proof}

\section{Convergence of the constrained free-energy}
\label{s:block-approximation}

Fix $a\in\mcl A^\circ$. We compare $\bar F_{N,a}(t,q)$ from \eqref{eq:enriched-fixed-a} with the free energy obtained by replacing $\nu_N$ with a product of Haar measures on smaller blocks. Lemma~\ref{lem:block-moment} controls $|B|^{-1}u_Bu_B^{\intercal}-I_D$ under both $\nu_N$ and the averaged Gibbs measure, where $u_B$ is the submatrix formed by the columns in a block $B$. Polar normalization then gives the comparison in Lemma~\ref{lem:block-comparison}, and convergence for the product measure gives Proposition~\ref{prop:block-limit}.

For $N\in\N$, write $\mathbf1_N:=(1,\ldots,1)^{\intercal}\in\R^N$.

\subsection{Rotation invariance}

For $N\in\N$ with $N\ge D$, define the subgroup of $\operatorname{Orth}(N)$ preserving $\mathbf1_N$ by
\begin{equation*}
K_N:=\{o\in\operatorname{Orth}(N):o\mathbf1_N=\mathbf1_N\}.
\end{equation*}
Recall the Gibbs bracket $\langle\cdot\rangle_{N,a,t,q}$ associated with \eqref{eq:fixed-a-Z}. For $o\in K_N$, we have $uo\mathbf1_N=u\mathbf1_N$ and $(uo)(vo)^{\intercal}=uv^{\intercal}$. Thus the deterministic terms in \eqref{eq:enriched-fixed-a-H} are unchanged by $u\mapsto uo$, and \eqref{eq:fixed-a-covariance-intro} and \eqref{eq:enriched-cascade-covariance} preserve the joint law of its Gaussian terms, conditionally on $\mathfrak R$. Together with Lemma~\ref{lem:stiefel-invariance}, this gives
\begin{equation}
\label{eq:block-K-invariance}
\E\langle\varphi(uo)\rangle_{N,a,t,q}=\E\langle\varphi(u)\rangle_{N,a,t,q},
\end{equation}
for every bounded measurable $\varphi:\mathbb U_N\to\R$ and every $o\in K_N$. For $B\subset\{1,\ldots,N\}$, write $u_B$ for the submatrix of $u$ formed by the columns indexed by $B$.

A probability measure $\lambda$ on $\mathbb U_N$ is \textbf{$K_N$-invariant} if $\int_{\mathbb U_N}\varphi(uo)\,\lambda(\mathrm d u)=\int_{\mathbb U_N}\varphi(u)\,\lambda(\mathrm d u)$ for every $o\in K_N$ and every bounded measurable $\varphi:\mathbb U_N\to\R$.

\begin{lemma}
\label{lem:block-moment}
There exists $C_D<\infty$ such that, for $N\in\N$ with $N\ge\max(D,3)$, every $K_N$-invariant probability measure $\lambda$ on $\mathbb U_N$, and every set $B\subset\{1,\ldots,N\}$ of cardinality $m\ge1$, we have
\begin{equation}
\label{eq:block-moment}
\int\left|\frac{u_Bu_B^{\intercal}}m-I_D\right|^2\lambda(\mathrm d u)\le\frac{C_D}{m}.
\end{equation}
\end{lemma}

\begin{proof}
\smallskip
\noindent\emph{Step 1: Averaging over $O\in K_N$.} Let $O$ have normalized Haar probability on $K_N$, characterized by $O\stackrel{\mathrm d}=o_1Oo_2$ for all $o_1,o_2\in K_N$, and write $\E_O$ for expectation with respect to $O$. The $K_N$-invariance of $\lambda$ and Tonelli's theorem give
\begin{equation}
\label{eq:block-rotation-law}
\int\left|\frac{u_Bu_B^\intercal}m-I_D\right|^2\lambda(\mathrm d u)=\int\E_O\left|\frac{(uO)_B(uO)_B^\intercal}m-I_D\right|^2\lambda(\mathrm d u).
\end{equation}
It therefore suffices to bound the expectation on the right-hand side of \eqref{eq:block-rotation-law} by $\frac{C_D}{m}$ for every fixed $u\in\mathbb U_N$, with $C_D$ independent of $u$, $N$, and $B$.

\smallskip
\noindent\emph{Step 2: Decomposition of $uO$.} Fix $u\in\mathbb U_N$ for the remainder of the proof, and define
\begin{equation*}
z:=\frac1N u\mathbf1_N=\frac1N\sum_{i=1}^N u_i\in\R^D,\qquad v:=u-z\mathbf1_N^\intercal\in\R^{D\times N},\qquad w:=vO.
\end{equation*}
Since $O\mathbf1_N=\mathbf1_N$ and $O^\intercal O=I_N$, we have $O^\intercal\mathbf1_N=O^\intercal O\mathbf1_N=\mathbf1_N$, and therefore $\mathbf1_N^\intercal O=\mathbf1_N^\intercal$. Thus we have
\begin{equation}
\label{eq:block-rotated-decomposition}
uO=(z\mathbf1_N^\intercal+v)O=z(\mathbf1_N^\intercal O)+vO=z\mathbf1_N^\intercal+w.
\end{equation}
Using $u\mathbf1_N=Nz$ and $\mathbf1_N^\intercal\mathbf1_N=N$, we compute $v\mathbf1_N=Nz-Nz=0$ and
\begin{equation*}
vv^\intercal=uu^\intercal-u\mathbf1_Nz^\intercal-z\mathbf1_N^\intercal u^\intercal+Nzz^\intercal=uu^\intercal-Nzz^\intercal.
\end{equation*}
Since $uu^\intercal=NI_D$ by \eqref{eq:stiefel-intro}, this proves
\begin{equation}
\label{eq:block-centered-identities}
v\mathbf1_N=0\qquad\text{and}\qquad vv^\intercal=N(I_D-zz^\intercal).
\end{equation}
The matrix $vv^\intercal$ is positive semidefinite, so $0\preccurlyeq I_D-zz^\intercal\preccurlyeq I_D$. The eigenvalues of $zz^\intercal$ are $|z|^2$ and zeros, which also gives $|z|\le1$.

\smallskip
\noindent\emph{Step 3: The law of $w^\intercal b$.} For $b\in\R^D$, put $r_b:=\sqrt{N b^\intercal(I_D-zz^\intercal)b}$. Equation \eqref{eq:block-centered-identities} gives
\begin{equation*}
(v^\intercal b)\cdot\mathbf1_N=b^\intercal v\mathbf1_N=0,\qquad |v^\intercal b|^2=b^\intercal vv^\intercal b\stackrel{\eqref{eq:block-centered-identities}}=r_b^2.
\end{equation*}
Let $\omega$ be uniform on $\{y\in\R^N:y\cdot\mathbf1_N=0,\ |y|=1\}$, independently of $O$, and write $\E_\omega$ and $\operatorname{Var}_\omega$ for expectation and variance with respect to $\omega$. We first verify that $O^\intercal x$ has the same law as $\omega$ for every fixed $x\in\mathbf1_N^\perp$ with $|x|=1$.

For any unit vectors $x,y\in\mathbf1_N^\perp$, extend $x$ and $y$ to orthonormal bases of $\mathbf1_N^\perp$. Mapping the basis containing $y$ to the basis containing $x$ and fixing $\mathbf1_N$ defines $o\in K_N$ with $oy=x$, hence $o^\intercal x=y$. For every bounded measurable $f:\R^N\to\R$, Haar invariance then gives
\begin{equation*}
\E_O f(O^\intercal y)=\E_O f(O^\intercal o^\intercal x)=\E_O f((oO)^\intercal x)=\E_O f(O^\intercal x).
\end{equation*}
Thus $\E_O f(O^\intercal x)$ is constant over unit vectors $x\in\mathbf1_N^\perp$. Average this equality over $x=\omega$. For every fixed $o\in K_N$, the definition of uniform probability on the sphere gives $o^\intercal\omega\stackrel{\mathrm d}=\omega$. Fubini's theorem therefore yields
\begin{equation}
\label{eq:block-Haar-sphere-test}
\E_O f(O^\intercal x)=\E_\omega\E_O f(O^\intercal\omega)=\E_O\E_\omega f(O^\intercal\omega)=\E_\omega f(\omega).
\end{equation}
If $r_b>0$, apply \eqref{eq:block-Haar-sphere-test} with $x=\frac{v^\intercal b}{r_b}$ and with $f(r_b\,\cdot)$ in place of $f$. If $r_b=0$, then $v^\intercal b=0$. In both cases, we obtain
\begin{equation}
\label{eq:block-rotated-sphere}
w^\intercal b=O^\intercal v^\intercal b\stackrel{\mathrm d}=r_b\omega,\qquad b\in\R^D.
\end{equation}
Equation \eqref{eq:block-rotated-sphere} specifies the distribution for each fixed $b$. It does not assert a simultaneous representation of all $w^\intercal b$ by one realization of $\omega$.

\smallskip
\noindent\emph{Step 4: The coordinates and moments of $\omega$.} Choose a deterministic orthonormal basis $\eta_1,\ldots,\eta_{N-1}$ of $\mathbf1_N^\perp$, and define the scalar random variables $\omega_i:=\eta_i\cdot\omega$ for $i\in\{1,\ldots,N-1\}$. Then we have
\begin{equation}
\label{eq:block-spherical-coordinates}
\omega=\sum_{i=1}^{N-1}\omega_i\eta_i,\qquad \omega_i\in\R,\qquad \sum_{i=1}^{N-1}\omega_i^2=1.
\end{equation}
The map $\omega\mapsto(\omega_1,\ldots,\omega_{N-1})$ is an isometry from $\mathbf1_N^\perp$ to $\R^{N-1}$. Thus $(\omega_1,\ldots,\omega_{N-1})$ is uniform on $\{s\in\R^{N-1}:\sum_{i=1}^{N-1}s_i^2=1\}$. The $\omega_i$ are the real coordinates of the single vector $\omega\in\R^N$. They have the same marginal distribution and are not independent.

Coordinate permutations and sign changes preserve the joint law of $(\omega_1,\ldots,\omega_{N-1})$. Hence $\E_\omega\omega_i\omega_j=0$ for $i\ne j$, and $1=\sum_{i=1}^{N-1}\E_\omega\omega_i^2=(N-1)\E_\omega\omega_1^2$. Moreover, a rotation in the first two coordinates gives $\frac{\omega_1+\omega_2}{\sqrt2}\stackrel{\mathrm d}=\omega_1$. Since $\E_\omega\omega_1^3\omega_2=\E_\omega\omega_1\omega_2^3=0$ by sign symmetry, expanding the fourth power yields
\begin{equation*}
\E_\omega\omega_1^4=\frac14\left(2\E_\omega\omega_1^4+6\E_\omega\omega_1^2\omega_2^2\right),\qquad \E_\omega\omega_1^4=3\E_\omega\omega_1^2\omega_2^2.
\end{equation*}
Squaring $\sum_{i=1}^{N-1}\omega_i^2=1$ from \eqref{eq:block-spherical-coordinates} and taking expectations gives
\begin{equation*}
1=(N-1)\E_\omega\omega_1^4+(N-1)(N-2)\E_\omega\omega_1^2\omega_2^2=(N-1)(N+1)\E_\omega\omega_1^2\omega_2^2.
\end{equation*}
Consequently, for distinct $i,j\in\{1,\ldots,N-1\}$, we get
\begin{equation}
\label{eq:block-spherical-coordinate-moments}
\E_\omega\omega_i^2=\frac1{N-1},\qquad \E_\omega\omega_i^4=\frac3{(N-1)(N+1)},\qquad \E_\omega\omega_i^2\omega_j^2=\frac1{(N-1)(N+1)}.
\end{equation}
For a symmetric operator $A:\mathbf1_N^\perp\to\mathbf1_N^\perp$, choose the basis $\eta_1,\ldots,\eta_{N-1}$ to consist of eigenvectors of $A$, with eigenvalues $\alpha_1,\ldots,\alpha_{N-1}$. Then $\omega^\intercal A\omega=\sum_{i=1}^{N-1}\alpha_i\omega_i^2$, so \eqref{eq:block-spherical-coordinate-moments} gives
\begin{equation*}
\E_\omega(\omega^\intercal A\omega)=\frac{\sum_{i=1}^{N-1}\alpha_i}{N-1},\qquad \E_\omega[(\omega^\intercal A\omega)^2]=\frac{3\sum_{i=1}^{N-1}\alpha_i^2+2\sum_{1\le i<j\le N-1}\alpha_i\alpha_j}{(N-1)(N+1)}.
\end{equation*}
Since $\sum_i\alpha_i=\tr A$, $\sum_i\alpha_i^2=\tr(A^2)$, and $2\sum_{i<j}\alpha_i\alpha_j=(\sum_i\alpha_i)^2-\sum_i\alpha_i^2$, subtracting the squared mean from the second moment yields
\begin{equation}
\label{eq:block-spherical-moments}
\E_\omega(\omega^\intercal A\omega)=\frac{\tr A}{N-1},\qquad \operatorname{Var}_\omega(\omega^\intercal A\omega)=\frac2{(N-1)(N+1)}\left(\tr(A^2)-\frac{(\tr A)^2}{N-1}\right).
\end{equation}

\smallskip
\noindent\emph{Step 5: The matrix $w_Bw_B^\intercal$.} Let $P_B:\R^N\to\R^N$ be the orthogonal projection onto $\{x\in\R^N:x_i=0\text{ for }i\notin B\}$, so $(P_Bx)_i=x_i$ for $i\in B$ and $(P_Bx)_i=0$ otherwise. Let $P:\R^N\to\R^N$ be the orthogonal projection onto $\mathbf1_N^\perp$, given by $P=I_N-\frac1N\mathbf1_N\mathbf1_N^\intercal$. Thus $P_B^\intercal=P_B^2=P_B$ and $P^\intercal=P^2=P$. The operator $PP_BP$ maps $\R^N$ into $\mathbf1_N^\perp$. To apply \eqref{eq:block-spherical-moments}, we define $A:=PP_B:\:\mathbf1_N^\perp\to\mathbf1_N^\perp$.
For $y\in\mathbf1_N^\perp$, we have $Py=y$ and thus $Ay=PP_BPy$. So, $A$ is the restriction of $PP_BP$ to $\mathbf1_N^\perp$. The vector $Ay=P(P_By)$ belongs to $\mathbf1_N^\perp$, even though $P_By$ need not belong to $\mathbf1_N^\perp$. Since $(PP_BP)^\intercal=PP_BP$, the operator $A$ is symmetric. Moreover, for $y\in\mathbf1_N^\perp$, we have
\begin{equation*}
y^\intercal Ay=y^\intercal PP_By=y^\intercal P_By=|P_By|^2\in[0,|y|^2].
\end{equation*}
Consequently, the eigenvalues of $A$ lie in $[0,1]$, and $\tr(A^2)\le\tr A$. To compute $\tr A$, use the orthogonal decomposition $\R^N=\mathbf1_N^\perp\oplus\operatorname{span}\{\mathbf1_N\}$. The operator $PP_BP$ equals $A$ on $\mathbf1_N^\perp$ and vanishes on $\operatorname{span}\{\mathbf1_N\}$, since $P\mathbf1_N=0$. Its trace on $\R^N$ therefore equals $\tr A$. Using cyclicity of the trace and $P^2=P$, we obtain
\begin{equation}
\label{eq:block-projection-traces}
\tr A=\tr(PP_BP)=\tr(P_BP^2)=\tr(P_BP)=\sum_{i\in B}P_{ii}=m\left(1-\frac1N\right),\qquad \tr(A^2)\le\tr A\le m.
\end{equation}

For $x\in\R^N$, write $x_B\in\R^m$ for the vector of standard coordinates of $x$ indexed by $B$, in increasing order. The vector $P_Bx\in\R^N$ has the same entries on $B$ and zeros outside $B$. In particular, $\omega_B\in\R^m$ and $P_B\omega\in\R^N$ have different dimensions but the same squared norm. Here $\omega_B$ uses the standard coordinates of $\omega$ in $\R^N$. The scalars $\omega_i=\eta_i\cdot\omega$ in \eqref{eq:block-spherical-coordinates} are its coordinates in the basis $\eta_1,\ldots,\eta_{N-1}$. Since $P_B^\intercal P_B=P_B$ and $P\omega=\omega$, we have
\begin{equation}
\label{eq:block-coordinate-projection}
|\omega_B|^2=|P_B\omega|^2=\omega^\intercal P_B\omega=\omega^\intercal PP_BP\omega=\omega^\intercal A\omega.
\end{equation}
Define $Q:=\frac1m w_Bw_B^\intercal-(I_D-zz^\intercal)$. Taking the standard coordinates indexed by $B$ in \eqref{eq:block-rotated-sphere} gives $w_B^\intercal b\stackrel{\mathrm d}=r_b\omega_B$ for each fixed $b\in\R^D$. Thus \eqref{eq:block-coordinate-projection} yields
\begin{equation*}
b^\intercal w_Bw_B^\intercal b=|w_B^\intercal b|^2\stackrel{\mathrm d}=r_b^2|\omega_B|^2=r_b^2\omega^\intercal A\omega.
\end{equation*}
Equations \eqref{eq:block-spherical-moments} and \eqref{eq:block-projection-traces} thus give
\begin{equation*}
\E_O[b^\intercal Qb]=\frac{r_b^2}{m}\frac{\tr A}{N-1}-b^\intercal(I_D-zz^\intercal)b=\frac{r_b^2}{N}-b^\intercal(I_D-zz^\intercal)b=0.
\end{equation*}
Using the variance formula in \eqref{eq:block-spherical-moments}, the inequality $\tr(A^2)\le m$, and $r_b^2\le N|b|^2$, we obtain
\begin{equation}
\label{eq:block-quadratic-error}
\E_O|b^\intercal Qb|^2=\frac{r_b^4}{m^2}\operatorname{Var}_\omega(\omega^\intercal A\omega)\le\frac{2N^2|b|^4}{(N-1)(N+1)m}\le\frac{3|b|^4}{m},
\end{equation}
where the last inequality uses $N\ge3$.

Let $e_1,\ldots,e_D$ be the standard basis of $\R^D$. For $i,j\in\{1,\ldots,D\}$, symmetry of $Q$ gives
\begin{equation*}
Q_{ii}=e_i^\intercal Qe_i,\qquad Q_{ij}=\frac12\left((e_i+e_j)^\intercal Q(e_i+e_j)-e_i^\intercal Qe_i-e_j^\intercal Qe_j\right)\quad(i\ne j).
\end{equation*}
By \eqref{eq:block-quadratic-error}, $\E_O|Q_{ii}|^2\le\frac3m$. For $i\ne j$, the inequality $|r+s+t|^2\le3(r^2+s^2+t^2)$ for $r,s,t\in\R$, together with $|e_i+e_j|^4=4$, gives
\begin{equation*}
\E_O|Q_{ij}|^2\le\frac34\sum_{b\in\{e_i,e_j,e_i+e_j\}}\E_O|b^\intercal Qb|^2\le\frac34\left(\frac3m+\frac3m+\frac{12}m\right)=\frac{27}{2m}.
\end{equation*}
Hence we get
\begin{equation}
\label{eq:block-matrix-error}
\E_O|Q|^2=\sum_{i,j=1}^D\E_O|Q_{ij}|^2\le\frac{3D+\frac{27}2D(D-1)}m\le\frac{14D^2}{m}.
\end{equation}

\smallskip
\noindent\emph{Step 6: The sum $\sum_{i\in B}w_i$ and completion.} Put $\mathbf1_B:=P_B\mathbf1_N$ and $\bar w_B:=\frac1m\sum_{i\in B}w_i=\frac1m w\mathbf1_B$. By \eqref{eq:block-spherical-coordinates}, \eqref{eq:block-spherical-coordinate-moments}, and $\E_\omega\omega_i\omega_j=0$ for $i\ne j$, we get
\begin{equation}
\label{eq:block-spherical-covariance}
\E_\omega\omega\omega^\intercal=\sum_{i,j=1}^{N-1}\E_\omega[\omega_i\omega_j]\eta_i\eta_j^\intercal=\frac1{N-1}\sum_{i=1}^{N-1}\eta_i\eta_i^\intercal=\frac{P}{N-1}.
\end{equation}
Here $\sum_{i=1}^{N-1}\eta_i\eta_i^\intercal=P$ because $\eta_1,\ldots,\eta_{N-1}$ is an orthonormal basis of $\mathbf1_N^\perp$. To compute $\E_\omega|\omega\cdot\mathbf1_B|^2$, first write
\begin{equation*}
|\omega\cdot\mathbf1_B|^2=(\mathbf1_B^\intercal\omega)(\omega^\intercal\mathbf1_B)=\mathbf1_B^\intercal(\omega\omega^\intercal)\mathbf1_B.
\end{equation*}
Since $\mathbf1_B$ is deterministic, expectation can be taken entrywise in $\omega\omega^\intercal$. Equation \eqref{eq:block-spherical-covariance} therefore gives
\begin{equation}
\label{eq:block-sphere-linear-moment}
\E_\omega|\omega\cdot\mathbf1_B|^2=\mathbf1_B^\intercal(\E_\omega\omega\omega^\intercal)\mathbf1_B=\frac1{N-1}\mathbf1_B^\intercal P\mathbf1_B=\frac1{N-1}|P\mathbf1_B|^2,
\end{equation}
where the last equality uses $P^\intercal P=P$. Applying \eqref{eq:block-rotated-sphere} to $b^\intercal w\mathbf1_B=(w^\intercal b)\cdot\mathbf1_B$ and then using \eqref{eq:block-sphere-linear-moment}, for every $b\in\R^D$ we obtain
\begin{equation}
\label{eq:block-linear-scalar-moment}
\E_O|b^\intercal w\mathbf1_B|^2=r_b^2\E_\omega|\omega\cdot\mathbf1_B|^2=\frac{r_b^2}{N-1}\mathbf1_B^\intercal P\mathbf1_B=\frac{r_b^2}{N-1}|P\mathbf1_B|^2.
\end{equation}
Since $\mathbf1_N\cdot\mathbf1_B=|\mathbf1_B|^2=m$, we have
\begin{equation*}
P\mathbf1_B=\mathbf1_B-\frac{m}{N}\mathbf1_N,\qquad |P\mathbf1_B|^2=m-\frac{2m^2}{N}+\frac{m^2}{N}=m-\frac{m^2}{N}.
\end{equation*}
Also, $\sum_{i=1}^D r_{e_i}^2=N\tr(I_D-zz^\intercal)$. Summing \eqref{eq:block-linear-scalar-moment} with $b=e_i$ over $i\in\{1,\ldots,D\}$ and dividing by $m^2$ gives
\begin{equation}
\label{eq:block-linear-error}
\E_O|\bar w_B|^2=\frac1{m^2}\sum_{i=1}^D r_{e_i}^2\frac{|P\mathbf1_B|^2}{N-1}=\frac{N(m-\frac{m^2}{N})}{(N-1)m^2}\tr(I_D-zz^\intercal)\le\frac{2D}{m}.
\end{equation}
Here we used $m-\frac{m^2}{N}\le m$, $\tr(I_D-zz^\intercal)\le D$, and $\frac{N}{N-1}\le2$.

By \eqref{eq:block-rotated-decomposition}, the $i$th column of $uO$ is $z+w_i$. Expanding the sum over $i\in B$ therefore gives
\begin{equation*}
\frac{(uO)_B(uO)_B^\intercal}m=\frac1m\sum_{i\in B}(z+w_i)(z+w_i)^\intercal=zz^\intercal+\frac1m w_Bw_B^\intercal+z\bar w_B^\intercal+\bar w_Bz^\intercal.
\end{equation*}
By the definition of $Q$, it follows that
\begin{equation*}
\frac{(uO)_B(uO)_B^\intercal}m-I_D=Q+z\bar w_B^\intercal+\bar w_Bz^\intercal.
\end{equation*}
Using $|A_1+A_2+A_3|^2\le3\sum_{i=1}^3|A_i|^2$ for matrices $A_i\in\R^{D\times D}$, $|z\bar w_B^\intercal|=|\bar w_Bz^\intercal|=|z|\,|\bar w_B|\le|\bar w_B|$, and \eqref{eq:block-matrix-error}--\eqref{eq:block-linear-error}, we obtain
\begin{equation*}
\E_O\left|\frac{(uO)_B(uO)_B^\intercal}m-I_D\right|^2\le3\E_O|Q|^2+6\E_O|\bar w_B|^2\le\frac{42D^2+12D}{m}.
\end{equation*}
This bound holds for every $u\in\mathbb U_N$. Taking $C_D:=42D^2+12D$ and integrating with respect to $\lambda$ in \eqref{eq:block-rotation-law} proves \eqref{eq:block-moment}.
\end{proof}

For $N,M\in\N$ with $N\ge M\ge\max(D,3)$, partition $\{1,\ldots,N\}$ into consecutive blocks $B_1,\ldots,B_J$ of sizes $m_j\in\{M,\ldots,2M-1\}$, with at most one $m_j\ne M$. With $C_D$ from Lemma~\ref{lem:block-moment}, define
\begin{equation}
\label{eq:block-shell}
A_{N,M}:=\left\{u\in\mathbb U_N:\sum_{j=1}^J\frac{m_j}{N}\left|\frac{u_{B_j}u_{B_j}^{\intercal}}{m_j}-I_D\right|\le\frac{2\sqrt{C_D}}{\sqrt M}\right\}.
\end{equation}
By Lemma~\ref{lem:stiefel-invariance} and \eqref{eq:block-K-invariance}, Lemma~\ref{lem:block-moment} applies both to $\lambda=\nu_N$ and to the marginal law of $u$ under $\E\langle\cdot\rangle_{N,a,t,q}$. In both cases, Cauchy--Schwarz gives
\begin{equation*}
\int\sum_{j=1}^J\frac{m_j}{N}\left|\frac{u_{B_j}u_{B_j}^{\intercal}}{m_j}-I_D\right|\,\lambda(\mathrm d u)\le\sum_{j=1}^J\frac{m_j}{N}\left(\frac{C_D}{m_j}\right)^{1/2}\le\frac{\sqrt{C_D}}{\sqrt M}.
\end{equation*}
Markov's inequality therefore yields
\begin{equation}
\label{eq:block-shell-probability}
\E\la\one_{A_{N,M}}\ra_{N,a,t,q}\ge\frac12\quad\text{and}\quad \nu_N(A_{N,M})\ge\frac12.
\end{equation}

\subsection{Comparison with a vector spin model}
\label{s:block-restriction}

For $M\in\N$ with $M\ge\max(D,3)$ and $n\in\N$, let $\bar F^{[M]}_{nM,a}(t,q)$ be defined by replacing $\nu_{nM}$ in \eqref{eq:fixed-a-Z}--\eqref{eq:enriched-fixed-a} with $\nu_M^{\otimes n}$. We identify $(u^1,\ldots,u^n)\in\mathbb U_M^n$ with the concatenated matrix $u\in\mathbb U_{nM}$. The identities $uu^{\intercal}=nMI_D$ and $|ah\cdot u\mathbf1_{nM}|\le nM|ah|$, together with Jensen's inequality, give
\begin{equation}
\label{eq:block-finite-bounds}
-|ah|\le\bar F^{[M]}_{nM,a}(t,q)\le t\xi_a(I_D)+q(1)\cdot I_D.
\end{equation}
For the lower bound in \eqref{eq:block-finite-bounds}, average the partition function over the Gaussian fields. The corrections cancel their variances, and the remaining integrand is at most $e^{nM|ah|}$. For the upper bound, apply Jensen's inequality inside the partition function and then average: the Gaussian terms have mean zero, and $\int u_i\,\nu_M^{\otimes n}(\mathrm d u)=0$ for every $i\in\{1,\ldots,nM\}$.

Recall the pairwise disjoint consecutive sets $B_1,\ldots,B_J$ used in \eqref{eq:block-shell}: their union is $\{1,\ldots,N\}$, and $|B_j|=m_j$ for $j\in\{1,\ldots,J\}$. Define $\bar F^{[M]}_{N,a}(t,q)$ by replacing $\nu_N$ in \eqref{eq:enriched-fixed-a} with $\bigotimes_{j=1}^J\nu_{m_j}$. The identity $uu^{\intercal}=NI_D$ holds under either reference measure, so the corrections agree.

\begin{lemma}
\label{lem:block-comparison}
For every $t\ge0$ and $q\in\mcl Q_\infty$, there exists $C<\infty$, independent of $M\in\N$ with $M\ge\max(D,3)$, such that
\begin{equation}
\label{eq:block-comparison-limsup}
\limsup_{N\to\infty}\left|\bar F_{N,a}(t,q)-\bar F^{[M]}_{N,a}(t,q)\right|\le C M^{-1/4}.
\end{equation}
For fixed $a$ and $D$, the constant $C$ is bounded when $t$, $|ah|$, and $|q|_{L^\infty}$ range over bounded sets.
\end{lemma}

\begin{proof}
\smallskip
\noindent\emph{Step 1: Restriction to $A_{N,M}$.} We first assume that $q$ is finite-step. For a probability measure $\eta_N$ on $\mathbb U_N$, write $Z_{N,a,\eta_N}(t,q)$ for the partition function obtained by replacing $\nu_N$ with $\eta_N$ in \eqref{eq:fixed-a-Z}. Lemma~\ref{lem:total-concentration-continuous} gives
\begin{equation}
\label{eq:block-total-concentration}
\E\left|\log Z_{N,a,\eta_N}(t,q)-\E\log Z_{N,a,\eta_N}(t,q)\right|\le C_{a,t,q,D}\sqrt N.
\end{equation}
Indeed, by $uu^{\intercal}=NI_D$ and \eqref{eq:can-concentration-variance}, we may take $V=2Nt\xi_a(I_D)+2N\tr q(1)$, independently of $\eta_N$. The terms $ah\cdot\sum_i u_i$, $-Nt\xi_a(I_D)$, and $-Nq(1)\cdot I_D$ are deterministic, and the expectation in \eqref{eq:block-total-concentration} includes the cascade weights.

Write $Z_N$ for $Z_{N,a,\nu_N}(t,q)$ and $Z_N^A$ for the same integral restricted to $A_{N,M}$, without renormalizing Haar measure. Set $Y_N:=\log Z_N^A-\log Z_N\le0$. By \eqref{eq:block-shell-probability}, we have $\E e^{Y_N}\ge1/2$. Since $e^{Y_N}\le1/4$ on $\{Y_N<-\log4\}$ and $e^{Y_N}\le1$ everywhere, we obtain
\begin{equation*}
\frac12\le\frac14+\frac34\P(Y_N\ge-\log4),
\end{equation*}
so $\P(Y_N\ge-\log4)\ge1/3$. Apply \eqref{eq:block-total-concentration} to $\nu_N$ and $\nu_N(\cdot\mid A_{N,M})$. The deterministic factor $\nu_N(A_{N,M})$ disappears upon centering. Consequently, $\E|Y_N-\E Y_N|\le C\sqrt N$, and
\begin{equation*}
\frac13\max\{0,-\log4-\E Y_N\}\le\E|Y_N-\E Y_N|\le C\sqrt N.
\end{equation*}
Since $\E Y_N\le0$, we conclude that
\begin{equation}
\label{eq:block-restriction-cost}
0\le\E\log Z_N-\E\log Z_N^A\le\log4+C\sqrt N.
\end{equation}

\smallskip
\noindent\emph{Step 2: The image of the conditioned Haar measure.} Recall that $u_{B_j}\in\R^{D\times m_j}$ is the matrix formed by the columns of $u$ indexed by $B_j$, in increasing order. The matrices $u_{B_j}u_{B_j}^\intercal$ are positive definite for $\nu_N$-almost every $u$. To verify this, let $G\in\R^{D\times N}$ have independent standard Gaussian entries. The matrix $\sqrt N(GG^\intercal)^{-1/2}G$ has law $\nu_N$ by right orthogonal invariance and Lemma~\ref{lem:stiefel-invariance}. Since $m_j\ge D$, the first $D$ columns of $G_{B_j}$ form a square Gaussian matrix with nonzero determinant almost surely. Thus $G_{B_j}$, and hence $u_{B_j}$, has row rank $D$ almost surely, for every $j\in\{1,\ldots,J\}$.

Define $T_{N,M}:\mathbb U_N\to\prod_{j=1}^J\mathbb U_{m_j}$ on this set of full $\nu_N$-measure by
\begin{equation}
\label{eq:block-normalization-map}
(T_{N,M}u)_{B_j}:=\bigg(\frac{u_{B_j}u_{B_j}^\intercal}{m_j}\bigg)^{-1/2}u_{B_j},\qquad j\in\{1,\ldots,J\}.
\end{equation}
Multiplying the expression in \eqref{eq:block-normalization-map} by its transpose gives $(T_{N,M}u)_{B_j}(T_{N,M}u)_{B_j}^\intercal=m_jI_D$, so $(T_{N,M}u)_{B_j}\in\mathbb U_{m_j}$. 

For each $j\in\{1,\ldots,J\}$, sample $O_j$ independently from Haar measure on $\operatorname{Orth}(m_j)$. Define $O:=\operatorname{diag}(O_1,\ldots,O_J)$. Then $O\in\operatorname{Orth}(N)$. For fixed $u$ satisfying $u_{B_j}u_{B_j}^\intercal\succ0$ for every $j\in\{1,\ldots,J\}$, we have
\begin{equation}
\label{eq:block-rotation-equivariance}
(uO)_{B_j}=u_{B_j}O_j,\qquad (uO)_{B_j}(uO)_{B_j}^\intercal=u_{B_j}u_{B_j}^\intercal,\quad\text{and}\quad T_{N,M}(uO)=(T_{N,M}u)O.
\end{equation}
For any fixed $v_j\in\mathbb U_{m_j}$, the rows of $m_j^{-1/2}v_j$ are orthonormal. Extending them to an orthonormal basis of $\R^{m_j}$ gives $o_j\in\operatorname{Orth}(m_j)$ with $v_j=\sqrt{m_j}(I_D\ 0)o_j$. Haar invariance and \eqref{eq:nu-pushforward} imply
\begin{equation*}
v_jO_j=\sqrt{m_j}(I_D\ 0)o_jO_j\stackrel{\mathrm d}=\sqrt{m_j}(I_D\ 0)O_j\sim\nu_{m_j}.
\end{equation*}
Take $v_j=(T_{N,M}u)_{B_j}$. Since the $O_j$ are independent, \eqref{eq:block-rotation-equivariance} gives, for every bounded measurable $\Phi:\prod_{j=1}^J\mathbb U_{m_j}\to\R$,
\begin{equation}
\label{eq:block-product-rotation-average}
\E_O\Phi(T_{N,M}(uO))=\int\Phi(v^1,\ldots,v^J)\prod_{j=1}^J\nu_{m_j}(\mathrm d v^j).
\end{equation}
In particular, this expectation does not depend on $u$.

By \eqref{eq:block-rotation-equivariance}, multiplication by any block diagonal $o=\operatorname{diag}(o_1,\ldots,o_J)$ preserves every $u_{B_j}u_{B_j}^\intercal$ and hence the event $A_{N,M}$ in \eqref{eq:block-shell}. Lemma~\ref{lem:stiefel-invariance} therefore implies that $\nu_N(\cdot\mid A_{N,M})$ is invariant under $u\mapsto uo$. Averaging this invariance over $O$ and using \eqref{eq:block-product-rotation-average}, we obtain
\begin{equation*}
\int\Phi(T_{N,M}u)\,\nu_N(\mathrm d u\mid A_{N,M})=\int\E_O\Phi(T_{N,M}(uO))\,\nu_N(\mathrm d u\mid A_{N,M})=\int\Phi\,\mathrm d\Big(\bigotimes_{j=1}^J\nu_{m_j}\Big).
\end{equation*}
Since this identity holds for every such $\Phi$, it proves
\begin{equation}
\label{eq:block-normalization-pushforward}
(T_{N,M})_\#\nu_N(\cdot\mid A_{N,M})=\bigotimes_{j=1}^J\nu_{m_j}.
\end{equation}
The following estimates are used only under $\nu_N(\cdot\mid A_{N,M})$, so we omit the null set where some $u_{B_j}u_{B_j}^\intercal$ is singular.

For $\nu_N$-almost every $u\in A_{N,M}$, write $u_{B_j}=(m_j^{-1}u_{B_j}u_{B_j}^{\intercal})^{1/2}(T_{N,M}u)_{B_j}$. Applying $(\sqrt s-1)^2\le|s-1|$ for $s\ge0$ to the eigenvalues of $m_j^{-1}u_{B_j}u_{B_j}^{\intercal}$ gives
\begin{equation}
\label{eq:block-distance}
\frac1N|u-T_{N,M}u|^2=\sum_{j=1}^J\frac{m_j}{N}\tr\left(\bigg(\Big(\frac{u_{B_j}u_{B_j}^{\intercal}}{m_j}\Big)^{1/2}-I_D\bigg)^2\right)\le\frac{2\sqrt{DC_D}}{\sqrt M}.
\end{equation}
The last inequality follows from Cauchy--Schwarz applied to the eigenvalues and \eqref{eq:block-shell}. Since $|u|=|T_{N,M}u|=\sqrt{ND}$, the decomposition $uv^{\intercal}-(T_{N,M}u)(T_{N,M}v)^{\intercal}=(u-T_{N,M}u)v^{\intercal}+(T_{N,M}u)(v-T_{N,M}v)^{\intercal}$ and \eqref{eq:block-distance} give
\begin{equation}
\label{eq:block-overlap-comparison}
\left|\frac{uv^{\intercal}}N-\frac{(T_{N,M}u)(T_{N,M}v)^{\intercal}}N\right|\le 2\sqrt2 D^{3/4}C_D^{1/4}M^{-1/4},\qquad\text{for }\nu_N^{\otimes2}\text{-a.e. }(u,v)\in A_{N,M}^2.
\end{equation}

\smallskip
\noindent\emph{Step 3: Gaussian interpolation and the limit in the path.} Conditionally on the cascade, take independent centered Gaussian fields $X_0,X_1$ on $A_{N,M}\times\operatorname{supp}\mathfrak R$ with covariances
\begin{align*}
 \E X_0(u,\alpha)X_0(v,\alpha')&=2Nt\xi_a(uv^{\intercal}/N)+2Nq(\alpha\wedge\alpha')\cdot uv^{\intercal}/N,\\
 \E X_1(u,\alpha)X_1(v,\alpha')&=2Nt\xi_a((T_{N,M}u)(T_{N,M}v)^{\intercal}/N)+2Nq(\alpha\wedge\alpha')\cdot (T_{N,M}u)(T_{N,M}v)^{\intercal}/N.
\end{align*}
In this step, put $T:=T_{N,M}$. Let $Z_s$ be the partition function with reference probability $\nu_N(\cdot\mid A_{N,M})$, Gaussian field $\sqrt{1-s}X_0+\sqrt sX_1$, linear term $ah\cdot((1-s)u+sTu)\mathbf1_N$, and correction $-Nt\xi_a(I_D)-Nq(1)\cdot I_D$. Write $\langle\cdot\rangle_s$ for its Gibbs bracket and, for two replicas, set $R:=uv^{\intercal}/N$ and $R':=(Tu)(Tv)^{\intercal}/N$. Then, the standard Gaussian interpolation computation gives
\begin{equation}
\label{eq:block-interpolation-derivative}
-\frac1N\frac{\mathrm d}{\mathrm d s}\E\log Z_s=\E\la t(\xi_a(R')-\xi_a(R))+q(\alpha\wedge\alpha')\cdot(R'-R)\ra_s-\frac1N\E\la ah\cdot(Tu-u)\mathbf1_N\ra_s.
\end{equation}
By Cauchy--Schwarz, $|\mathbf1_N|=\sqrt N$, and \eqref{eq:block-distance}, we have
\begin{equation*}
\frac1N\left|\E\la ah\cdot(Tu-u)\mathbf1_N\ra_s\right|\le |ah|\E\la N^{-1/2}|Tu-u|\ra_s\le\sqrt2|ah|(DC_D)^{1/4}M^{-1/4}.
\end{equation*}
Both $R$ and $R'$ have operator norm at most one. By \eqref{eq:mixture-summability} and $|aRa|\le |a|^2|R|_{\mathrm{op}}\le D$, the gradient of $\xi_a$ is bounded on this convex operator-norm ball. Using \eqref{eq:block-overlap-comparison}, we obtain
\begin{equation}
\label{eq:block-interpolation-bound}
\sup_{0<s<1}\frac1N\left|\frac{\mathrm d}{\mathrm d s}\E\log Z_s\right|\le C_{a,D}(t+|q|_{L^\infty}+|ah|)M^{-1/4}.
\end{equation}
The Gaussian covariances at $s=0,1$ and the pushforward identity \eqref{eq:block-normalization-pushforward} give
\begin{equation}
\label{eq:block-interpolation-endpoints}
\E\log Z_0=\E\log Z_N^A-\log\nu_N(A_{N,M})\quad\text{and}\quad \E\log Z_1=-N\bar F^{[M]}_{N,a}(t,q).
\end{equation}
Combining \eqref{eq:block-interpolation-bound}--\eqref{eq:block-interpolation-endpoints} with \eqref{eq:block-shell-probability} and \eqref{eq:block-restriction-cost}, we obtain
\begin{equation}
\label{eq:block-comparison}
\left|\bar F_{N,a}(t,q)-\bar F^{[M]}_{N,a}(t,q)\right|\le C M^{-1/4}+\frac{C_{a,t,q,D}}{\sqrt N}+\frac{\log8}{N}.
\end{equation}
Here $C$ is bounded when $t$, $|ah|$, and $|q|_{L^\infty}$ range over bounded sets. The constant $C_{a,t,q,D}$ may also depend on the steps of $q$, but both constants are independent of $N$ and $M$.

For bounded $q\in\mcl Q$, put $q_k(s):=q(2^{-k}\lfloor2^ks\rfloor)$ for $s\in[0,1)$. Then $q_k\in\mcl Q_\infty$ is finite-step, $|q_k|_{L^\infty}\le|q|_{L^\infty}$, and $\lim_{k\to\infty}|q_k-q|_{L^1}=0$ by dominated convergence. The argument of Lemma~\ref{lem:fixed-path-Lipschitz} applies to both reference measures, since $uu^{\intercal}=NI_D$. Applying \eqref{eq:fixed-enriched-q-Lipschitz} twice gives
\begin{equation*}
\limsup_{N\to\infty}|\bar F_{N,a}(t,q)-\bar F^{[M]}_{N,a}(t,q)|\le 2\sqrt D|q-q_k|_{L^1}+CM^{-1/4}.
\end{equation*}
The constant $C$ in \eqref{eq:block-comparison} is uniform in $k$, by the bound on $|q_k|_{L^\infty}$. Taking $N\to\infty$ at fixed $k$ and then $k\to\infty$ proves \eqref{eq:block-comparison-limsup}. No concentration bound uniform in the number of steps is needed.
\end{proof}

\subsection{Existence of limit free energy}

We apply \cite[Theorem~1.1]{ho2026convergence} for compactly supported reference probabilities. See also \cite[Theorem~1.1]{chenHoIssa2026formula}. Assumption \eqref{eq:mixture-summability} supplies the required absolute convergence on the whole matrix space. We identify the reference probability, Hamiltonian, and enrichment before applying these results.

For $R\in\R^{DM\times DM}$ and $i,j\in\{1,\ldots,M\}$, let $R_{ij}\in\R^{D\times D}$ be the submatrix with row indices $(i-1)D+1,\ldots,iD$ and column indices $(j-1)D+1,\ldots,jD$. Define $L_M:\R^{DM\times DM}\to\R^{D\times D}$ and $\boldsymbol{\xi}_M:\R^{DM\times DM}\to\R$ by
\begin{equation}
\label{eq:block-contraction}
L_M(R):=\frac1M\sum_{i=1}^M R_{ii}\quad\text{and}\quad \boldsymbol{\xi}_M(R):=M\xi_a(L_M(R)).
\end{equation}
Define the linear isometry $\iota_M:\R^{D\times M}\to\R^{DM}$ by stacking the columns:
\begin{equation}
\label{eq:block-vectorization}
\iota_M(u):=(u_1^\intercal,\ldots,u_M^\intercal)^\intercal,\qquad (\iota_M(u))_{(i-1)D+d}=u_{d,i}\quad(i\in\{1,\ldots,M\},\ d\in\{1,\dots,D\}).
\end{equation}
Its inverse splits a vector in $\R^{DM}$ into $M$ consecutive groups of $D$ entries. In particular, $\iota_M(\mathbb U_M)$ is compact and $|\iota_M(u)|^2=\tr(uu^\intercal)=DM$ for $u\in\mathbb U_M$.

For $M\in\N$ with $M\ge\max(D,3)$, define a probability measure $\mu_{M,ah}$ on $\R^{DM}$ by
\begin{equation}
\label{eq:block-vector-prior}
Z_{M,ah}:=\int_{\mathbb U_M}e^{\sum_{i=1}^M ah\cdot u_i}\,\nu_M(\mathrm d u)\quad\text{and}\quad \mu_{M,ah}:=(\iota_M)_\#\left(Z_{M,ah}^{-1}e^{\sum_{i=1}^M ah\cdot u_i}\nu_M(\mathrm d u)\right).
\end{equation}
The bound $|\sum_i ah\cdot u_i|\le M|ah|$ gives $0<Z_{M,ah}<\infty$.

For $r\in[1,\infty]$, let $\mcl Q_r^{DM}$ be the analogue of $\mcl Q_r$ with values in $\S_+^{DM}$. For $\mathbf{q}\in\mcl Q_\infty^{DM}$, put $\mathbf{q}(1):=\lim_{s\uparrow1}\mathbf{q}(s)$. Conditionally on $\mathfrak R$, let $\zeta^{\mathbf{q}}:\supp\mathfrak R\to\R^{DM}$ be a centered Gaussian field with covariance $\E[\zeta^{\mathbf{q}}(\alpha)\zeta^{\mathbf{q}}(\alpha')^\intercal\mid\mathfrak R]=\mathbf{q}(\alpha\wedge\alpha')$. Define $\widetilde\psi_M:\mcl Q_\infty^{DM}\to\R$ by
\begin{equation}
\label{eq:block-vector-initial}
\widetilde\psi_M(\mathbf{q}):=-\E\log\iint\exp\big(\sqrt2\zeta^{\mathbf{q}}(\alpha)\cdot v-\mathbf{q}(1)\cdot vv^\intercal\big)\mu_{M,ah}(\mathrm d v)\mathfrak R(\mathrm d\alpha).
\end{equation}
We use its $L^1$-continuous extension to $\mcl Q_1^{DM}$, supplied by \cite[Proposition~5.1]{chen2025free}. The constants may depend on $M$ because $|v|^2=DM$ on $\supp\mu_{M,ah}$.

\begin{lemma}
\label{lem:block-product-convergence}
For every $M\in\N$ with $M\ge\max(D,3)$, the limit
\begin{equation}
\label{eq:block-product-full-limit}
\bar F_{\infty,a}^{[M]}(t,q):=\lim_{N\to\infty}\bar F^{[M]}_{N,a}(t,q)
\end{equation}
exists for all $t\ge0$ and $q\in\mcl Q_\infty$, and defines a function $\bar F_{\infty,a}^{[M]}:[0,\infty)\times\mcl Q_\infty\to\R$.
\end{lemma}

\begin{proof}
\smallskip
\noindent\emph{Step 1: The reference probability and Hamiltonian for $N=nM$.} Fix $M\in\N$ with $M\ge\max(D,3)$. We use the product reference probability $\mu_{M,ah}^{\otimes n}$ on $(\R^{DM})^n$. For a configuration $\sigma=(\sigma^1,\ldots,\sigma^n)$, put $u^b:=\iota_M^{-1}(\sigma^b)$ for $b\in\{1,\ldots,n\}$ and concatenate these matrices as $U:=(u^1\mid\cdots\mid u^n)\in\R^{D\times nM}$. Define the centered Gaussian Hamiltonian $\widetilde H_n:\R^{DM\times n}\to\R$ by $\widetilde H_n(\sigma):=H_{nM}(aU)$.

For another configuration $\tau=(\tau^1,\ldots,\tau^n)$, let $v^b:=\iota_M^{-1}(\tau^b)$ and $V:=(v^1\mid\cdots\mid v^n)$. Define $R_n:\R^{DM\times n}\times\R^{DM\times n}\to\R^{DM\times DM}$ by $R_n(\sigma,\tau):=\frac1n\sum_{b=1}^n\sigma^b(\tau^b)^\intercal$. Its $(i,j)$ block is $\frac1n\sum_bu_i^b(v_j^b)^\intercal$, so \eqref{eq:block-contraction} gives
\begin{equation}
\label{eq:block-superspin-overlap}
L_M(R_n(\sigma,\tau))=\frac{UV^\intercal}{nM}=\frac1{nM}\sum_{b=1}^n\sum_{i=1}^M u_i^b(v_i^b)^\intercal\quad\text{and}\quad \E\widetilde H_n(\sigma)\widetilde H_n(\tau)=n\boldsymbol{\xi}_M(R_n(\sigma,\tau)).
\end{equation}
The covariance equality follows from \eqref{e.EHH=xi} and \eqref{eq:block-contraction}. Under $\mu_{M,ah}^{\otimes n}$, every $u^b$ belongs to $\mathbb U_M$, so $UU^\intercal=nMI_D$ and $L_M(R_n(\sigma,\sigma))=I_D$.

\smallskip
\noindent\emph{Step 2: The enriched Hamiltonian and partition function.} For $\mathbf{q}\in\mcl Q_\infty^{DM}$, conditionally on $\mathfrak R$, let $\zeta_b^{\mathbf{q}}:\supp\mathfrak R\to\R^{DM}$, $b\in\{1,\ldots,n\}$, be independent copies of $\zeta^{\mathbf{q}}$, independent of $\widetilde H_n$. Their covariance is
\begin{equation*}
\E[\zeta_b^{\mathbf{q}}(\alpha)\zeta_c^{\mathbf{q}}(\alpha')^\intercal\mid\mathfrak R]=\one_{\{b=c\}}\mathbf{q}(\alpha\wedge\alpha'),\qquad b,c\in\{1,\ldots,n\}.
\end{equation*}
For $t\ge0$, define the real-valued Hamiltonian on $\R^{DM\times n}\times\supp\mathfrak R$, its partition function, and its free energy by
\begin{align}
\widetilde H_n^{t,\mathbf{q}}(\sigma,\alpha)&:=\sqrt{2t}\widetilde H_n(\sigma)-nt\boldsymbol{\xi}_M(R_n(\sigma,\sigma))+\sqrt2\sum_{b=1}^n\zeta_b^{\mathbf{q}}(\alpha)\cdot\sigma^b-n\mathbf{q}(1)\cdot R_n(\sigma,\sigma),\label{eq:block-vector-Hamiltonian}\\
\widetilde Z_n(t,\mathbf{q})&:=\iint e^{\widetilde H_n^{t,\mathbf{q}}(\sigma,\alpha)}\mu_{M,ah}^{\otimes n}(\mathrm d\sigma)\mathfrak R(\mathrm d\alpha),\label{eq:block-vector-partition}\\
\widetilde F_n(t,\mathbf{q})&:=-\frac1n\E\log\widetilde Z_n(t,\mathbf{q}).\label{eq:block-vector-free-energy}
\end{align}
Now take $\mathbf{q}(s):=\operatorname{diag}(q(s),\ldots,q(s))$ for $s\in[0,1)$. For any $R\in\R^{DM\times DM}$, the block definitions give
\begin{equation}
\label{eq:block-enrichment-contraction}
\mathbf{q}(s)\cdot R=\sum_{i=1}^M q(s)\cdot R_{ii}=Mq(s)\cdot L_M(R).
\end{equation}
We can realize $\zeta_b^{\mathbf{q}}$ by stacking the $D$-dimensional fields $z_{(b-1)M+i}^q$, $i\in\{1,\ldots,M\}$, from \eqref{eq:enriched-cascade-covariance}. Consequently, we have
\begin{gather*}
\sum_{b=1}^n\zeta_b^{\mathbf{q}}(\alpha)\cdot\sigma^b=\sum_{b=1}^n\sum_{i=1}^M z_{(b-1)M+i}^q(\alpha)\cdot u_i^b,
\\
nt\boldsymbol{\xi}_M(R_n(\sigma,\sigma))=nMt\xi_a(I_D)\quad\text{and}\quad n\mathbf{q}(1)\cdot R_n(\sigma,\sigma)=nMq(1)\cdot I_D.
\end{gather*}
Finally, the density in \eqref{eq:block-vector-prior} contributes $\sum_{b,i}ah\cdot u_i^b-n\log Z_{M,ah}$ to the exponent. Comparing \eqref{eq:block-vector-Hamiltonian}--\eqref{eq:block-vector-partition} with \eqref{eq:enriched-fixed-a-H}--\eqref{eq:fixed-a-Z}, we obtain
\begin{equation}
\label{eq:block-vector-finite-energy}
\bar F^{[M]}_{nM,a}(t,q)=\frac{\widetilde F_n(t,\mathbf{q})-\log Z_{M,ah}}M.
\end{equation}
Thus the reference probability and every term of the Hamiltonian have been identified exactly.

\smallskip
\noindent\emph{Step 3: Convergence for $M\mid N$.} The reference probability $\mu_{M,ah}$ has compact support, and \eqref{eq:block-superspin-overlap} realizes $\boldsymbol{\xi}_M$ as a Gaussian covariance. Also, $\boldsymbol{\xi}_M(0)=0$. Since $|L_M(R)|\le|R|$, the degree-$k$ term of $\boldsymbol{\xi}_M$ satisfies
\begin{equation*}
M\left|\langle\Delta_k,(aL_M(R)a)^{\otimes k}\rangle\right|\le M|\Delta_k|\big(|a|_{\mathrm{op}}^2|R|\big)^k,\qquad R\in\R^{DM\times DM}.
\end{equation*}
Assumption \eqref{eq:mixture-summability} therefore gives absolute convergence at every $R$. By \cite[Theorem~1.1]{ho2026convergence}, the limit $\lim_{n\to\infty}\widetilde F_n(t,\mathbf{q})$ exists for every $t\ge0$ and $\mathbf{q}\in\mcl Q_\infty^{DM}$. This result applies directly to compactly supported reference probabilities. At $t=0$, \eqref{eq:block-vector-initial} and the product identity \cite[(1.6)]{chen2025free} give $\widetilde F_n(0,\mathbf{q})=\widetilde\psi_M(\mathbf{q})$ for every $n\in\N$. For $\mathbf{q}(s)=\operatorname{diag}(q(s),\ldots,q(s))$, equation \eqref{eq:block-vector-finite-energy} yields
\begin{equation}
\label{eq:block-normalized-free-energy}
\lim_{n\to\infty}\bar F^{[M]}_{nM,a}(t,q)=\frac{\lim_{n\to\infty}\widetilde F_n(t,\mathbf{q})-\log Z_{M,ah}}M\quad\text{and}\quad \psi_{M,ah}(q)=\frac{\widetilde\psi_M(\mathbf{q})-\log Z_{M,ah}}M.
\end{equation}

\smallskip
\noindent\emph{Step 4: Convergence for $M\nmid N$.} If $N$ is not a multiple of $M$, the partition in \eqref{eq:block-shell} has $n$ blocks of size $M$ and one block of size $r\in(M,2M)$, where $N=nM+r$. For fixed $M$ and sufficiently large $N$, we have $n\ge1$. Let $R$ and $R'$ be the overlaps of these two groups, normalized by $nM$ and $r$, respectively, so $|R|_{\mathrm{op}},|R'|_{\mathrm{op}}\le1$. If $L_a$ is a Lipschitz constant of $\xi_a$ on the operator-norm unit ball, then
\begin{align}
\left|N\xi_a\left(\frac{N-r}{N}R+\frac rN R'\right)-(N-r)\xi_a(R)\right|&\le r\sup_{|S|_{\mathrm{op}}\le1}|\xi_a(S)|+rL_a|R'-R|\notag\\
&\le r\sup_{|S|_{\mathrm{op}}\le1}|\xi_a(S)|+2r\sqrt D L_a.
\label{eq:block-exceptional-covariance}
\end{align}
On $\nu_M^{\otimes n}\otimes\nu_r$, compare the size-$N$ model to an independent size-$nM$ model that ignores the last $r$ columns. By \eqref{eq:block-exceptional-covariance}, the interaction covariance difference is at most $Cr$. The enrichment covariance difference is $2r q(\alpha\wedge\alpha')\cdot R'$, bounded by $2r\sqrt D|q|_{L^\infty}$. The corrections differ by $rt\xi_a(I_D)+rq(1)\cdot I_D$, and the omitted term $ah\cdot\sum_{i\in B_{\mathrm{last}}}u_i$ has absolute value at most $r|ah|$, where $B_{\mathrm{last}}$ denotes the last block of $r$ columns. Gaussian interpolation, with independent fields as in the proof of \eqref{eq:block-interpolation-derivative}, gives
\begin{equation}
\label{eq:block-exceptional-logarithm}
\left|N\bar F^{[M]}_{N,a}(t,q)-(N-r)\bar F^{[M]}_{N-r,a}(t,q)\right|\le Cr\le2CM.
\end{equation}
Divide \eqref{eq:block-exceptional-logarithm} by $N$ and use convergence along multiples of $M$ from \eqref{eq:block-normalized-free-energy}. This proves \eqref{eq:block-product-full-limit}.
\end{proof}

\begin{proposition}
\label{prop:block-limit}
For every $a\in\mcl A^\circ$, $t\ge0$, and $q\in\mcl Q_\infty$, the limit 
\begin{equation}
\label{eq:block-fixed-limit}
\bar F_{\infty,a}(t,q):=\lim_{N\to\infty}\bar F_{N,a}(t,q)
\end{equation}
exists. Moreover, for some $C=C(a,t,q,h,D)<\infty$, we have 
\begin{equation}
\label{eq:block-limit-error}
|\bar F_{\infty,a}(t,q)-\bar F_{\infty,a}^{[M]}(t,q)|\le C M^{-1/4}.
\end{equation}
\end{proposition}

\begin{proof}
Lemma~\ref{lem:block-product-convergence} and \eqref{eq:block-comparison-limsup} give, for every $M\ge\max(D,3)$,
\begin{equation*}
\limsup_{N\to\infty}\bar F_{N,a}(t,q)\le \bar F_{\infty,a}^{[M]}(t,q)+CM^{-1/4}\quad\text{and}\quad \liminf_{N\to\infty}\bar F_{N,a}(t,q)\ge \bar F_{\infty,a}^{[M]}(t,q)-CM^{-1/4}.
\end{equation*}
Subtracting these inequalities, we obtain
\begin{equation*}
0\le\limsup_{N\to\infty}\bar F_{N,a}(t,q)-\liminf_{N\to\infty}\bar F_{N,a}(t,q)\le2CM^{-1/4}.
\end{equation*}
Letting $M\to\infty$ proves \eqref{eq:block-fixed-limit}. Taking the limit in \eqref{eq:block-comparison-limsup} then gives \eqref{eq:block-limit-error}.
\end{proof}

\section{Critical-point representation}
\label{s:block-critical}

Proposition~\ref{prop:block-limit} gives the limit $\bar F_{\infty,a}(t,q)$ appearing in~\eqref{eq:block-fixed-limit}. We prove its critical-point representation in Proposition~\ref{prop:block-critical} at the end of this section. Lemma~\ref{lem:block-initial-gradient} first bounds the derivative $\partial_q\psi_{M,x}$ in \eqref{eq:block-gradient-bound}. We then state the critical identities for $\bar F_{\infty,a}^{[M]}$, and finally prove an estimate that permits passage of these identities to the limit $M\to\infty$.

Recall $\mathbb U_M$ from \eqref{eq:stiefel-intro}, with $N=M$, and its invariant probability measure $\nu_M$ from \eqref{eq:nu-pushforward}. For $M\in\N$ with $M\ge D$ and $x\in\R^D$, the function $\psi_{M,x}:\mcl Q_1\to\R$ is defined by \eqref{eq:psi-N-stiefel}, with $N=M$ and $v=x$, and extended by \eqref{eq:psi-Lipschitz}. For $q\in\mcl Q_\infty$, we write $\langle\cdot\rangle_{M,x,q}$ for the corresponding Gibbs expectation on $(u,\alpha)$. Thus the exponent is $\sqrt2\sum_{i=1}^M z_i^q(\alpha)\cdot u_i+x\cdot\sum_{i=1}^M u_i-Mq(1)\cdot I_D$, and $(u,\alpha)$ and $(u',\alpha')$ are independent replicas conditionally on the Gaussian fields and cascade. The expectation $\E$ averages these fields and the cascade. We use the continuous representation from Section~\ref{s:enriched-properties}, so $\alpha\wedge\alpha'$ takes values in $[0,1]$.
Recall also that $\mcl Q_{\infty,\preccurlyeq I_D}$ in \eqref{e.Q_infty,<=} consists of the paths $p\in\mcl Q_\infty$ with $p(s)\preccurlyeq I_D$ for every $s\in[0,1)$. 

Following \cite[Definition~2.1]{chen2025free} and \cite[Section~2.2]{chenHoIssa2026formula}, a function $G:\mcl Q_2\to\R$ is \textbf{Fr\'echet differentiable} relative to $\mcl Q_2$ at $q\in\mcl Q_2$ if there is a unique $p\in L^2([0,1),\S^D)$ such that 
\begin{equation*}
\lim_{\substack{q'\in\mcl Q_2\setminus\{q\}\\|q'-q|_{L^2}\to0}}\frac{|G(q')-G(q)-\langle p,q'-q\rangle_{L^2}|}{|q'-q|_{L^2}}=0.
\end{equation*}
We call $p$ the \textbf{Fr\'echet derivative} and write $\partial_qG(q):=p$. The same definition applies to $\mcl Q_2^{DM}$ with $\S^D$ replaced by $\S^{DM}$.

\begin{lemma}[Derivative of $\psi_{M,x}$]
\label{lem:block-initial-gradient}
For $M\in\N$ with $M\ge D$ and $x\in\R^D$, the function $\psi_{M,x}$ is Fr\'echet differentiable relative to $\mcl Q_2$, and $\partial_q\psi_{M,x}(q)\in\mcl Q_{\infty,\preccurlyeq I_D}$ for every $q\in\mcl Q_2$. For $q\in\mcl Q_\infty$, if $p_{M,x}:=\partial_q\psi_{M,x}(q)$, then, for every $b\in\R^D$ and almost every $s\in[0,1)$, we have
\begin{equation}
\label{eq:block-gradient-bound}
b^{\intercal}p_{M,x}(s)b=\frac1M\E\left\langle (u^{\intercal}b)\cdot((u')^{\intercal}b)\,\middle|\,\alpha\wedge\alpha'=s\right\rangle_{M,x,q}\le |b|^2.
\end{equation}
The conditional expression denotes disintegration of the averaged two-replica law $\E\langle\cdot\rangle_{M,x,q}$. Moreover, there exists $C_{M,D}<\infty$, independent of $x$, such that
\begin{equation}
\label{eq:block-gradient-continuity}
|\partial_q\psi_{M,x}(q')-\partial_q\psi_{M,x}(q)|_{L^2}\le C_{M,D}|q'-q|_{L^2},\qquad q,q'\in\mcl Q_2.
\end{equation}
\end{lemma}

\begin{proof}
\smallskip
\noindent\emph{Step 1: Differentiability and monotonicity.} For this proof, replace $ah$ by $x$ in the definitions of $Z_{M,ah}$, $\mu_{M,ah}$, and $\widetilde\psi_M$ in \eqref{eq:block-vector-prior} and \eqref{eq:block-vector-initial}. These definitions apply for every $M\ge D$. The vector $\iota_M(u)=(u_1^\intercal,\ldots,u_M^\intercal)^\intercal$ has squared norm $DM$, so $\mu_{M,x}$ has compact support. Apply \cite[Proposition~5.1, (5.6)--(5.8)]{chen2025free} with one vector spin of dimension $DM$. Compact support suffices, as explained in \cite[Section~1.6]{chen2025free}. The constants depend on the radius $\sqrt{DM}$, but not on the density of $\mu_{M,x}$. Thus $\widetilde\psi_M$ is Fr\'echet differentiable on $\mcl Q_2^{DM}$, its derivative belongs to $\mcl Q_\infty^{DM}$, and its derivative is Lipschitz in $L^2$ with a constant depending only on $M$ and $D$.

For $q\in\mcl Q_2$, put $\mathbf{q}(s):=\operatorname{diag}(q(s),\ldots,q(s))$ and $\mathbf{p}:=\partial_{\mathbf{q}}\widetilde\psi_M(\mathbf{q})$. The second identity in \eqref{eq:block-normalized-free-energy}, with $ah=x$, extends from bounded paths to $\mcl Q_2$ by $L^1$ continuity. Its derivative in a direction $\kappa\in L^2([0,1),\S^D)$ is represented by
\begin{equation}
\label{eq:block-gradient-contraction}
\langle\partial_q\psi_{M,x}(q),\kappa\rangle_{L^2}=\frac1M\langle \mathbf{p},\operatorname{diag}(\kappa,\ldots,\kappa)\rangle_{L^2}=\langle L_M(\mathbf{p}),\kappa\rangle_{L^2}.
\end{equation}
Here $L_M(\mathbf{p})(s)=M^{-1}\sum_{i=1}^M\mathbf{p}_{ii}(s)$, as in \eqref{eq:block-contraction}. Each diagonal block of a positive-semidefinite matrix is positive semidefinite. Consequently, $\mathbf{p}(v)-\mathbf{p}(u)\succcurlyeq0$ implies $L_M(\mathbf{p})(v)-L_M(\mathbf{p})(u)\succcurlyeq0$, and $\partial_q\psi_{M,x}(q)=L_M(\mathbf{p})$ belongs to $\mcl Q_\infty$. The bounds $|\operatorname{diag}(\kappa,\ldots,\kappa)|_{L^2}=\sqrt M|\kappa|_{L^2}$ and $|L_M(\mathbf{p})|_{L^2}\le M^{-1/2}|\mathbf{p}|_{L^2}$ transfer the derivative continuity from $\widetilde\psi_M$ to \eqref{eq:block-gradient-continuity}.

\smallskip
\noindent\emph{Step 2: The conditional identity and the bound by $I_D$.} Assume that $q\in\mcl Q_\infty$. By \cite[(5.7)]{chen2025free} and \eqref{eq:block-gradient-contraction}, we have 
\begin{equation}
\label{eq:block-gradient-pairing}
\int_0^1\kappa(s)\cdot p_{M,x}(s)\,\mathrm d s=\frac1M\E\la\kappa(\alpha\wedge\alpha')\cdot u(u')^\intercal\ra_{M,x,q},\qquad \kappa\in L^2([0,1),\S^D).
\end{equation}
The factor $M^{-1}$ comes from the normalization per column in \eqref{eq:psi-N-stiefel}. The sum of the diagonal blocks of $\iota_M(u)\iota_M(u')^\intercal$ is $\sum_{i=1}^Mu_i(u_i')^\intercal=u(u')^\intercal$, which gives the matrix on the right side of \eqref{eq:block-gradient-pairing}.

By \cite[Proposition~4.8]{chen2025free}, the averaged law of $\alpha\wedge\alpha'$ under this Gibbs measure is uniform on $[0,1]$. Take $\kappa(s)=\varphi(s)bb^\intercal$ in \eqref{eq:block-gradient-pairing}, where $\varphi:[0,1)\to\R$ is bounded and measurable. Disintegration gives
\begin{equation*}
\int_0^1\varphi(s)b^\intercal p_{M,x}(s)b\,\mathrm d s=\int_0^1\varphi(s)\frac1M\E\left\langle (u^\intercal b)\cdot((u')^\intercal b)\,\middle|\,\alpha\wedge\alpha'=s\right\rangle_{M,x,q}\,\mathrm d s.
\end{equation*}
Since this holds for every $\varphi$, the two integrands agree almost everywhere. For every $u,u'\in\mathbb U_M$, Cauchy--Schwarz and $uu^\intercal=u'(u')^\intercal=MI_D$ give
\begin{equation*}
\frac1M(u^{\intercal}b)\cdot((u')^{\intercal}b)\le\frac1M|u^{\intercal}b|\,|(u')^{\intercal}b|=|b|^2.
\end{equation*}
Conditional expectation preserves this inequality, proving \eqref{eq:block-gradient-bound}. Apply it to the countable set $b\in\Q^D$ and remove the union of the corresponding null sets. Continuity in $b$ then gives $p_{M,x}(s)\preccurlyeq I_D$ for almost every $s$. For any $s<1$, choose $s_n\downarrow s$ outside this null set. Right continuity gives $p_{M,x}(s)=\lim_{n\to\infty}p_{M,x}(s_n)\preccurlyeq I_D$.

\smallskip
\noindent\emph{Step 3: Paths in $\mcl Q_2$.} For $q\in\mcl Q_2$, define $q_k:[0,1)\to\S_+^D$ by $q_k(s):=q(\min\{s,1-1/k\})$ for $k\in\{2,3,\ldots\}$. These paths are bounded and increasing. Also, $0\preccurlyeq q_k(s)\preccurlyeq q(s)$ and $q_k(s)=q(s)$ for all sufficiently large $k$, for every fixed $s<1$. Dominated convergence gives $\lim_{k\to\infty}|q_k-q|_{L^2}=0$. Equation \eqref{eq:block-gradient-continuity} implies $\lim_{k\to\infty}|\partial_q\psi_{M,x}(q_k)-\partial_q\psi_{M,x}(q)|_{L^2}=0$. Pass to an almost-everywhere convergent subsequence. The inequalities $0\preccurlyeq\partial_q\psi_{M,x}(q_k)\preccurlyeq I_D$ pass to the limit, and right continuity again extends them to every $s<1$. Thus $\partial_q\psi_{M,x}(q)\in\mcl Q_{\infty,\preccurlyeq I_D}$.
\end{proof}

Define $\theta:\S_+^D\to\R$ by 
\begin{equation}
\label{eq:theta-definition}
\theta(r):=r\cdot\nabla\xi(r)-\xi(r).
\end{equation}
For $a\in\mcl A^\circ$ and $r\in\S_+^D$, the chain rule applied to \eqref{eq:fixed-a-xi-intro} gives $\nabla\xi_a(r)=a\nabla\xi(ara)a$. Using $\theta$ from \eqref{eq:theta-definition}, define $\theta_a:\S_+^D\to\R$ by
\begin{equation}
\label{eq:theta-a-definition}
\theta_a(r):=r\cdot\nabla\xi_a(r)-\xi_a(r)=\theta(ara),\qquad r\in\S_+^D.
\end{equation}

\begin{lemma}[Critical-point representation for $\bar F_{\infty,a}^{[M]}(t,q)$]
\label{lem:block-product-critical}
For every $a\in\mcl A^\circ$, $M\in\N$ with $M\ge\max(D,3)$, $t\ge0$, and $q\in\mcl Q_\infty$, there exists $p_M\in\mcl Q_{\infty,\preccurlyeq I_D}$ such that, with $\pi_M:=q+t\nabla\xi_a(p_M)$, we have
\begin{equation}
\label{eq:block-product-critical}
p_M=\partial_q\psi_{M,ah}(\pi_M)\quad\text{and}\quad \bar F_{\infty,a}^{[M]}(t,q)=\psi_{M,ah}(\pi_M)-t\int_0^1\theta_a(p_M(s))\,\mathrm d s.
\end{equation}
\end{lemma}

\begin{proof}
\smallskip
\noindent\emph{Step 1: The critical identities in dimension $DM$.} Recall from \eqref{eq:block-contraction} that $L_M(R)=M^{-1}\sum_{i=1}^MR_{ii}$ and $\boldsymbol{\xi}_M(R)=M\xi_a(L_M(R))$. The probability $\mu_{M,ah}$ in \eqref{eq:block-vector-prior} is the image under $\iota_M$ of $Z_{M,ah}^{-1}e^{ah\cdot\sum_i u_i}\nu_M(\mathrm d u)$, and $\widetilde\psi_M$ is its initial condition in \eqref{eq:block-vector-initial}. Put $\mathbf{q}(s):=\operatorname{diag}(q(s),\ldots,q(s))$ for $s\in[0,1)$. Step~3 of the proof of Lemma~\ref{lem:block-product-convergence} gives convergence of $\widetilde F_n(t,\mathbf{q})$ for every $t\ge0$ and $\mathbf{q}\in\mcl Q_\infty^{DM}$. The uniform $L^1$ bound in \cite[Proposition~3.1]{chen2025free}, at fixed $M$, extends this convergence to $\mathbf{q}\in\mcl Q_2^{DM}$. We can therefore apply \cite[Proposition~7.3]{chen2025free} to $\boldsymbol{\xi}_M$ and $\mu_{M,ah}$. Compact support is sufficient by \cite[Section~1.6]{chen2025free}. 
It gives $\mathbf{p}_M,\mathbf{q}_M'\in\mcl Q_\infty^{DM}$ such that
\begin{equation}
\label{eq:block-full-critical-paths}
\mathbf{p}_M=\partial_{\mathbf{q}}\widetilde\psi_M(\mathbf{q}_M')\quad\text{and}\quad \mathbf{q}_M'=\mathbf{q}+t\nabla\boldsymbol{\xi}_M(\mathbf{p}_M).
\end{equation}
The paths $\mathbf{p}_M$ and $\mathbf{q}_M'$ take values in $\S_+^{DM}$ and need not initially be block diagonal. By \eqref{eq:block-vector-finite-energy} and \eqref{eq:block-product-full-limit}, the limit free energy per vector spin is $M\bar F_{\infty,a}^{[M]}(t,q)+\log Z_{M,ah}$. Substituting $\mathbf{q}-\mathbf{q}_M'=-t\nabla\boldsymbol{\xi}_M(\mathbf{p}_M)$ into the critical value from \cite[Proposition~7.3]{chen2025free} gives
\begin{equation}
\label{eq:block-full-critical-value}
M\bar F_{\infty,a}^{[M]}(t,q)+\log Z_{M,ah}=\widetilde\psi_M(\mathbf{q}_M')-t\int_0^1\big(\mathbf{p}_M(s)\cdot\nabla\boldsymbol{\xi}_M(\mathbf{p}_M(s))-\boldsymbol{\xi}_M(\mathbf{p}_M(s))\big)\,\mathrm d s.
\end{equation}

\smallskip
\noindent\emph{Step 2: The path $p_M=L_M(\mathbf{p}_M)$.} Define $p_M:[0,1)\to\S_+^D$ by $p_M(s):=L_M(\mathbf{p}_M(s))$. For $R,V\in\S^{DM}$, the definition of $\boldsymbol{\xi}_M$ gives
\begin{equation*}
\left.\frac{\mathrm d}{\mathrm d\varepsilon}\boldsymbol{\xi}_M(R+\varepsilon V)\right|_{\varepsilon=0}=M\nabla\xi_a(L_M(R))\cdot L_M(V)=\sum_{i=1}^M\nabla\xi_a(L_M(R))\cdot V_{ii}.
\end{equation*}
Since the last expression contains only the diagonal blocks of $V$, we obtain
\begin{equation}
\label{eq:block-critical-contraction}
\nabla\boldsymbol{\xi}_M(\mathbf{p}_M)=\operatorname{diag}\big(\nabla\xi_a(p_M),\ldots,\nabla\xi_a(p_M)\big).
\end{equation}
Equations \eqref{eq:block-full-critical-paths} and \eqref{eq:block-critical-contraction} show that $\mathbf{q}_M'=\operatorname{diag}(\pi_M,\ldots,\pi_M)$ with $\pi_M:=q+t\nabla\xi_a(p_M)$. In particular, $\pi_M\in\mcl Q_\infty$ because it is a diagonal block of $\mathbf{q}_M'$. Differentiating the second identity in \eqref{eq:block-normalized-free-energy}, exactly as in \eqref{eq:block-gradient-contraction}, yields
\begin{equation*}
\partial_q\psi_{M,ah}(\pi_M)=\frac1M\sum_{i=1}^M\big(\partial_{\mathbf{q}}\widetilde\psi_M(\mathbf{q}_M')\big)_{ii}=L_M(\mathbf{p}_M)=p_M.
\end{equation*}
Lemma~\ref{lem:block-initial-gradient} now gives $p_M\in\mcl Q_{\infty,\preccurlyeq I_D}$. This proves the first identity in \eqref{eq:block-product-critical} and ensures that $\theta_a(p_M)$ is defined.

\smallskip
\noindent\emph{Step 3: The value identity.} Recall $\theta_a(r)=r\cdot\nabla\xi_a(r)-\xi_a(r)$ from \eqref{eq:theta-a-definition}. 
By \eqref{eq:block-critical-contraction}, the contraction $\mathbf{p}_M\cdot\nabla\boldsymbol{\xi}_M(\mathbf{p}_M)$ equals $\sum_i(\mathbf{p}_M)_{ii}\cdot\nabla\xi_a(p_M)$. Since $\sum_i(\mathbf{p}_M)_{ii}=Mp_M$, we have
\begin{equation}
\label{eq:block-critical-trace}
\mathbf{p}_M\cdot\nabla\boldsymbol{\xi}_M(\mathbf{p}_M)-\boldsymbol{\xi}_M(\mathbf{p}_M)=M\big(p_M\cdot\nabla\xi_a(p_M)-\xi_a(p_M)\big)=M\theta_a(p_M).
\end{equation}
Substitute \eqref{eq:block-critical-trace} into \eqref{eq:block-full-critical-value} and divide by $M$. The result is
\begin{equation*}
\bar F_{\infty,a}^{[M]}(t,q)=\frac{\widetilde\psi_M(\mathbf{q}_M')-\log Z_{M,ah}}M-t\int_0^1\theta_a(p_M(s))\,\mathrm d s.
\end{equation*}
The second identity in \eqref{eq:block-normalized-free-energy}, applied at $\pi_M$, identifies the first term on the right with $\psi_{M,ah}(\pi_M)$. This proves the second identity in \eqref{eq:block-product-critical}.
\end{proof}

\begin{lemma}[Local semiconcavity of $\psi_{M,x}$]
\label{lem:block-semiconcavity}
Fix $c>0$. Suppose that $\pi\in\mcl Q_\infty$ satisfies $\pi(0)\succcurlyeq cI_D$, $\pi(v)-\pi(u)\succcurlyeq c(v-u)I_D$ for $0\le u\le v<1$, and $|\pi|_{L^\infty}\le c^{-1}$. For every $\kappa\in C^1([0,1],\S^D)$, $M\in\N$ with $M\ge D$, $x\in\R^D$, and $s\in\R$ satisfying $|s|(|\kappa(0)|+|\kappa'|_{L^\infty})\le c/2$, we have
\begin{equation}
\label{eq:block-semiconcavity}
\psi_{M,x}(\pi+s\kappa)\le\psi_{M,x}(\pi)+s\langle\partial_q\psi_{M,x}(\pi),\kappa\rangle_{L^2}+\frac{s^2}{c}(|\kappa(0)|+|\kappa'|_{L^\infty})^2.
\end{equation}
\end{lemma}

\begin{proof}
The constants in \cite[Propositions~3.7 and~3.8]{chen2025free} are already independent of the number of spins. The assumptions in \cite[(3.21) and~(3.32)]{chen2025free} additionally require a bound on the ratio of the largest and smallest eigenvalues of every path increment. We adapt the square-root argument in \cite[proof of Proposition~3.7]{chen2025free}. The improvement is in Step~2 below: a term involving the second derivative of the square root is nonpositive, so we retain its sign instead of bounding its absolute value.

\smallskip
\noindent\emph{Step 1: Finite-step approximation and square roots.} Put $L:=|\kappa(0)|+|\kappa'|_{L^\infty}$. If $L=0$, then $\kappa=0$ and \eqref{eq:block-semiconcavity} is an equality. Assume $L>0$. For $|s|L\le c/2$, we have $\pi(0)+s\kappa(0)\succcurlyeq(c/2)I_D$ and
\begin{equation*}
(\pi+s\kappa)(v)-(\pi+s\kappa)(u)\succcurlyeq(c-|s|L)(v-u)I_D\succcurlyeq\frac c2(v-u)I_D,\qquad 0\le u\le v<1.
\end{equation*}
Thus $\pi+s\kappa\in\mcl Q_\infty$ throughout this interval of $s$.

Fix $k\in\N$, put $m_j:=j/(k+1)$ for $0\le j\le k+1$, and let $\pi^{(k)}$ and $\kappa^{(k)}$ equal $\pi(m_j)$ and $\kappa(m_j)$ on $[m_j,m_{j+1})$. Define the positive-definite square roots by
\begin{align*}
\gamma_0(s)&:=\big(\pi(0)+s\kappa(0)\big)^{1/2},\\
\gamma_j(s)&:=\big(\pi(m_j)-\pi(m_{j-1})+s(\kappa(m_j)-\kappa(m_{j-1}))\big)^{1/2},\qquad 1\le j\le k.
\end{align*}
Their squares depend affinely on $s$, so differentiating twice gives
\begin{equation}
\label{eq:block-square-root-identities}
\pi^{(k)}+s\kappa^{(k)}=\sum_{j=0}^k\gamma_j(s)^2\one_{[m_j,1)}\qquad\text{and}\qquad\gamma_j\gamma_j''+\gamma_j''\gamma_j=-2(\gamma_j')^2.
\end{equation}
Here and below, primes denote derivatives in $s$. The preceding lower bounds give $\gamma_0(s)^2\succcurlyeq(c/2)I_D$ and $\gamma_j(s)^2\succcurlyeq cI_D/(2(k+1))$ for $j\ge1$. By \cite[(3.23)]{chen2025free}, or by differentiating $\gamma_j^2$ in an orthonormal eigenbasis of $\gamma_j$, we have
\begin{equation*}
|\gamma_j'(s)|^2\le\frac{\big|\frac{\mathrm d}{\mathrm d s}(\gamma_j(s)^2)\big|^2}{4\lambda_{\min}(\gamma_j(s)^2)}.
\end{equation*}
Since $|\kappa(m_j)-\kappa(m_{j-1})|\le L/(k+1)$, summing this inequality yields
\begin{equation}
\label{eq:block-square-root-mesh-bounds}
\sum_{j=0}^k|\gamma_j'(s)|^2\le\frac{|\kappa(0)|^2}{2c}+\frac{k+1}{2c}\sum_{j=1}^k|\kappa(m_j)-\kappa(m_{j-1})|^2\le\frac{L^2}{2c}\left(1+\frac{k}{k+1}\right)\le\frac{L^2}{c}.
\end{equation}

\smallskip
\noindent\emph{Step 2: The second derivative along $\gamma(s)$.} For $\gamma=(\gamma_0,\ldots,\gamma_k)\in(\S^D)^{k+1}$, define
\begin{equation*}
G(\gamma):=\psi_{M,x}\bigg(\sum_{j=0}^k\gamma_j^2\one_{[m_j,1)}\bigg).
\end{equation*}
Following \cite[proof of Proposition~3.7, Step~1, from (3.27) to (3.29)]{chen2025free}, with $N=M$, $t=0$, and reference probability proportional to $e^{x\cdot\sum_i u_i}\nu_M(\mathrm d u)$, and using $uu^\intercal=MI_D$, we obtain, for every $b=(b_0,\ldots,b_k)\in(\S^D)^{k+1}$,
\begin{equation}
\label{eq:block-square-root-concavity}
(b\cdot\nabla)^2G(\gamma)\le2\sum_{j=0}^k|b_j|^2.
\end{equation}
The finite Gaussian recursion justifies these derivatives, as in \cite[proof of Proposition~3.7, Step~1]{chen2025free}.

Along $\gamma(s)$, put $p_s:=\partial_q\psi_{M,x}(\pi^{(k)}+s\kappa^{(k)})$. Lemma~\ref{lem:block-initial-gradient} gives $0\preccurlyeq p_s\preccurlyeq I_D$. Differentiating the argument of $\psi_{M,x}$ in the direction $\gamma''(s)$ and then using \eqref{eq:block-square-root-identities}, we obtain
\begin{align}
\label{eq:block-square-root-sign}
\gamma''(s)\cdot\nabla G(\gamma(s))=\sum_{j=0}^k\int_{m_j}^1p_s(v)\cdot\big(\gamma_j\gamma_j''+\gamma_j''\gamma_j\big)\,\mathrm d v=-2\sum_{j=0}^k\int_{m_j}^1p_s(v)\cdot(\gamma_j')^2\,\mathrm d v\le0.
\end{align}
The last inequality uses $p_s(v)\succcurlyeq0$ and $(\gamma_j')^2\succcurlyeq0$. In \cite[proof of Proposition~3.7, Step~2]{chen2025free}, the corresponding term is bounded in absolute value by $2\sum_j|\gamma_j|\,|\gamma_j''|$. That estimate introduces the ratio of the largest and smallest eigenvalues of the increments. Equation~\eqref{eq:block-square-root-sign} removes this term from the upper bound. Only the lower eigenvalue bounds used in \eqref{eq:block-square-root-mesh-bounds} remain.

The chain rule, \eqref{eq:block-square-root-concavity}, and \eqref{eq:block-square-root-sign} now give
\begin{equation*}
\frac{\mathrm d^2}{\mathrm d s^2}G(\gamma(s))=(\gamma'(s)\cdot\nabla)^2G(\gamma(s))+\gamma''(s)\cdot\nabla G(\gamma(s))\le2\sum_{j=0}^k|\gamma_j'(s)|^2\le\frac{2L^2}{c}.
\end{equation*}
Since $G(\gamma(s))=\psi_{M,x}(\pi^{(k)}+s\kappa^{(k)})$, its derivative at $s=0$ is $\langle\partial_q\psi_{M,x}(\pi^{(k)}),\kappa^{(k)}\rangle_{L^2}$. Taylor's formula thus proves \eqref{eq:block-semiconcavity} for $\pi^{(k)}$ and $\kappa^{(k)}$, with remainder $c^{-1}L^2s^2$.

\smallskip
\noindent\emph{Step 3: Passage to the original paths.} The bounded monotone path $\pi$ is continuous almost everywhere, so dominated convergence gives $\lim_{k\to\infty}|\pi^{(k)}-\pi|_{L^2}=0$. Also, $|\kappa^{(k)}-\kappa|_{L^\infty}\le L/(k+1)$. Equation~\eqref{eq:psi-Lipschitz} therefore gives convergence of both function values in the finite-step inequality. At fixed $M$ and $x$, Lemma~\ref{lem:block-initial-gradient} gives
\begin{align*}
\left|\langle\partial_q\psi_{M,x}(\pi^{(k)}),\kappa^{(k)}\rangle_{L^2}-\langle\partial_q\psi_{M,x}(\pi),\kappa\rangle_{L^2}\right|\le C_{M,D}|\pi^{(k)}-\pi|_{L^2}|\kappa|_{L^2}+\sqrt D\,|\kappa^{(k)}-\kappa|_{L^1}.
\end{align*}
The right-hand side tends to zero as $k\to\infty$. The remainder $c^{-1}L^2s^2$ is independent of $k$, $M$, and $x$, so this limit proves \eqref{eq:block-semiconcavity}.
\end{proof}

Recall from \eqref{eq:block-fixed-limit} that $\bar F_{\infty,a}(t,q)=\lim_{N\to\infty}\bar F_{N,a}(t,q)$. For $a\in\mcl A^\circ$, $t\ge0$, and $q\in\mcl Q_2$, define $\sP_{a,t,q}:\mcl Q_{\infty,\preccurlyeq I_D}\to\R$ by
\begin{equation*}
\sP_{a,t,q}(p):=\psi_{ah}\big(q+t\nabla\xi_a(p)\big)-t\int_0^1\theta_a(p(s))\,\mathrm d s.
\end{equation*}
Here $q+t\nabla\xi_a(p)\in\mcl Q_2$ by Lemma~\ref{lem:mixture-monotonicity} and the chain rule for $\xi_a$, and $\theta_a$ is defined in \eqref{eq:theta-a-definition}. We suppress the dependence on $h$.
The subdifferential $\partial\psi_{ah}$ is defined by \eqref{eq:initial-subdifferential} relative to $\mcl Q_2$ and its $L^2$ pairing.

\begin{proposition}[Critical-point representation for the constrained free energy]
\label{prop:block-critical}
For every $a\in\mcl A^\circ$, $t\ge0$, and $q\in\mcl Q_\infty$, there exists $p\in\mcl Q_{\infty,\preccurlyeq I_D}$ such that
\begin{equation}
\label{eq:block-critical-conclusion}
p\in\partial\psi_{ah}\big(q+t\nabla\xi_a(p)\big)\quad\text{and}\quad \lim_{N\to\infty}\bar F_{N,a}(t,q)=\sP_{a,t,q}(p).
\end{equation}
\end{proposition}

\begin{proof}
\smallskip
\noindent\emph{Step 1: Compactness of $p_M$ and $\pi_M$.} First assume that, for some $c>0$, we have $q(0)\succcurlyeq cI_D$ and $q(v)-q(u)\succcurlyeq c(v-u)I_D$ for $0\le u\le v<1$. For each $M\in\N$ with $M\ge\max(D,3)$, choose $p_M$ from Lemma~\ref{lem:block-product-critical} and put $\pi_M:=q+t\nabla\xi_a(p_M)$. We will pass both identities in \eqref{eq:block-product-critical} to the limit.

By \cite[Lemma~3.4]{chen2025free}, $L^\infty$-bounded subsets of $\mcl Q_\infty$ are precompact in $L^2$. Since $|p_M|_{L^\infty}\le\sqrt D$, there are $M_j\to\infty$ and $p\in\mcl Q_\infty$ such that $p_{M_j}$ converges to $p$ in $L^2$ and almost everywhere. The inequalities $0\preccurlyeq p_M\preccurlyeq I_D$ pass to the limit almost everywhere and then, by right continuity, at every $s<1$. Thus $p\in\mcl Q_{\infty,\preccurlyeq I_D}$.

The Hessian of $\xi_a$ is bounded on $\{r\in\S^D:0\preccurlyeq r\preccurlyeq I_D\}$ by \eqref{eq:mixture-summability}. Hence, for a constant $C_a<\infty$, we have $|\nabla\xi_a(p_M)-\nabla\xi_a(p)|_{L^2}\le C_a|p_M-p|_{L^2}$. Define $\pi:=q+t\nabla\xi_a(p)$. We conclude that
\begin{equation}
\label{eq:block-critical-compactness}
\lim_{j\to\infty}|p_{M_j}-p|_{L^2}=0\quad\text{and}\quad \lim_{j\to\infty}|\pi_{M_j}-\pi|_{L^2}=0.
\end{equation}
By Lemma~\ref{lem:mixture-monotonicity}, the paths $\nabla\xi_a(p_M)$ and $\nabla\xi_a(p)$ are nonnegative and increasing. Therefore $\pi_M$ and $\pi$ satisfy the same inequalities at $0$ and on increments as $q$. Put $C_0:=|q|_{L^\infty}+t\sup_{0\preccurlyeq r\preccurlyeq I_D}|\nabla\xi_a(r)|$. They also satisfy
\begin{equation}
\label{eq:block-critical-path-bound}
\max\{|\pi_M|_{L^\infty},|\pi|_{L^\infty}\}\le C_0<\infty.
\end{equation}
Decrease $c$ if necessary so that $c\le1$ and $C_0\le c^{-1}$. The lower bounds at $0$ and on increments remain valid, and Lemma~\ref{lem:block-semiconcavity} now applies with this common $c$.

\smallskip
\noindent\emph{Step 2: Passage to the value identity.} The positive-semidefinite paths in \eqref{eq:block-critical-path-bound} belong to $\mcl Q_{\infty,\preccurlyeq C_0I_D}$. By Corollary~\ref{cor:initial-uniform} and \eqref{eq:psi-Lipschitz}, we have
\begin{equation*}
|\psi_{M_j,ah}(\pi_{M_j})-\psi_{ah}(\pi)|\le\sup_{r\in\mcl Q_{\infty,\preccurlyeq C_0I_D}}|\psi_{M_j,ah}(r)-\psi_{ah}(r)|+\sqrt D|\pi_{M_j}-\pi|_{L^1}.
\end{equation*}
Both terms on the right have limit zero by \eqref{eq:initial-uniform} and \eqref{eq:block-critical-compactness}. The function $\theta_a$ is Lipschitz on $\{0\preccurlyeq r\preccurlyeq I_D\}$, so, for a constant $C_a<\infty$, we also have
\begin{equation*}
\left|\int_0^1\theta_a(p_{M_j}(s))\,\mathrm d s-\int_0^1\theta_a(p(s))\,\mathrm d s\right|\le C_a|p_{M_j}-p|_{L^1}.
\end{equation*}
Consequently, every term of the value identity in \eqref{eq:block-product-critical} has a limit. Equations \eqref{eq:block-limit-error} and \eqref{eq:block-critical-compactness} give
\begin{equation}
\label{eq:block-limiting-value}
\bar F_{\infty,a}(t,q)=\psi_{ah}(\pi)-t\int_0^1\theta_a(p(s))\,\mathrm d s.
\end{equation}

\smallskip
\noindent\emph{Step 3: Derivatives of $\psi_{ah}$ along smooth paths.} Fix $\kappa\in C^1([0,1],\S^D)$, and put $L:=|\kappa(0)|+|\kappa'|_{L^\infty}$. For $|s|L\le c/2$, the bounds $|\kappa(0)|\le L$ and $|\kappa(v)-\kappa(u)|\le L(v-u)$ show that $\pi_{M_j}+s\kappa$ and $\pi+s\kappa$ have value at least $(c/2)I_D$ at $0$ and increments at least $(c/2)(v-u)I_D$. Thus they belong to $\mcl Q_\infty$. Their norms are at most $C_0+c/2$. The convergence of $\psi_{M_j,ah}(\pi_{M_j}+s\kappa)$ to $\psi_{ah}(\pi+s\kappa)$ follows from the same inequality as in Step~2, now with this larger bound. Moreover,
\begin{equation*}
|\langle p_{M_j}-p,\kappa\rangle_{L^2}|\le|p_{M_j}-p|_{L^2}|\kappa|_{L^2}.
\end{equation*}
Apply \eqref{eq:block-semiconcavity} with $M=M_j$, $x=ah$, and $\pi=\pi_{M_j}$. Its remainder is $c^{-1}L^2s^2$, independent of $j$, and \eqref{eq:block-product-critical} identifies the derivative with $p_{M_j}$. Taking $j\to\infty$ gives, for both signs of $s$,
\begin{equation}
\label{eq:block-semiconcavity-limit}
\psi_{ah}(\pi+s\kappa)\le\psi_{ah}(\pi)+s\langle p,\kappa\rangle_{L^2}+C_\kappa s^2.
\end{equation}
Here $C_\kappa:=c^{-1}L^2$ is independent of $s$. By Proposition~\ref{prop:stiefel-convexity}, the function $s\mapsto\psi_{ah}(\pi+s\kappa)$ is convex. For $s>0$ with $sL\le c/2$, its left and right secant slopes satisfy
\begin{equation*}
\frac{\psi_{ah}(\pi)-\psi_{ah}(\pi-s\kappa)}s\le\frac{\psi_{ah}(\pi+s\kappa)-\psi_{ah}(\pi)}s.
\end{equation*}
Use \eqref{eq:block-semiconcavity-limit} at $-s$ to bound the left side below, and at $s$ to bound the right side above. This yields
\begin{equation}
\label{eq:block-directional-limit}
\langle p,\kappa\rangle_{L^2}-C_\kappa s\le\frac{\psi_{ah}(\pi+s\kappa)-\psi_{ah}(\pi)}s\le\langle p,\kappa\rangle_{L^2}+C_\kappa s.
\end{equation}
In particular, the limit of this difference quotient as $s\downarrow0$ is $\langle p,\kappa\rangle_{L^2}$.

\smallskip
\noindent\emph{Step 4: The inequality for every $r\in\mcl Q_2$.} Fix $r\in\mcl Q_2$. By density, choose $\kappa_k\in C^1([0,1],\S^D)$ such that $\lim_{k\to\infty}|\kappa_k-(r-\pi)|_{L^1}=0$. For example, extend $r-\pi$ by zero and convolve its matrix entries with a smooth approximate identity. These approximations need not be increasing. The inequalities satisfied by $\pi$ in Step~1 ensure that $\pi+s\kappa_k\in\mcl Q_\infty$ whenever $|s|(|\kappa_k(0)|+|\kappa_k'|_{L^\infty})\le c/2$. On the other hand, $\pi+s(r-\pi)\in\mcl Q_2$ for $s\in[0,1]$ by convexity of the path cone.

For each fixed $k$, choose $s_k\in(0,1]$ such that both paths belong to $\mcl Q_2$ for $0\le s\le s_k$. Equation \eqref{eq:psi-Lipschitz} gives
\begin{equation}
\label{eq:block-subgradient-approximation}
\left|\frac{\psi_{ah}(\pi+s(r-\pi))-\psi_{ah}(\pi)}s-\frac{\psi_{ah}(\pi+s\kappa_k)-\psi_{ah}(\pi)}s\right|\le\sqrt D|r-\pi-\kappa_k|_{L^1},\qquad 0<s\le s_k.
\end{equation}
Convexity of $\psi_{ah}$ gives $\psi_{ah}(\pi+s(r-\pi))\le(1-s)\psi_{ah}(\pi)+s\psi_{ah}(r)$ for $s\in[0,1]$. Therefore, by \eqref{eq:block-subgradient-approximation}, we obtain
\begin{equation*}
\psi_{ah}(r)-\psi_{ah}(\pi)\ge\frac{\psi_{ah}(\pi+s\kappa_k)-\psi_{ah}(\pi)}s-\sqrt D|r-\pi-\kappa_k|_{L^1},\qquad 0<s\le s_k.
\end{equation*}
First let $s\downarrow0$ at fixed $k$, using \eqref{eq:block-directional-limit} with $\kappa_k$. Then let $k\to\infty$. The pairings converge since $|\langle p,\kappa_k-(r-\pi)\rangle_{L^2}|\le\sqrt D|\kappa_k-(r-\pi)|_{L^1}$. We conclude that
\begin{equation}
\label{eq:block-limiting-subgradient}
\psi_{ah}(r)\ge\psi_{ah}(\pi)+\langle p,r-\pi\rangle_{L^2},\qquad r\in\mcl Q_2.
\end{equation}
By \eqref{eq:initial-subdifferential}, this is $p\in\partial\psi_{ah}(\pi)$. Together with \eqref{eq:block-limiting-value} and $\pi=q+t\nabla\xi_a(p)$, it proves \eqref{eq:block-critical-conclusion} under the assumption of Step~1.

\smallskip
\noindent\emph{Step 5: Passage from $q+\varepsilon(1+\cdot)I_D$ to $q$.} For arbitrary $q\in\mcl Q_\infty$ and $\varepsilon\in(0,1]$, define $q_\varepsilon:[0,1)\to\S_+^D$ by $q_\varepsilon(s):=q(s)+\varepsilon(1+s)I_D$. Its value at $0$ is at least $\varepsilon I_D$, and its increment from $u$ to $v$ is at least $\varepsilon(v-u)I_D$. Steps~1--4 therefore supply $p_\varepsilon\in\mcl Q_{\infty,\preccurlyeq I_D}$ and $\pi_\varepsilon:=q_\varepsilon+t\nabla\xi_a(p_\varepsilon)$ such that
\begin{equation}
\label{eq:block-epsilon-critical}
p_\varepsilon\in\partial\psi_{ah}(\pi_\varepsilon)\quad\text{and}\quad \bar F_{\infty,a}(t,q_\varepsilon)=\psi_{ah}(\pi_\varepsilon)-t\int_0^1\theta_a(p_\varepsilon(s))\,\mathrm d s.
\end{equation}
Apply the compactness argument of Step~1 to a sequence $\varepsilon\downarrow0$. There are $\varepsilon_j\downarrow0$ and $p\in\mcl Q_{\infty,\preccurlyeq I_D}$ such that $\lim_{j\to\infty}|p_{\varepsilon_j}-p|_{L^2}=0$. Put $\pi:=q+t\nabla\xi_a(p)$. Since $|q_\varepsilon-q|_{L^\infty}\le2\sqrt D\varepsilon$ and $\nabla\xi_a$ and $\theta_a$ are Lipschitz on $\{0\preccurlyeq r\preccurlyeq I_D\}$, we have
\begin{equation}
\label{eq:block-epsilon-limits}
\lim_{j\to\infty}|\pi_{\varepsilon_j}-\pi|_{L^2}=0\quad\text{and}\quad \lim_{j\to\infty}\int_0^1\theta_a(p_{\varepsilon_j}(s))\,\mathrm d s=\int_0^1\theta_a(p(s))\,\mathrm d s.
\end{equation}
For fixed $r\in\mcl Q_2$, the inclusion in \eqref{eq:block-epsilon-critical} means $\psi_{ah}(r)\ge\psi_{ah}(\pi_{\varepsilon_j})+\langle p_{\varepsilon_j},r-\pi_{\varepsilon_j}\rangle_{L^2}$. The function values converge by \eqref{eq:psi-Lipschitz}. For the pairings, we have
\begin{equation*}
\left|\langle p_{\varepsilon_j},r-\pi_{\varepsilon_j}\rangle_{L^2}-\langle p,r-\pi\rangle_{L^2}\right|\le |p_{\varepsilon_j}-p|_{L^2}|r-\pi|_{L^2}+\sqrt D|\pi_{\varepsilon_j}-\pi|_{L^2}.
\end{equation*}
Taking $j\to\infty$ yields \eqref{eq:block-limiting-subgradient}, and hence $p\in\partial\psi_{ah}(\pi)$. Finally, taking $N\to\infty$ in \eqref{eq:fixed-enriched-q-Lipschitz} gives
\begin{equation*}
|\bar F_{\infty,a}(t,q_{\varepsilon_j})-\bar F_{\infty,a}(t,q)|\le\sqrt D|q_{\varepsilon_j}-q|_{L^1}\le2D\varepsilon_j.
\end{equation*}
We may therefore pass to the limit in the value identity in \eqref{eq:block-epsilon-critical}, using \eqref{eq:block-epsilon-limits} and \eqref{eq:psi-Lipschitz}. This gives $\bar F_{\infty,a}(t,q)=\sP_{a,t,q}(p)$ and completes the proof of \eqref{eq:block-critical-conclusion}.
\end{proof}

\section{Identifying the constrained free energy}
\label{s:fixed-a-free-energy}

We identify the constrained limit with the Hopf value $f_a$ defined in \eqref{eq:Hopf-fixed-a-intro}. Proposition~\ref{prop:HJ-lower} gives one inequality by \cite[Theorem~3.4]{mourrat2023upper} and \cite[Theorem~4.7(3)]{chenXia2025HJ}. The critical relation in Proposition~\ref{prop:block-critical} gives the reverse inequality. Corollary~\ref{cor:fixed-a-limit} gives the limits of both $\bar F_{N,a}$ and
\begin{equation}
\label{eq:fixed-a-free-energy-intro}
\mathsf F_{N,a}:=-\frac1N\E\log\int_{\mathbb U_N}\exp\Big(H_N(au)+ah\cdot\sum_{i=1}^Nu_i\Big)\nu_N(\mathrm d u).
\end{equation}

For $a\in\mcl A^\circ$, recall $f_a(t,q)$ from~\eqref{eq:Hopf-fixed-a-intro} and $\mcl J_{a,t,q}$ from~\eqref{eq:J-fixed-a-intro}. Since $\xi_a\ge0$ on $\S_+^D$, choosing $p\in\partial\psi_{ah}(q)\cap\mcl Q_{\infty,\preccurlyeq I_D}$ by Lemma~\ref{lem:initial-subgradient} gives the first inequality below, and testing $q'=q$ for every $p$ gives the second, in both cases by \eqref{eq:Hopf-fixed-a-intro}: 
\begin{equation*}
\psi_{ah}(q)\le f_a(t,q)\le\psi_{ah}(q)+t\max_{0\preccurlyeq r\preccurlyeq I_D}\xi_a(r)<\infty,\qquad t\ge0,\quad q\in\mcl Q_2.
\end{equation*}
In particular, $f_a(0,q)=\psi_{ah}(q)$. Since the dependence of $\mcl J_{a,t,q}(q',p)$ on $q$ is through $\langle q,p\rangle_{L^2}$ and $|p|_{L^\infty}\le\sqrt D$, taking the infimum over $q'$ and then the supremum over $p$ gives
\begin{equation}
\label{eq:Hopf-fixed-path-Lipschitz}
|f_a(t,q)-f_a(t,\widetilde q)|\le\sqrt D\,|q-\widetilde q|_{L^1},\qquad t\ge0,\quad q,\widetilde q\in\mcl Q_2.
\end{equation}

Proposition~\ref{prop:block-limit} gives convergence of $\bar F_{N,a}(t,q)$ for bounded $q$. The uniform $L^1$ estimate \eqref{eq:fixed-enriched-q-Lipschitz} extends this convergence to every $q\in\mcl Q_2$ by bounded-path approximation. We therefore write a limit in the following bound.

\begin{proposition}[Hamilton--Jacobi lower bound]
\label{prop:HJ-lower}
For every $a\in\mcl A^\circ$, $t\ge0$, and $q\in\mcl Q_2$, we have
\begin{equation}
\label{eq:HJ-liminf}
f_a(t,q)\le\lim_{N\to\infty}\bar F_{N,a}(t,q).
\end{equation}
\end{proposition}

\begin{proof}
Propositions~\ref{prop:stiefel-initial} and~\ref{prop:stiefel-convexity}, together with Lemma~\ref{lem:initial-subgradient}, give the limit initial condition $\psi_{ah}$, its convexity, and its $\mcl Q_2^*$-monotonicity. Since \eqref{eq:mixture-summability} gives absolute convergence of $\xi$ everywhere, the results of \cite{mourrat2023upper,chenXia2025HJ} apply to our normalization by a fixed scaling. Thus \cite[Theorem~3.4]{mourrat2023upper}, in the form of \cite[Theorem~4.14]{chenXia2025HJ}, and the Hopf formula in \cite[Theorem~4.7(3)]{chenXia2025HJ} give
\begin{equation*}
\lim_{N\to\infty}\bar F_{N,a}(t,q)\ge\sup_{p\in\mcl Q_\infty}\inf_{q'\in\mcl Q_\infty}\left\{\psi_{ah}(q')+\langle q-q',p\rangle_{L^2}+t\int_0^1\xi_a(p(s))\,\mathrm d s\right\}.
\end{equation*}
For fixed $p\in\mcl Q_\infty$, the expression in braces is $L^1$-continuous in $q'$ by \eqref{eq:psi-Lipschitz}, so the inner infimum is unchanged if $q'\in\mcl Q_2$. Restricting the supremum to $p\in\mcl Q_{\infty,\preccurlyeq I_D}$ then gives \eqref{eq:HJ-liminf} by \eqref{eq:Hopf-fixed-a-intro}.
\end{proof}

\begin{corollary}[Limits of the constrained free energies]
\label{cor:fixed-a-limit}
For every $a\in\mcl A^\circ$, $t\ge0$, and $q\in\mcl Q_2$, we have 
\begin{gather}
\label{eq:fixed-a-negative-limit}
\lim_{N\to\infty}\bar F_{N,a}(t,q)=f_a(t,q)\qquad\text{and}\qquad
\lim_{N\to\infty}\mathsf F_{N,a}=f_a\Ll(\tfrac12,0\Rr)-\tfrac12\xi(a^2).
\end{gather}
\end{corollary}

\begin{proof}
First suppose that $q\in\mcl Q_\infty$. Choose $p^+\in\mcl Q_{\infty,\preccurlyeq I_D}$ as in Proposition~\ref{prop:block-critical}, and set
\begin{equation}
\label{eq:qplus}
q^+:=q+t\nabla\xi_a(p^+).
\end{equation}
Then $q^+\in\mcl Q_\infty$ by Lemma~\ref{lem:mixture-monotonicity}, and $p^+\in\partial\psi_{ah}(q^+)$ means that
\begin{equation}
\label{eq:subgradient-qplus}
\psi_{ah}(q')\ge\psi_{ah}(q^+)+\int_0^1p^+(s)\cdot\Ll(q'(s)-q^+(s)\Rr)\,\mathrm d s,\qquad q'\in\mcl Q_2.
\end{equation}
Recall $\theta_a$ from \eqref{eq:theta-a-definition}. For every $q'\in\mcl Q_2$, we obtain
\begin{align*}
\mcl J_{a,t,q}(q',p^+)&\stackrel{\eqref{eq:subgradient-qplus}}{\ge}\psi_{ah}(q^+)+\langle q-q^+,p^+\rangle_{L^2}+t\int_0^1\xi_a(p^+(s))\,\mathrm d s\\
&\stackrel{\eqref{eq:qplus}}{=}\psi_{ah}(q^+)-t\int_0^1\theta_a(p^+(s))\,\mathrm d s=\sP_{a,t,q}(p^+).
\end{align*}
Equality holds at $q'=q^+$. Hence $q^+$ attains the inner infimum at $p^+$, and \eqref{eq:Hopf-fixed-a-intro} gives
\begin{equation}
\label{eq:P-less-f}
\inf_{q'\in\mcl Q_2}\mcl J_{a,t,q}(q',p^+)=\sP_{a,t,q}(p^+)\le f_a(t,q).
\end{equation}
Combining this with Propositions~\ref{prop:HJ-lower} and~\ref{prop:block-critical}, we obtain
\begin{equation*}
f_a(t,q)\stackrel{\eqref{eq:HJ-liminf}}{\le}\lim_{N\to\infty}\bar F_{N,a}(t,q)\stackrel{\eqref{eq:block-critical-conclusion}}{=}\sP_{a,t,q}(p^+)\stackrel{\eqref{eq:P-less-f}}{\le}f_a(t,q).
\end{equation*}
This proves the first identity in \eqref{eq:fixed-a-negative-limit} for bounded $q$, including $t=0$.

For $q\in\mcl Q_2$ and $n\ge2$, set $q_n(s):=q(\min\{s,1-n^{-1}\})$. Then $q_n\in\mcl Q_\infty$ and $\lim_{n\to\infty}|q_n-q|_{L^1}=0$. The convergence already proved at $q_n$, together with \eqref{eq:fixed-enriched-q-Lipschitz} and \eqref{eq:Hopf-fixed-path-Lipschitz}, gives
\begin{equation*}
\limsup_{N\to\infty}|\bar F_{N,a}(t,q)-f_a(t,q)|\le2\sqrt D\,|q_n-q|_{L^1}.
\end{equation*}
Sending $n\to\infty$ proves the first identity in \eqref{eq:fixed-a-negative-limit} for every $q\in\mcl Q_2$.

Finally, $\xi_a(I_D)=\xi(a^2)$ by \eqref{eq:fixed-a-xi-intro}. Setting $t=1/2$ and $q=0$ in \eqref{eq:enriched-fixed-a} therefore gives
\begin{equation}
\label{eq:xia-I}
\bar F_{N,a}\Ll(\tfrac12,0\Rr)=\mathsf F_{N,a}+\tfrac12\xi(a^2).
\end{equation}
Taking $N\to\infty$ and using the first identity in \eqref{eq:fixed-a-negative-limit} proves the second.
\end{proof}

\section{Optimization over the constraint parameter}
\label{s:main-proof}

Section~\ref{s:fixed-a-free-energy} identifies the limit for each $a\in\mcl A^\circ$. We now integrate over $a$ to obtain the unrestricted free energies. Lemma~\ref{lem:disintegration} gives the exact integral and its determinant weight. We then prove convergence uniform on sets where $\det(a^2)\ge\kappa$, compute the exponential rate of the restricted integrals, and bound the contribution from $\det(a^2)<\kappa$. These estimates give Propositions~\ref{prop:global-enriched} and~\ref{prop:outer-a-free-energy}, including attainment of the extrema.

\subsection{Disintegration}
\label{s:outer-a} 

Throughout this section, assume $N>D$. Recall $S_N^D$, $P_N=\mu_N^{\otimes D}$, and $\mcl A$ from Section~\ref{s:introduction}. The diagonal entries of $N^{-1}\sigma\sigma^{\intercal}$ equal $1$ for $\sigma\in S_N^D$. Thus its possible values form $\mcl B:=\{b\in\S_+^D:b_{dd}=1\text{ for }d\in\{1,\dots,D\}\}=\{a^2:a\in\mcl A\}$. The set $\mcl B^\circ$ from the proof outline is $\mcl B\cap\S_{++}^D$, and $\d b_{\mathrm{off}}$ is Lebesgue measure in the $D(D-1)/2$ off-diagonal entries of $b$. Thus $b\in\mcl B^\circ$ corresponds bijectively to $a=b^{1/2}\in\mcl A^\circ$. We retain $b$ as the integration variable in this section so that the Wishart density does not acquire a change-of-variables Jacobian. When $D=1$, the measure is the unit measure on $\{I_1\}=\{1\}$ and $\mcl A=\mcl A^\circ=\{I_1\}$. The boundary arguments below are needed only for $D\ge2$. 
Let $\Gamma(x):=\int_0^\infty t^{x-1}e^{-t}\,\mathrm d t$, $x>0$, be the gamma function. Define the normalizing constant by
\begin{equation}
\label{eq:disintegration-constant}
c_{N,D}:=\frac{\Gamma(N/2)^D}{\pi^{D(D-1)/4}\prod_{j=1}^D\Gamma((N-j+1)/2)}.
\end{equation}

\begin{lemma}
\label{lem:matrix-density}
For $N>D$, the law of $\frac1N\sigma\sigma^{\intercal}$ under $P_N$ has density
\begin{equation}
\label{eq:matrix-density}
c_{N,D}(\det b)^{(N-D-1)/2}\one_{\{b\in\mcl B^\circ\}}\,\d b_{\mathrm{off}}.
\end{equation}
Moreover, we have
\begin{equation}
\label{eq:normalization-subexp}
\lim_{N\to\infty}\frac1N\log c_{N,D}=0.
\end{equation}
\end{lemma}

\begin{proof}
When $D=1$, $N^{-1}\sigma\sigma^{\intercal}=I_1$ and $c_{N,1}=1$, so both conclusions hold. Assume $D\ge2$.

\smallskip
\noindent\emph{Step 1: Representing $N^{-1}\sigma\sigma^{\intercal}$ by Gaussian vectors.} Let $g\in\R^{D\times N}$ have independent standard Gaussian entries, and write $g_c\in\R^N$ for its $c$th row. Rotational invariance of $g_c$ implies that $\sigma_c:=\sqrt N\,g_c/|g_c|$ is uniform on $S_N$. Since each $\sigma_c$ depends only on $g_c$, these rows remain independent. Hence the matrix $\sigma$ with rows $\sigma_c$ has law $P_N$. Put
\begin{equation*}
w:=gg^{\intercal},\qquad s_c:=w_{cc}=|g_c|^2\quad(c\in\{1,\dots,D\}),\qquad s:=\operatorname{diag}(s_1,\ldots,s_D).
\end{equation*}
The first $D$ columns of $g$ form an invertible matrix almost surely: their determinant is a nonzero polynomial, whose zero set has Lebesgue measure zero, and their joint law has a density. Thus $w\in\S_{++}^D$ almost surely. The row normalization gives
\begin{equation*}
\sigma=\sqrt N\,s^{-1/2}g,\qquad b:=\frac1N\sigma\sigma^{\intercal}=s^{-1/2}ws^{-1/2}\in\mcl B^\circ,\qquad w=s^{1/2}bs^{1/2}.
\end{equation*}
We will compute the law of $b$ from the joint law of $(s_1,\ldots,s_D,b)$.

\smallskip
\noindent\emph{Step 2: The change of variables and its Jacobian.} The matrix $w=gg^{\intercal}$ has the Wishart distribution with $N$ degrees of freedom and covariance $I_D$. By \cite[Theorem~3.2.1]{muirhead1982aspects}, its density with respect to $\prod_{c\le d}\mathrm d w_{cd}$ is
\begin{equation*}
\frac{e^{-\tr w/2}(\det w)^{(N-D-1)/2}}{2^{ND/2}\pi^{D(D-1)/4}\prod_{j=1}^D\Gamma((N-j+1)/2)}\one_{\{w\in\S_{++}^D\}}.
\end{equation*}
Consider the map
\begin{equation*}
\Phi:(0,\infty)^D\times\mcl B^\circ\to\S_{++}^D,\qquad \Phi(s_1,\ldots,s_D,b):=s^{1/2}bs^{1/2}.
\end{equation*}
Its inverse is $s_c=w_{cc}$ and $b_{cd}=w_{cd}/\sqrt{w_{cc}w_{dd}}$, so $\Phi$ is a smooth bijection with smooth inverse. To compute its Jacobian, list the output coordinates as $(w_{11},\ldots,w_{DD},(w_{cd})_{c<d})$ and the input coordinates as $(s_1,\ldots,s_D,(b_{cd})_{c<d})$, using the same lexicographic order for the pairs $(c,d)$ in both lists. Since $b_{cc}=1$, the coordinate identities are
\begin{equation*}
w_{cc}=s_c\quad(c\in\{1,\dots,D\}),\qquad w_{cd}=\sqrt{s_cs_d}\,b_{cd}\quad(c<d).
\end{equation*}
The derivative matrix of $(w_{11},\ldots,w_{DD})$ with respect to $(s_1,\ldots,s_D)$ is $I_D$, and all derivatives of these outputs with respect to the off-diagonal entries of $b$ vanish. For $c<d$, the output $w_{cd}$ depends on the off-diagonal entries of $b$ only through $b_{cd}$, with derivative $\sqrt{s_cs_d}$. Thus the Jacobian matrix is
\begin{equation*}
\nabla\Phi=\begin{pmatrix}
I_D&0\\
\big(\partial w_{cd}/\partial s_j\big)_{c<d,\,j\in\{1,\dots,D\}}&\operatorname{diag}\big((\sqrt{s_cs_d})_{c<d}\big)
\end{pmatrix}.
\end{equation*}
The lower-left block does not enter the determinant of this block triangular matrix. Each $s_c$ appears in exactly $D-1$ of the factors indexed by $c<d$, so
\begin{equation}
\label{eq:disintegration-Jacobian}
|\det \nabla\Phi|=\prod_{c<d}\sqrt{s_cs_d}=\prod_{c=1}^Ds_c^{(D-1)/2}.
\end{equation}

\smallskip
\noindent\emph{Step 3: Integrating out $s_1,\ldots,s_D$.} Substitute $w=\Phi(s_1,\ldots,s_D,b)$ into the Wishart density and multiply by \eqref{eq:disintegration-Jacobian}. The identities
\begin{equation*}
\det w=(\det b)\prod_{c=1}^Ds_c,\qquad \tr w=\sum_{c=1}^Ds_c,\qquad \frac{N-D-1}{2}+\frac{D-1}{2}=\frac N2-1
\end{equation*}
show that the joint density of $(s_1,\ldots,s_D,b)$ with respect to $\prod_{c=1}^D\mathrm d s_c\,\mathrm d b_{\mathrm{off}}$ is
\begin{equation*}
\frac{(\det b)^{(N-D-1)/2}\prod_{c=1}^D s_c^{N/2-1}e^{-s_c/2}}{2^{ND/2}\pi^{D(D-1)/4}\prod_{j=1}^D\Gamma((N-j+1)/2)},\qquad b\in\mcl B^\circ,\quad s_c>0.
\end{equation*}
For each $c\in\{1,\dots,D\}$, the substitution $s_c=2t$ gives
\begin{equation*}
\int_0^\infty s_c^{N/2-1}e^{-s_c/2}\,\mathrm d s_c=2^{N/2}\int_0^\infty t^{N/2-1}e^{-t}\,\mathrm d t=2^{N/2}\Gamma(N/2).
\end{equation*}
Tonelli's theorem permits these $D$ integrations. Their product contributes $2^{ND/2}\Gamma(N/2)^D$, which cancels the factor $2^{ND/2}$ in the denominator and leaves $c_{N,D}$ from \eqref{eq:disintegration-constant}. This proves \eqref{eq:matrix-density}.

\smallskip
\noindent\emph{Step 4: The asymptotics of $c_{N,D}$.} For fixed $\beta\ge0$, Stirling's formula and $\log(x-\beta)=\log x-\beta/x+O(x^{-2})$ give, as $x\to\infty$,
\begin{align*}
\log\Gamma(x)-\log\Gamma(x-\beta)
&=(x-\tfrac12)\log x-(x-\beta-\tfrac12)\log(x-\beta)-\beta+O(x^{-1})\\
&=\beta\log x+O(x^{-1}).
\end{align*}
Exponentiating and setting $x=N/2$ and $\beta=(j-1)/2$ yields
\begin{equation*}
\frac{\Gamma(N/2)}{\Gamma((N-j+1)/2)}=(N/2)^{(j-1)/2}\big(1+O_D(N^{-1})\big),\qquad j\in\{1,\dots,D\}.
\end{equation*}
Since $D$ is fixed, multiplying these $D$ identities gives
\begin{align*}
c_{N,D}\stackrel{\eqref{eq:disintegration-constant}}{=}\pi^{-D(D-1)/4}\prod_{j=1}^D\frac{\Gamma(N/2)}{\Gamma((N-j+1)/2)}&=\pi^{-D(D-1)/4}(N/2)^{\sum_{j=1}^D(j-1)/2}\big(1+O_D(N^{-1})\big)\\
&=(N/(2\pi))^{D(D-1)/4}\big(1+O_D(N^{-1})\big).
\end{align*}
Thus
\begin{equation*}
\frac1N\log c_{N,D}=\frac{D(D-1)}{4N}\log\frac{N}{2\pi}+O_D(N^{-2}),
\end{equation*}
whose right-hand side tends to zero as $N\to\infty$. This proves \eqref{eq:normalization-subexp}.
\end{proof}

Lemma~\ref{lem:matrix-density} gives the law of $b=N^{-1}\sigma\sigma^{\intercal}$. To integrate a function of $\sigma$, we also need its average over configurations with the same $b$. Right-orthogonal invariance of $P_N$ identifies this average with integration over $\sigma=b^{1/2}u$, $u\in\mathbb U_N$, against $\nu_N$. The next lemma combines this identity with the density \eqref{eq:matrix-density}.

\begin{lemma}
\label{lem:disintegration}
For every integer $N\ge D+1$ and every nonnegative measurable function $G$ on $S_N^D$, we have 
\begin{equation}
\label{eq:measure-disintegration}
\int_{S_N^D}G(\sigma)P_N(\d\sigma) =c_{N,D}\iint_{\mathbb U_N\times\mcl B^\circ} G(b^{1/2}u)\nu_N(\d u)(\det b)^{(N-D-1)/2}\,\d b_{\mathrm{off}}.
\end{equation}
\end{lemma}

\begin{proof}
\smallskip
\noindent\emph{Step 1: Averaging over orthogonal rotations.} Let $O$ be Haar distributed on $\operatorname{Orth}(N)$, independently of $\sigma$ under $P_N$, and let $\E_O$ denote expectation with respect to $O$. For every deterministic $o\in\operatorname{Orth}(N)$, right multiplication by $o$ preserves the uniform law of each row of $\sigma$ and their independence. Hence $\sigma o$ has law $P_N$. Integrating this identity over $o$ with Haar probability and applying Tonelli's theorem gives
\begin{equation}
\label{eq:disintegration-Haar-average}
\int_{S_N^D}G(\sigma)P_N(\mathrm d\sigma)=\int_{S_N^D}\E_O G(\sigma O)P_N(\mathrm d\sigma).
\end{equation}

\smallskip
\noindent\emph{Step 2: The law of $uO$ for fixed $u\in\mathbb U_N$.} Fix any deterministic $u\in\mathbb U_N$. Then $uO\in\mathbb U_N$, since $(uO)(uO)^{\intercal}=uu^{\intercal}=NI_D$. For every deterministic $o\in\operatorname{Orth}(N)$, right invariance of Haar probability gives $Oo\stackrel{\mathrm D}=O$. Therefore, for every bounded measurable $\varphi:\mathbb U_N\to\R$,
\begin{equation*}
\E_O\varphi((uO)o)=\E_O\varphi(u(Oo))=\E_O\varphi(uO).
\end{equation*}
Thus the probability law of $uO$ is invariant under right multiplication by every $o\in\operatorname{Orth}(N)$. The uniqueness assertion in Lemma~\ref{lem:stiefel-invariance} applies to this law and identifies it with $\nu_N$. In particular, for every nonnegative measurable $\varphi$,
\begin{equation}
\label{eq:fixed-rotation-law}
\E_O\varphi(uO)=\int_{\mathbb U_N}\varphi(v)\nu_N(\mathrm d v),\qquad u\in\mathbb U_N.
\end{equation}
The extension from bounded to nonnegative $\varphi$ follows by applying the bounded identity to $\min\{\varphi,m\}$ and then letting $m\to\infty$. In \eqref{eq:fixed-rotation-law}, only $O$ is random. No distributional assumption on $u$ is needed.

\smallskip
\noindent\emph{Step 3: Integrating over $b=N^{-1}\sigma\sigma^{\intercal}$.} By Lemma~\ref{lem:matrix-density}, $b:=N^{-1}\sigma\sigma^{\intercal}\in\mcl B^\circ$ for $P_N$-almost every $\sigma$. For each such $\sigma$, put $u:=b^{-1/2}\sigma$. Then
\begin{equation*}
uu^{\intercal}=b^{-1/2}\sigma\sigma^{\intercal}b^{-1/2}=NI_D,\qquad \sigma=b^{1/2}u.
\end{equation*}
For this fixed $\sigma$, apply \eqref{eq:fixed-rotation-law} with $\varphi(v)=G(b^{1/2}v)$. This gives
\begin{equation*}
\E_O G(\sigma O)=\E_O G(b^{1/2}uO)=\int_{\mathbb U_N}G(b^{1/2}v)\nu_N(\mathrm d v).
\end{equation*}
The last expression is a nonnegative measurable function of $b$: the map $(b,v)\mapsto b^{1/2}v$ is continuous, and integration against $\nu_N$ preserves measurability. Substituting this expression into \eqref{eq:disintegration-Haar-average} and using the density \eqref{eq:matrix-density}, we obtain
\begin{equation*}
\int_{S_N^D}G(\sigma)P_N(\mathrm d\sigma)=c_{N,D}\int_{\mcl B^\circ}\left(\int_{\mathbb U_N}G(b^{1/2}v)\nu_N(\mathrm d v)\right)(\det b)^{(N-D-1)/2}\,\mathrm d b_{\mathrm{off}}.
\end{equation*}
Tonelli's theorem now gives \eqref{eq:measure-disintegration}. 
\end{proof}

\subsection{Formula for \texorpdfstring{$\bar F_N$}{F̄\_N}}
\label{s:global-enriched}

We prove Proposition~\ref{prop:global-enriched} by integrating the constrained partition functions. For $\kappa\in(0,1)$, we first prove convergence uniform on $\det(a^2)\ge\kappa$ and then bound the integral over $\det(a^2)<\kappa$.

Recall $f_a$ from \eqref{eq:Hopf-fixed-a-intro}, $q_a(s)=aq(s)a$ from \eqref{eq:qa-intro}, and the continuous cascade $\mathfrak R$ and Gaussian fields $z_i^q$ from Section~\ref{s:enriched-properties}. For $t\ge0$ and $q\in\mcl Q_\infty$, define $H_N^{t,q}:S_N^D\times\supp\mathfrak R\to\R$ by
\begin{equation}
\label{eq:global-enriched-H}
H_N^{t,q}(\sigma,\alpha):=\sqrt{2t}H_N(\sigma)-Nt\xi\Ll(N^{-1}\sigma\sigma^{\intercal}\Rr)+\sqrt2\sum_{i=1}^Nz_i^q(\alpha)\cdot\sigma_i+h\cdot\sum_{i=1}^N\sigma_i-q(1)\cdot\sigma\sigma^{\intercal},
\end{equation}
and let
\begin{equation}
\label{eq:global-enriched-Z}
Z_N(t,q):=\iint_{S_N^D\times\supp\mathfrak R}\exp\big(H_N^{t,q}(\sigma,\alpha)\big)P_N(\mathrm d\sigma)\mathfrak R(\mathrm d\alpha).
\end{equation}
Equations \eqref{eq:global-enriched-intro} and \eqref{eq:global-enriched-Z} give $\bar F_N(t,q)=-N^{-1}\E\log Z_N(t,q)$. For finite-step $q$, the discrete cascade representation has the same averaged logarithm by \cite[Proposition~4.6]{chen2025free}, as recalled in Section~\ref{s:enriched-properties}.

We use the same interaction tensors, cascade, and fields $z_i^q$ for every $a\in\mcl A^\circ$. In $Z_{N,a}(t,q_a)$ from \eqref{eq:fixed-a-Z}, realize the Gaussian field with covariance path $q_a$ as $az_i^q$. Indeed, \eqref{eq:enriched-cascade-covariance} gives, for $a,a'\in\mcl A^\circ$ and $i,j\in\{1,\ldots,N\}$,
\begin{equation}
\label{eq:enrichment-fiber-cross-covariance}
\E\big[(az_i^q(\alpha))(a'z_j^q(\alpha'))^{\intercal}\mid\mathfrak R\big]=\one_{\{i=j\}}aq(\alpha\wedge\alpha')a'.
\end{equation}
With $a'=a$, this is $\one_{\{i=j\}}q_a(\alpha\wedge\alpha')$. Thus each $Z_{N,a}(t,q_a)$ has its prescribed law, and \eqref{eq:enrichment-fiber-cross-covariance} also specifies the joint law for different $a$. The interaction fields have cross-covariance
\begin{equation}
\label{eq:fiber-cross-covariance}
\E H_N(au)H_N(a'u')=N\xi\Ll(a(N^{-1}u(u')^{\intercal})a'\Rr),\qquad u,u'\in\mathbb U_N.
\end{equation}
Under this realization, the Hamiltonians satisfy
\begin{equation}
\label{eq:enriched-restriction-H}
H_N^{t,q}(au,\alpha)=H_{N,a}^{t,q_a}(u,\alpha),\qquad a\in\mcl A^\circ,\quad u\in\mathbb U_N,\quad\alpha\in\supp\mathfrak R.
\end{equation}
To verify \eqref{eq:enriched-restriction-H}, use $uu^{\intercal}=NI_D$, $\xi_a(I_D)=\xi(a^2)$, and the identities
\begin{equation*}
z_i^q(\alpha)\cdot au_i=(az_i^q(\alpha))\cdot u_i,\qquad h\cdot au_i=ah\cdot u_i,\qquad\text{and}\qquad q(1)\cdot a^2=\tr(aq(1)a)=q_a(1)\cdot I_D.
\end{equation*}
Substitution into \eqref{eq:global-enriched-H} gives exactly \eqref{eq:enriched-fixed-a-H} with path $q_a$ and fields $az_i^q$.

\begin{lemma}
\label{lem:global-enriched-disintegration}
For every $N\in\N$ with $N\ge D+1$, $t\ge0$, and $q\in\mcl Q_\infty$, we have almost surely
\begin{equation}
\label{eq:global-enriched-disintegration}
Z_N(t,q)=c_{N,D}\int_{\mcl B^\circ}Z_{N,b^{1/2}}(t,q_{b^{1/2}})(\det b)^{(N-D-1)/2}\,\mathrm d b_{\mathrm{off}},
\end{equation}
where all constrained partition functions use the realization in \eqref{eq:enriched-restriction-H}.
\end{lemma}

\begin{proof}
Fix a realization of the jointly measurable Gaussian fields and cascade. Define $G:S_N^D\to[0,\infty]$ by $G(\sigma):=\int\exp(H_N^{t,q}(\sigma,\alpha))\mathfrak R(\mathrm d\alpha)$. Then $G$ is nonnegative and measurable, and $Z_N(t,q)=\int G(\sigma)P_N(\mathrm d\sigma)$ by \eqref{eq:global-enriched-Z}. Lemma~\ref{lem:disintegration} gives
\begin{equation*}
Z_N(t,q)=c_{N,D}\int_{\mcl B^\circ}\left(\int_{\mathbb U_N}G(b^{1/2}u)\nu_N(\mathrm d u)\right)(\det b)^{(N-D-1)/2}\,\mathrm d b_{\mathrm{off}}.
\end{equation*}
Tonelli's theorem permits us to interchange the $u$ and cascade integrals. Hence, we have
\begin{equation*}
\int_{\mathbb U_N}G(b^{1/2}u)\nu_N(\mathrm d u)\stackrel{\eqref{eq:enriched-restriction-H}}{=}\iint_{\mathbb U_N\times\supp\mathfrak R}\exp\big(H_{N,b^{1/2}}^{t,q_{b^{1/2}}}(u,\alpha)\big)\nu_N(\mathrm d u)\mathfrak R(\mathrm d\alpha)\stackrel{\eqref{eq:fixed-a-Z}}{=}Z_{N,b^{1/2}}(t,q_{b^{1/2}}).
\end{equation*}
Substituting this identity proves \eqref{eq:global-enriched-disintegration} for the same realization of the fields and cascade.
\end{proof}

The next lemma bounds the dependence of $H_N(vu)$ on the factor matrix $v$. It will also be used in Section~\ref{s:positive-outer}. For $K\in\N$, let $H_N^{(K)}:\R^{D\times N}\to\R$ be the sum in \eqref{eq:mixed-Hamiltonian} restricted to $k\in\{1,\ldots,K\}$. For $v,b\in\R^{D\times D}$, the notation $(b\cdot\nabla_v)H_N(vu)$ denotes the derivative in $v$ in direction $b$.

\begin{lemma}
\label{lem:matrix-interaction-derivatives}
For every $R>0$, there exists $C_R<\infty$, independent of $N$, such that, for every $N\in\N$ with $N\ge D$, we have 
\begin{equation}
\label{eq:matrix-interaction-derivative-sup}
    \E\sup_{\substack{|v|\le R,\ |b|\le1\\u\in\mathbb U_N}}|(b\cdot\nabla_v)H_N(vu)|\le C_RN.
\end{equation}
As $K\to\infty$, the partial sums $H_N^{(K)}(vu)$ and their derivatives in $v$ converge to $H_N(vu)$ and its derivatives in $L^1$ of the uniform norm on $\{|v|\le R\}\times\mathbb U_N$, for each fixed $N$.
\end{lemma}

\begin{proof}
\smallskip
\noindent\emph{Step 1: The $L^2(\P)$ bounds.} Let $H_{N,k}$ be the degree-$k$ summand of \eqref{eq:mixed-Hamiltonian}. We work first with $H_N^{(K)}$, uniformly in $K$. These finite sums can be differentiated term by term. Step 3 justifies this operation for the full random series. Cauchy--Schwarz gives
\begin{equation}
\label{eq:matrix-derivative-summability}
\sum_{k\ge1}k^2(1\vee R)^k|\Delta_k|_{\mathrm{op}}^{1/2}\le\Big(\sum_{k\ge1}4^{-k}\Big)^{1/2}\Big(\sum_{k\ge1}k^4(2(1\vee R))^{2k}|\Delta_k|\Big)^{1/2}<\infty.
\end{equation}
For the last inequality, use $k^4\le C2^k$ and apply \eqref{eq:mixture-summability} with radius $8(1\vee R)^2$. Differentiating replaces one of the $k$ factors $vu$ by $bu$. By \eqref{eq:g-covariance}, each resulting contraction has $L^2(\P)$ norm at most $|\Delta_k|_{\mathrm{op}}^{1/2}$ times the product of its factor norms. Summing these $k$ bounds gives
\begin{equation}
\label{eq:matrix-derivative-L2}
\big|(b\cdot\nabla_v)H_{N,k}(vu)\big|_{L^2(\P)}\le k|\Delta_k|_{\mathrm{op}}^{1/2}N^{-(k-1)/2}|vu|^{k-1}|bu|=k|\Delta_k|_{\mathrm{op}}^{1/2}\sqrt N\,|v|^{k-1}|b|.
\end{equation}
Here $|vu|^2=\tr(vuu^{\intercal}v^{\intercal})=N|v|^2$ and $|bu|^2=N|b|^2$. To compare triples $(v,b,u)$ and $(v',b',u')$ with $|v|,|v'|\le R$ and $|b|,|b'|\le1$, we use
\begin{equation*}
|vu-v'u'|\le\sqrt N|v-v'|+R|u-u'|\quad\text{and}\quad|bu-b'u'|\le\sqrt N|b-b'|+|u-u'|.
\end{equation*}
Telescope each of the $k$ differentiated products through its $k$ factors, as in \eqref{eq:tensor-telescoping}. There are at most $k^2$ resulting terms. After multiplying by $N^{-(k-1)/2}$, each term is bounded in $L^2(\P)$ by
\begin{equation*}
|\Delta_k|_{\mathrm{op}}^{1/2}(1\vee R)^k\sqrt N\left(|v-v'|+|b-b'|+N^{-1/2}|u-u'|\right).
\end{equation*}
Summing over $k\in\{1,\ldots,K\}$ and using \eqref{eq:matrix-derivative-summability} gives
\begin{equation}
\label{eq:matrix-derivative-distance}
\big|(b\cdot\nabla_v)H_N^{(K)}(vu)-(b'\cdot\nabla_{v'})H_N^{(K)}(v'u')\big|_{L^2(\P)}\le C_R\sqrt N\left(|v-v'|+|b-b'|+N^{-1/2}|u-u'|\right).
\end{equation}

\smallskip
\noindent\emph{Step 2: From finite nets to the supremum.} Use the distance $|v-v'|+|b-b'|+N^{-1/2}|u-u'|$ on the triples. The rescaled triples $(v,b,N^{-1/2}u)$ lie in Euclidean balls of dimensions $D^2$, $D^2$, and $DN$, with radii $R$, $1$, and $\sqrt D$. The volume comparison for disjoint balls of radius $\varepsilon/2$ therefore gives an $\varepsilon$-net with at most $(C_R/\varepsilon)^{DN+2D^2}$ points, for $\varepsilon\in(0,1]$.

For centered Gaussian variables $X_1,\ldots,X_m$ with variances at most $\tau^2$ and any $\lambda>0$, we have
\begin{equation*}
\E\max_{\ell\le m}|X_\ell|\le\lambda^{-1}\log\sum_{\ell=1}^m\E\big(e^{\lambda X_\ell}+e^{-\lambda X_\ell}\big)\le\frac{\log(2m)}\lambda+\frac{\lambda\tau^2}{2}.
\end{equation*}
Minimizing over $\lambda$ gives $\E\max_{\ell\le m}|X_\ell|\le\tau\sqrt{2\log(2m)}$, without requiring independence.

Choose nets at distances $2^{-j}$ for $j\in\{0,1,\ldots\}$. For each triple, choose an approximating point in each net. Two successive approximations are at distance at most $3\cdot2^{-j}$ for $j\ge1$. By \eqref{eq:matrix-derivative-distance}, their Gaussian difference has variance at most $C_RN2^{-2j}$. The number of possible pairs is at most $(C_R2^j)^{2(DN+2D^2)}$, so the expected maximum of these absolute differences is at most
\begin{equation*}
C_R\sqrt N\,2^{-j}\sqrt{(DN+2D^2)(j+1)}.
\end{equation*}
The expected maximum on the initial net is at most $C_R\sqrt{N(DN+2D^2)}$, by \eqref{eq:matrix-derivative-L2} and \eqref{eq:matrix-derivative-summability}. For fixed $K$, the derivative field is continuous, so telescoping the successive net approximations gives
\begin{equation*}
\E\sup_{\substack{|v|\le R,\ |b|\le1\\u\in\mathbb U_N}}|(b\cdot\nabla_v)H_N^{(K)}(vu)|\le C_R\sqrt{N(DN+2D^2)}\bigg(1+\sum_{j\ge1}2^{-j}\sqrt{j+1}\bigg)\le C_RN.
\end{equation*}
The last inequality uses $N\ge D$, fixed $D$, and convergence of the numerical series.

\smallskip
\noindent\emph{Step 3: Passing to the full Hamiltonian.} Apply Steps 1 and 2 to $H_N^{(L)}-H_N^{(K)}$, where $L>K$. The standard-deviation bounds in Step 1 are then controlled by $\sum_{k>K}k^2(1\vee R)^k|\Delta_k|_{\mathrm{op}}^{1/2}$. Consequently, we have
\begin{equation*}
\E\sup_{\substack{|v|\le R,\ |b|\le1\\u\in\mathbb U_N}}|(b\cdot\nabla_v)(H_N^{(L)}-H_N^{(K)})(vu)|\le C_RN\sum_{k>K}k^2(1\vee R)^k|\Delta_k|_{\mathrm{op}}^{1/2}.
\end{equation*}
The right-hand side tends to zero as $K\to\infty$, uniformly in $L>K$, by \eqref{eq:matrix-derivative-summability}. Since every partial sum vanishes at $v=0$, the fundamental theorem of calculus along $\lambda v$, $\lambda\in[0,1]$, also gives
\begin{equation*}
\sup_{\substack{|v|\le R\\u\in\mathbb U_N}}|(H_N^{(L)}-H_N^{(K)})(vu)|\le R\sup_{\substack{|v|\le R,\ |b|\le1\\u\in\mathbb U_N}}|(b\cdot\nabla_v)(H_N^{(L)}-H_N^{(K)})(vu)|.
\end{equation*}
Thus the partial sums and their derivatives are Cauchy in $L^1$ of the uniform norm. A subsequence converges uniformly almost surely for both quantities. Passing to the limit in the line-integral identity identifies the derivative of the limit. One may work on the ball of radius $R+1$ to justify derivatives at $|v|=R$. The pointwise $L^2(\P)$ limit is the Hamiltonian in \eqref{eq:mixed-Hamiltonian}, so this constructs its continuously differentiable version in $v$. Passing to the limit in Step 2 proves \eqref{eq:matrix-interaction-derivative-sup}.
\end{proof}

For $\kappa\in(0,1)$, set 
\begin{equation}
\label{eq:Skappa-enriched}
\mcl A_\kappa:=\{a\in\mcl A:\det(a^2)\ge\kappa\}\quad\text{and}\quad \mcl B_\kappa:=\{b\in\mcl B:\det b\ge\kappa\}.
\end{equation}
The map $b\mapsto b^{1/2}$ is a bijection from $\mcl B_\kappa$ to $\mcl A_\kappa$. These sets are compact and contain $I_D$. Indeed, $\tr b=D$ gives $\lambda_{\max}(b)\le D$, and $\det b\ge\kappa$ gives $\lambda_{\min}(b)\ge\kappa D^{-(D-1)}>0$. We will first integrate over $\mcl B_\kappa$ and then estimate its complement. The next estimates yield convergence uniform on $\mcl A_\kappa$, which will justify replacing the random integrand by its limit.

\begin{lemma}
\label{lem:enriched-a-regularity}
Fix $t\ge0$, $q\in\mcl Q_\infty$, and $\kappa\in(0,1)$. There exists $C_{t,q,h,\kappa}<\infty$, independent of $N$, such that, for every $N\in\N$ with $N\ge D$ and every $\delta\in(0,1)$, we have
\begin{equation}
\label{eq:enriched-a-modulus}
    \E\sup_{\substack{a,a'\in \mcl A_\kappa\\|a-a'|\le\delta}}\left|F_{N,a}(t,q_a)-F_{N,a'}(t,q_{a'})\right|\le C_{t,q,h,\kappa}\delta.
\end{equation}
Moreover, for every $a\in\mcl A^\circ$, we have
\begin{equation}
\label{eq:enriched-a-concentration}
\E\left|F_{N,a}(t,q_a)-\bar F_{N,a}(t,q_a)\right|\le C_{t,q}N^{-1/2}.
\end{equation}
\end{lemma}

\begin{proof}
\smallskip
\noindent\emph{Step 1: Proof of \eqref{eq:convex-matrix-modulus}.} Fix $R>2\sqrt D$ and put $B_R:=\{w\in\R^{D\times D}:|w|\le R\}$. We first prove that every convex function $g:[-2R,2R]^{D^2}\to\R$ satisfies
\begin{equation}
\label{eq:convex-matrix-modulus}
\operatorname{Lip}_{B_R}g\le C_{D,R}\bigg(\sum_{w\in\{\pm2R\}^{D^2}}|g(w)|+|g(0)|\bigg).
\end{equation}
Put $M:=\max_{v\in\{\pm2R\}^{D^2}}g(v)$. For $w$ in the cube, write $w=\sum_v\theta_vv$ as a convex combination of its vertices, where $\theta_v\ge0$ and $\sum_v\theta_v=1$. Convexity gives
\begin{equation*}
g(w)\le\sum_v\theta_vg(v)\le M\qquad\text{and}\qquad g(0)=g\Ll(\tfrac12w+\tfrac12(-w)\Rr)\le\tfrac12g(w)+\tfrac12g(-w)\le\tfrac12g(w)+\tfrac12M.
\end{equation*}
Rearranging the second inequality and using the first gives
\begin{equation}
\label{eq:convex-cube-bounds}
2g(0)-M\le g(w)\le M,\qquad w\in[-2R,2R]^{D^2}.
\end{equation}
For $w\in B_R$ and $|e|=1$, both $w-Re$ and $w+Re$ belong to the cube. The secant slopes of the convex function $s\mapsto g(w+se)$ satisfy, for $s\in(0,R]$,
\begin{equation*}
\frac{g(w)-g(w-Re)}R\le\frac{g(w+se)-g(w)}s\le\frac{g(w+Re)-g(w)}R.
\end{equation*}
The difference of any two values in \eqref{eq:convex-cube-bounds} has absolute value at most $2(M-g(0))$. Thus both outer quotients have absolute value at most $2(M-g(0))/R$, and
\begin{equation}
\label{eq:convex-directional-modulus}
|g(w+se)-g(w)|\le\frac{2(M-g(0))}{R}s,\qquad w\in B_R,\quad |e|=1,\quad s\in(0,R].
\end{equation}
For $w,w'\in B_R$ with $0<|w'-w|\le R$, take $e=(w'-w)/|w'-w|$ and $s=|w'-w|$ in \eqref{eq:convex-directional-modulus}. If $|w'-w|>R$, their midpoint belongs to $B_R$, and both half-segments have length $|w'-w|/2\le R$. Applying the same bound to the two half-segments and adding gives
\begin{equation*}
|g(w')-g(w)|\le\frac{2(M-g(0))}{R}\left(\frac{|w'-w|}{2}+\frac{|w'-w|}{2}\right)=\frac{2(M-g(0))}{R}|w'-w|.
\end{equation*}
The case $w=w'$ is immediate. Combining these estimates gives
\begin{equation*}
\operatorname{Lip}_{B_R}g\le\frac{2(M-g(0))}{R}\le\frac2R\bigg(\sum_{v\in\{\pm2R\}^{D^2}}|g(v)|+|g(0)|\bigg),
\end{equation*}
which proves \eqref{eq:convex-matrix-modulus}.

\smallskip
\noindent\emph{Step 2: Definition of $L_N$ and proof of \eqref{eq:matrix-enrichment-vertex-bound}.} Define $L_N:\R^{D\times D}\times\R^{D\times D}\to\R$ by
\begin{equation*}
L_N(v,w):=\frac1N\log\iint\exp\left(\sqrt{2t}H_N(vu)+\sqrt2\sum_{i=1}^Nz_i^q(\alpha)\cdot wu_i+h\cdot\sum_{i=1}^Nvu_i\right)\nu_N(\mathrm d u)\mathfrak R(\mathrm d\alpha).
\end{equation*}
The integral is over $\mathbb U_N\times\supp\mathfrak R$. The matrix $v$ varies the interaction and magnetic-field terms, whereas $w$ varies the Gaussian linear term. For fixed $v$, H\"older's inequality gives convexity of $w\mapsto L_N(v,w)$, since the exponent is affine in $w$. We first bound $L_N(0,w)$ at the vertices needed in \eqref{eq:convex-matrix-modulus}.

By \eqref{eq:stiefel-right-invariance} with $o=-I_N$, the measure $\nu_N$ is invariant under $u\mapsto-u$. Consequently, $\int u_i\nu_N(\mathrm d u)=0$ for each $i\in\{1,\ldots,N\}$. For each $\alpha\in\supp\mathfrak R$, Jensen's inequality therefore gives
\begin{equation*}
\int_{\mathbb U_N}\exp\left(\sqrt2\sum_{i=1}^Nz_i^q(\alpha)\cdot wu_i\right)\nu_N(\mathrm d u)\ge\exp\left(\sqrt2\sum_{i=1}^Nz_i^q(\alpha)\cdot w\int_{\mathbb U_N}u_i\nu_N(\mathrm d u)\right)=1.
\end{equation*}
Integrating over the probability measure $\mathfrak R$ gives $L_N(0,w)\ge0$ almost surely. For fixed $u\in\mathbb U_N$ and $\alpha\in\supp\mathfrak R$, we also have
\begin{equation*}
\E\left[\left(\sqrt2\sum_{i=1}^Nz_i^q(\alpha)\cdot wu_i\right)^2\bigg|\mathfrak R\right]=2\sum_{i=1}^Nu_i^{\intercal}w^{\intercal}q(1)wu_i=2N\tr(w^{\intercal}q(1)w).
\end{equation*}
The last equality uses $uu^{\intercal}=NI_D$. Gaussian integration and Tonelli's theorem thus give $\E\exp(NL_N(0,w))=\exp(N\tr(w^{\intercal}q(1)w))$. Hence, we have
\begin{equation}
\label{eq:matrix-enrichment-vertex-bound}
0\le\E|L_N(0,w)|=\E L_N(0,w)\le N^{-1}\log\E e^{NL_N(0,w)}=\operatorname{tr}(w^{\intercal}q(1)w).
\end{equation}
This bounds the expectations at the finitely many vertices of $[-2R,2R]^{D^2}$. It also gives finiteness there almost surely. Convexity gives finiteness throughout the cube. Applying this argument on cubes with integer radii gives finiteness of $L_N(0,w)$ for every $w$ on a common probability-one event.

\smallskip
\noindent\emph{Step 3: Proof of \eqref{eq:enriched-a-modulus}.} A uniform bound on the difference of two exponents bounds the difference of the logarithms of their integrals. Applying this fact to $L_N$, we obtain
\begin{equation*}
\sup_w|L_N(v,w)-L_N(v',w)|\le\frac{\sqrt{2t}}N\sup_{u\in\mathbb U_N}|H_N(vu)-H_N(v'u)|+|h|\,|v-v'|.
\end{equation*}
For the magnetic-field term, we used $|N^{-1}\sum_i u_i|\le N^{-1}|u|_{\mathrm{op}}|\mathbf1_N|=1$. Since $H_N(0)=0$, integration along the segment from $0$ to $v$ gives
\begin{equation}
\label{eq:matrix-factor-value-comparison}
\sup_{\substack{v\in B_R\\w\in\R^{D\times D}}}|L_N(v,w)-L_N(0,w)|\le\frac{R\sqrt{2t}}N\sup_{\substack{v\in B_R,\ |b|\le1\\u\in\mathbb U_N}}|(b\cdot\nabla_v)H_N(vu)|+R|h|.
\end{equation}
In particular, Lemma~\ref{lem:matrix-interaction-derivatives} and Step 2 give finiteness of $L_N(v,w)$ for all $v,w$ almost surely. Apply \eqref{eq:convex-matrix-modulus} to $g(w)=L_N(v,w)$. For each vertex $w$ and also for $w=0$, we have
\begin{equation*}
\sup_{v\in B_R}|L_N(v,w)|\le |L_N(0,w)|+\sup_{\substack{v\in B_R\\w'\in\R^{D\times D}}}|L_N(v,w')-L_N(0,w')|.
\end{equation*}
The expectation of the first term is bounded by \eqref{eq:matrix-enrichment-vertex-bound}. The second is bounded by \eqref{eq:matrix-factor-value-comparison} and \eqref{eq:matrix-interaction-derivative-sup}. Since there are finitely many vertices, we obtain
\begin{equation}
\label{eq:matrix-enrichment-uniform-lipschitz}
\E\sup_{v\in B_R}\operatorname{Lip}_{B_R}L_N(v,\cdot)\le C_{t,q,h,R}.
\end{equation}
For $v,v',w,w'\in B_R$, split the difference as
\begin{equation*}
|L_N(v,w)-L_N(v',w')|\le|L_N(v,w)-L_N(v',w)|+|L_N(v',w)-L_N(v',w')|.
\end{equation*}
Integration along the segment from $v$ to $v'$ and \eqref{eq:matrix-interaction-derivative-sup} control the first term. Equation \eqref{eq:matrix-enrichment-uniform-lipschitz} controls the second. Thus, for every $\delta>0$, we have
\begin{equation}
\label{eq:matrix-factor-modulus}
\E\sup_{\substack{v,v',w,w'\in B_R\\|v-v'|+|w-w'|\le\delta}}|L_N(v,w)-L_N(v',w')|\le C_{t,q,h,R}\delta.
\end{equation}
For $a\in\mcl A^\circ$, the realization in \eqref{eq:enriched-restriction-H} gives
\begin{equation*}
F_{N,a}(t,q_a)=-L_N(a,a)+t\xi(a^2)+q(1)\cdot a^2.
\end{equation*}
Every $a\in\mcl A$ satisfies $|a|=\sqrt D<R$. Moreover, we have
\begin{equation*}
|a^2-(a')^2|=|a(a-a')+(a-a')a'|\le2\sqrt D\,|a-a'|,\qquad a,a'\in\mcl A.
\end{equation*}
The function $\xi$ is Lipschitz on the compact set of positive-semidefinite matrices with trace at most $D$. More precisely, with $C:=t\sup_{r\in\S_+^D,\,\tr r\le D}|\nabla\xi(r)|+|q(1)|$, we have
\begin{equation*}
|t\xi(a^2)+q(1)\cdot a^2-t\xi((a')^2)-q(1)\cdot(a')^2|\le2\sqrt D\,C|a-a'|.
\end{equation*}
Applying \eqref{eq:matrix-factor-modulus} with $v=w=a$, $v'=w'=a'$, and $2\delta$ in place of $\delta$, then adding this deterministic bound, proves \eqref{eq:enriched-a-modulus}.

\smallskip
\noindent\emph{Step 4: Proof of \eqref{eq:enriched-a-concentration}.} For fixed $u\in\mathbb U_N$ and $\alpha\in\supp\mathfrak R$, conditional independence of the interaction and enrichment fields gives
\begin{equation*}
\E\left[\left(\sqrt{2t}H_N(au)+\sqrt2\sum_{i=1}^Nz_i^q(\alpha)\cdot au_i\right)^2\Bigm|\mathfrak R\right]=2N\big(t\xi(a^2)+\tr(aq(1)a)\big)\le C_{t,q}N.
\end{equation*}
Here $\tr(aq(1)a)\le D|q(1)|_{\mathrm{op}}$, and $\xi$ is bounded on $\{a^2:a\in\mcl A\}$. Apply Lemma~\ref{lem:total-concentration-continuous} with reference measure $\nu_N$, interaction field $\sqrt{2t}H_N(au)$, and enrichment fields $\sqrt2az_i^q$. Its deterministic function $A:\mathbb U_N\to\R$ is $A(u):=ah\cdot\sum_i u_i-Nt\xi(a^2)-Nq_a(1)\cdot I_D$, which is bounded for each $N$. Consequently, we have
\begin{equation*}
\E|F_{N,a}(t,q_a)-\E F_{N,a}(t,q_a)|\stackrel{\eqref{eq:can-total-concentration}}{\le}\frac C N\sqrt{1+C_{t,q}N}\le C_{t,q}N^{-1/2}.
\end{equation*}
Since $\E F_{N,a}(t,q_a)=\bar F_{N,a}(t,q_a)$, this proves \eqref{eq:enriched-a-concentration}.
\end{proof}

Recall from \eqref{eq:enriched-fixed-a} that $F_{N,a}(t,q_a)=-N^{-1}\log Z_{N,a}(t,q_a)$ and $\bar F_{N,a}(t,q_a)=\E F_{N,a}(t,q_a)$. Thus Corollary~\ref{cor:fixed-a-limit} gives $\lim_{N\to\infty}\E F_{N,a}(t,q_a)=f_a(t,q_a)$. To control the random error $F_{N,a}(t,q_a)-f_a(t,q_a)$, we use the identity
\begin{equation}
F_{N,a}(t,q_a)-f_a(t,q_a)=\big(F_{N,a}(t,q_a)-\E F_{N,a}(t,q_a)\big)+\big(\bar F_{N,a}(t,q_a)-f_a(t,q_a)\big).
\label{eq:enriched-error-decomposition}
\end{equation}
Lemma~\ref{lem:enriched-a-regularity} bounds the first term in $L^1$ and supplies a uniform modulus in $a$. Corollary~\ref{cor:fixed-a-limit} makes the second term tend to zero for each $a$. A finite-net argument then gives the following uniform limits on $\mcl A_\kappa$.
\begin{lemma}
\label{lem:enriched-local-uniform}
For every $t\ge0$, $q\in\mcl Q_\infty$, and $\kappa\in(0,1)$, we have
\begin{equation}
\label{eq:enriched-local-uniform-deterministic}
\lim_{N\to\infty}\sup_{a\in \mcl A_\kappa} |\bar F_{N,a}(t,q_a)-f_a(t,q_a)|=0\quad\text{and}\quad
\lim_{N\to\infty}\E\sup_{a\in \mcl A_\kappa} \left|F_{N,a}(t,q_a)-f_a(t,q_a)\right|=0.
\end{equation}
\end{lemma}

\begin{proof}
\smallskip
\noindent\emph{Step 1: Continuity of $a\mapsto f_a(t,q_a)$.} Let $a,a'\in\mcl A_\kappa$ satisfy $|a-a'|<1$. If $a=a'$, the bound below is immediate. Otherwise, apply \eqref{eq:enriched-a-modulus} with $\delta=|a-a'|\in(0,1)$. Since $\bar F_{N,a}(t,q_a)=\E F_{N,a}(t,q_a)$, the triangle inequality for the expectation gives
\begin{equation*}
|\bar F_{N,a}(t,q_a)-\bar F_{N,a'}(t,q_{a'})|\le\E\left|F_{N,a}(t,q_a)-F_{N,a'}(t,q_{a'})\right|\stackrel{\eqref{eq:enriched-a-modulus}}{\le}C_{t,q,h,\kappa}|a-a'|.
\end{equation*}
By Corollary~\ref{cor:fixed-a-limit}, the two deterministic terms converge to $f_a(t,q_a)$ and $f_{a'}(t,q_{a'})$, respectively. Hence, we have
\begin{equation}
\label{eq:enriched-limit-a-modulus}
|f_a(t,q_a)-f_{a'}(t,q_{a'})|\le C_{t,q,h,\kappa}|a-a'|,\qquad a,a'\in \mcl A_\kappa,\quad |a-a'|<1.
\end{equation}
In particular, the limiting function is continuous on $\mcl A_\kappa$.

\smallskip
\noindent\emph{Step 2: Approximation by finitely many $a_j$.} Fix $\delta\in(0,1)$ and choose a finite $\delta$-net $a_1,\ldots,a_m$ of the compact set $\mcl A_\kappa$. For each $a\in \mcl A_\kappa$, choose $j\in\{1,\ldots,m\}$ with $|a-a_j|\le\delta$. Apply \eqref{eq:enriched-error-decomposition} at $a_j$ and use $\E F_{N,a_j}=\bar F_{N,a_j}$. The triangle inequality gives
\begin{align*}
\left|F_{N,a}(t,q_a)-f_a(t,q_a)\right|
&\le\left|F_{N,a}(t,q_a)-F_{N,a_j}(t,q_{a_j})\right|+\left|F_{N,a_j}(t,q_{a_j})-\bar F_{N,a_j}(t,q_{a_j})\right|\\
&\quad+|\bar F_{N,a_j}(t,q_{a_j})-f_{a_j}(t,q_{a_j})|+|f_{a_j}(t,q_{a_j})-f_a(t,q_a)|.
\end{align*}
Take the supremum over $a$ and then expectation. The first term is controlled by \eqref{eq:enriched-a-modulus}, and the last by \eqref{eq:enriched-limit-a-modulus}. For the second term, bound the maximum over $j$ by the sum and use \eqref{eq:enriched-a-concentration}. We obtain
\begin{align}
\label{eq:enriched-net-error}
\E\sup_{a\in \mcl A_\kappa}\left|F_{N,a}(t,q_a)-f_a(t,q_a)\right|\le2C_{t,q,h,\kappa}\delta+\max_{j\le m}|\bar F_{N,a_j}(t,q_{a_j})-f_{a_j}(t,q_{a_j})|+C_{t,q}mN^{-1/2}.
\end{align}
For fixed $\delta$, both $m$ and the points $a_j$ are independent of $N$. Corollary~\ref{cor:fixed-a-limit} therefore makes the maximum in \eqref{eq:enriched-net-error} tend to zero, and $mN^{-1/2}$ tends to zero as well. Thus we have
\begin{equation*}
\limsup_{N\to\infty}\E\sup_{a\in \mcl A_\kappa}\left|F_{N,a}(t,q_a)-f_a(t,q_a)\right|\le2C_{t,q,h,\kappa}\delta.
\end{equation*}
Letting $\delta\downarrow0$ proves the second identity in \eqref{eq:enriched-local-uniform-deterministic}. Since $\bar F_{N,a}(t,q_a)-f_a(t,q_a)=\E\big(F_{N,a}(t,q_a)-f_a(t,q_a)\big)$, the triangle inequality for the expectation gives
\begin{equation*}
\sup_{a\in \mcl A_\kappa}|\bar F_{N,a}(t,q_a)-f_a(t,q_a)|\le\E\sup_{a\in \mcl A_\kappa}\left|F_{N,a}(t,q_a)-f_a(t,q_a)\right|,
\end{equation*}
which proves the first identity in \eqref{eq:enriched-local-uniform-deterministic}.
\end{proof}

For $\kappa\in(0,1)$, let $Z_{N}^{\geq\kappa}(t,q)$ and $Z_{N}^{<\kappa}(t,q)$ denote the parts of \eqref{eq:global-enriched-disintegration} over $\mcl B_\kappa$ and $\{b\in\mcl B^\circ:\det b<\kappa\}$, respectively. Thus $Z_N=Z_{N}^{\geq\kappa}+Z_{N}^{<\kappa}$. The next estimate controls the second integral through its expectation. 

\begin{lemma}
\label{lem:enriched-boundary}
Fix $t\ge0$ and $q\in\mcl Q_\infty$. For every $a\in\mcl A^\circ$ and $N\in\N$ with $N\ge D$, we have
\begin{equation}
\label{eq:enriched-annealed-one}
\E Z_{N,a}(t,q_a)=\int_{\mathbb U_N}e^{ah\cdot\sum_{i=1}^Nu_i}\nu_N(\mathrm d u)\le e^{N|ah|}\le e^{N\sqrt D\,|h|}.
\end{equation}
For every $\kappa\in(0,1)$, we have
\begin{equation}
\label{eq:enriched-boundary-estimate}
\limsup_{N\to\infty}\frac1N\E\log Z_{N}^{<\kappa}(t,q) \le\sqrt D\,|h|+\frac12\log\kappa.
\end{equation}
\end{lemma}

\begin{proof}
\smallskip
\noindent\emph{Step 1: Computing $\E Z_{N,a}(t,q_a)$.} Recall the interaction covariance \eqref{e.EHH=xi} and the enrichment covariance \eqref{eq:enriched-cascade-covariance}. Gaussian integration therefore gives
\begin{equation*}
\E\exp\Ll(\sqrt{2t}H_N(au)-Nt\xi(a^2)\Rr)=1\qquad\text{and}\qquad\E\exp\Ll(\sqrt2\sum_{i=1}^Nz_i^q(\alpha)\cdot au_i-Nq_a(1)\cdot I_D\Rr)=1.
\end{equation*}
Integrate first over the Gaussian fields conditional on the cascade. Only $e^{ah\cdot\sum_i u_i}$ remains. Since $\mathfrak R$ has mass one, Tonelli's theorem gives the equality in \eqref{eq:enriched-annealed-one}. Also, $uu^\intercal=NI_D$ gives $|N^{-1}\sum_i u_i|\le1$, and $|a|_{\mathrm{op}}\le\sqrt D$. These two bounds give the inequalities in \eqref{eq:enriched-annealed-one}.

\smallskip
\noindent\emph{Step 2: The integral over $\det b<\kappa$.} When $D=1$, this integral is zero, so \eqref{eq:enriched-boundary-estimate} holds with $\log0=-\infty$. Suppose $D\ge2$. The set $\mcl B$ has finite, positive off-diagonal Lebesgue volume because its off-diagonal entries belong to $[-1,1]$ and $I_D$ is in its relative interior. For $N>D+1$, the exponent $(N-D-1)/2$ is positive, so $(\det b)^{(N-D-1)/2}\le\kappa^{(N-D-1)/2}$ on $\det b<\kappa$. Tonelli's theorem and \eqref{eq:enriched-annealed-one} therefore give
\begin{equation*}
\E Z_{N}^{<\kappa}(t,q)\le c_{N,D}e^{N\sqrt D\,|h|}\int_{\substack{b\in\mcl B^\circ\\\det b<\kappa}}(\det b)^{(N-D-1)/2}\,\mathrm d b_{\mathrm{off}}\le c_{N,D}e^{N\sqrt D\,|h|}\operatorname{vol}(\mcl B)\kappa^{(N-D-1)/2}.
\end{equation*}
Jensen's inequality yields
\begin{equation*}
\frac1N\E\log Z_{N}^{<\kappa}(t,q)\le\sqrt D\,|h|+\frac{\log c_{N,D}+\log\operatorname{vol}(\mcl B)}N+\frac{N-D-1}{2N}\log\kappa.
\end{equation*}
For fixed $\kappa$, the term $N^{-1}\log c_{N,D}$ tends to zero by \eqref{eq:normalization-subexp}, $N^{-1}\log\operatorname{vol}(\mcl B)$ tends to zero, and $(N-D-1)/(2N)$ tends to $1/2$. Taking the upper limit proves \eqref{eq:enriched-boundary-estimate}.
\end{proof}

We can now combine the local uniform convergence in Lemma~\ref{lem:enriched-local-uniform} with the boundary estimate in Lemma~\ref{lem:enriched-boundary} to prove the main result of this section. Recall $\bar F_N$ from \eqref{eq:global-enriched-intro}, $f_a$ from \eqref{eq:Hopf-fixed-a-intro}, $q_a$ from \eqref{eq:qa-intro}, and $\mcl A^\circ$ from \eqref{eq:elliptope}.

\begin{proposition}
\label{prop:global-enriched}
For every $t\ge0$ and $q\in\mcl Q_2$, we have 
\begin{equation}
\label{eq:global-enriched-limit}
\lim_{N\to\infty}\bar F_N(t,q) =\inf_{a\in\mcl A^\circ}\left\{f_a(t,q_a)-\log\det a\right\}.
\end{equation}
The infimum is attained in $\mcl A^\circ$.
\end{proposition}

\begin{proof}
When $D=1$, $\mcl A^\circ=\{I_1\}$, $q_{I_1}=q$, and $P_N=\nu_N$. The conclusion, including attainment, is therefore Corollary~\ref{cor:fixed-a-limit} for every $q\in\mcl Q_2$. Assume $D\ge2$. We first prove the result for $q\in\mcl Q_\infty$ and then pass to $q\in\mcl Q_2$.

\smallskip
\noindent\emph{Step 1: Replacing $F_{N,a}(t,q_a)$ by $f_a(t,q_a)$.} Fix $\kappa\in(0,1)$, and recall $\mcl A_\kappa$ and $\mcl B_\kappa$ from \eqref{eq:Skappa-enriched}, and $Z_N^{\geq\kappa}(t,q)$ from the paragraph before Lemma~\ref{lem:enriched-boundary}. Lemma~\ref{lem:enriched-local-uniform} gives continuity of $a\mapsto-f_a(t,q_a)+\log\det a$ on $\mcl A_\kappa$. For this comparison only, put $r_N:=\sup_{a\in\mcl A_\kappa}|F_{N,a}(t,q_a)-f_a(t,q_a)|$. Then $\lim_{N\to\infty}\E r_N=0$ by the second identity in \eqref{eq:enriched-local-uniform-deterministic}, and for every $a\in\mcl A_\kappa$ we have
\begin{equation*}
e^{-Nr_N}e^{-Nf_a(t,q_a)}\le Z_{N,a}(t,q_a)=e^{-NF_{N,a}(t,q_a)}\le e^{Nr_N}e^{-Nf_a(t,q_a)}.
\end{equation*}
Multiply these inequalities by the determinant density in \eqref{eq:global-enriched-disintegration} and integrate over $\mcl B_\kappa$. Thus $Z_N^{\geq\kappa}(t,q)$ lies between $e^{-Nr_N}$ and $e^{Nr_N}$ times the same integral with $Z_{N,b^{1/2}}(t,q_{b^{1/2}})$ replaced by $e^{-Nf_{b^{1/2}}(t,q_{b^{1/2}})}$ (see~\eqref{eq:global-enriched-disintegration} with $\mcl B^\circ$ replaced by $\mcl B_\kappa$). Taking logarithms and expectations gives
\begin{align}
\label{eq:enriched-Laplace-error}
\left|\frac1N\E\log Z_N^{\geq\kappa}(t,q)-\frac1N\log\left(c_{N,D}\int_{\mcl B_\kappa}e^{-Nf_{b^{1/2}}(t,q_{b^{1/2}})}(\det b)^{(N-D-1)/2}\,\d b_{\mathrm{off}}\right)\right|\notag\\
\le\E\sup_{a\in\mcl A_\kappa}|F_{N,a}(t,q_a)-f_a(t,q_a)|.
\end{align}
The right-hand side tends to zero by the second identity in \eqref{eq:enriched-local-uniform-deterministic}.

\smallskip
\noindent\emph{Step 2: The exponential rate on $\mcl B_\kappa$.} For $b\in\mcl B_\kappa\subset\mcl B=\{a^2:a\in\mcl A\}$, \eqref{eq:elliptope} gives $b_{dd}=1$ for every $d\in\{1,\dots,D\}$. Hadamard's inequality therefore gives $\det b\le1$. Since $\log\det(b^{1/2})=\frac12\log\det b$, the deterministic integrand in \eqref{eq:enriched-Laplace-error} satisfies
\begin{equation*}
e^{-Nf_{b^{1/2}}(t,q_{b^{1/2}})}(\det b)^{(N-D-1)/2}=\exp\left(N\left[-f_{b^{1/2}}(t,q_{b^{1/2}})+\frac12\log\det b\right]\right)(\det b)^{-(D+1)/2}.
\end{equation*}
On $\mcl B_\kappa$, the last factor satisfies
\begin{equation*}
1\le(\det b)^{-(D+1)/2}\le\kappa^{-(D+1)/2}.
\end{equation*}
We first compute the exponential rate after removing this bounded factor. The upper bound is
\begin{equation}
\label{eq:compact-Laplace-upper}
\frac1N\log\int_{\mcl B_\kappa}e^{-Nf_{b^{1/2}}(t,q_{b^{1/2}})}(\det b)^{N/2}\,\d b_{\mathrm{off}}
\le\sup_{a\in\mcl A_\kappa}\{-f_a(t,q_a)+\log\det a\}+\frac1N\log\operatorname{vol}(\mcl B_\kappa).
\end{equation}
For the lower bound, choose a maximizer $a_*$ of $a\mapsto-f_a(t,q_a)+\log\det a$ on $\mcl A_\kappa$ and put $b_*:=a_*^2$. We first move $b_*$ slightly toward $I_D$ so that a neighborhood of the resulting matrix has determinant greater than $\kappa$. For $b\in\mcl B_\kappa$ and $\rho\in(0,1)$, concavity of $\log\det$ gives
\begin{equation*}
\log\det((1-\rho)b+\rho I_D)\ge(1-\rho)\log\det b>\log\kappa.
\end{equation*}
Here the strict inequality holds because $\log\det b\ge\log\kappa$ and $\log\kappa<0$. Fix $\ep>0$. The function $b\mapsto-f_{b^{1/2}}(t,q_{b^{1/2}})+\frac12\log\det b$ is continuous on $\mcl B_\kappa$. Hence, for sufficiently small $\rho>0$, we can choose an open ball $U$ in the off-diagonal coordinates, centered at $(1-\rho)b_*+\rho I_D$, such that $U\subset\mcl B_\kappa$ and
\begin{equation*}
-f_{b^{1/2}}(t,q_{b^{1/2}})+\frac12\log\det b\ge-f_{a_*}(t,q_{a_*})+\log\det a_*-\ep,\qquad b\in U.
\end{equation*}
The containment is possible because positive definiteness and $\det b>\kappa$ are strict at the center. The ball $U$ has strictly positive off-diagonal Lebesgue volume. Restricting the integral to $U$ gives
\begin{equation}
\label{eq:compact-Laplace-lower}
\frac1N\log\int_{\mcl B_\kappa}e^{-Nf_{b^{1/2}}(t,q_{b^{1/2}})}(\det b)^{N/2}\,\d b_{\mathrm{off}}
\ge-f_{a_*}(t,q_{a_*})+\log\det a_*-\ep+\frac1N\log\operatorname{vol}(U).
\end{equation}
The sets $\mcl B_\kappa$ and $U$ have finite, positive volume independent of $N$. Sending first $N\to\infty$ and then $\ep\downarrow0$ in \eqref{eq:compact-Laplace-upper}--\eqref{eq:compact-Laplace-lower} gives $\sup_{a\in\mcl A_\kappa}\{-f_a(t,q_a)+\log\det a\}$ as the limit. Restoring the factor $(\det b)^{-(D+1)/2}$ changes the normalized logarithm by at most $(D+1)|\log\kappa|/(2N)$, and restoring $c_{N,D}$ changes it by $N^{-1}\log c_{N,D}$. Both changes tend to zero, by \eqref{eq:normalization-subexp}. Combining this calculation with \eqref{eq:enriched-Laplace-error} proves
\begin{equation}
\label{eq:enriched-restricted-Laplace}
\lim_{N\to\infty}\frac1N\E\log Z_N^{\geq\kappa}(t,q)=\sup_{a\in\mcl A_\kappa}\{-f_a(t,q_a)+\log\det a\}.
\end{equation}

\smallskip
\noindent\emph{Step 3: Concentration and $\log(Z_{N}^{\geq\kappa}+Z_{N}^{<\kappa})$.} To pass from the two restricted integrals to their sum, we need concentration of both logarithms. The two reference measures are $P_N$ restricted to $\det(N^{-1}\sigma\sigma^\intercal)\ge\kappa$ and to its complement. Both have positive mass when $D\ge2$ and $N>D$. Indeed, a relative neighborhood of $I_D$ has $\det b>\kappa$, while the matrix with upper-left block $\left(\begin{smallmatrix}1&r\\r&1\end{smallmatrix}\right)$ and remaining diagonal block $I_{D-2}$ (omitted when $D=2$) has determinant $1-r^2\in(0,\kappa)$ for $r\in(\sqrt{1-\kappa},1)$. A relative neighborhood of this matrix also has determinant in $(0,\kappa)$, and the density \eqref{eq:matrix-density} is positive on both neighborhoods. Normalize each restricted measure to a probability measure. Under either normalized reference measure, $|\sigma|^2=ND$, and the variance of the Gaussian exponent is bounded by
\begin{equation*}
2N\left(t\max_{a\in\mcl A}\xi(a^2)+D|q(1)|_{\mathrm{op}}\right).
\end{equation*}
The correction terms and $h\cdot\sum_i\sigma_i$ are bounded deterministic functions of $\sigma$. Apply Lemma~\ref{lem:total-concentration-continuous} with $\Omega$ equal to either restriction, $G_0(\sigma)=\sqrt{2t}H_N(\sigma)$, $v_i(\sigma)=\sqrt2\sigma_i$, and covariance paths $q_i=q$. Its deterministic function $A:\Omega\to\R$ is $A(\sigma):=h\cdot\sum_i\sigma_i-Nt\xi(N^{-1}\sigma\sigma^{\intercal})-q(1)\cdot\sigma\sigma^{\intercal}$. Each $v_i$ and $A$ is bounded for fixed $N$, so the lemma applies. Restoring the deterministic mass of the reference measure adds a constant to its logarithm and does not affect its centered value. Consequently, we have
\begin{equation}
\label{eq:enriched-restricted-concentration}
\E\left|N^{-1}\log Z_{N}^{\geq\kappa}-N^{-1}\E\log Z_{N}^{\geq\kappa}\right|+\E\left|N^{-1}\log Z_{N}^{<\kappa}-N^{-1}\E\log Z_{N}^{<\kappa}\right|\le C_{t,q}N^{-1/2}.
\end{equation}
Here the arguments $(t,q)$ are suppressed only in this display.

For integrable random variables $X,Y$, we have
\begin{equation*}
\left|\E\max(X,Y)-\max(\E X,\E Y)\right|\le\E|X-\E X|+\E|Y-\E Y|.
\end{equation*}
Apply this inequality with $X=\log Z_{N}^{\geq\kappa}(t,q)$ and $Y=\log Z_{N}^{<\kappa}(t,q)$. Since $\max(X,Y)\le\log(e^X+e^Y)\le\max(X,Y)+\log2$, dividing by $N$ and using \eqref{eq:enriched-restricted-concentration} gives
\begin{equation}
\label{eq:two-part-log-sum}
\frac1N\E\log\Ll(Z_{N}^{\geq\kappa}(t,q)+Z_{N}^{<\kappa}(t,q)\Rr)
=\max\left\{\frac1N\E\log Z_{N}^{\geq\kappa}(t,q),\frac1N\E\log Z_{N}^{<\kappa}(t,q)\right\}+O_{t,q}(N^{-1/2}).
\end{equation}

\smallskip
\noindent\emph{Step 4: Choosing $\kappa$ and attaining the infimum.} Jensen's inequality and \eqref{eq:enriched-annealed-one} give $\bar F_{N,a}(t,q_a)\ge-|ah|\ge-\sqrt D\,|h|$, so Corollary~\ref{cor:fixed-a-limit} implies
\begin{equation}
\label{eq:enriched-fiber-nonnegative}
f_a(t,q_a)\ge-|ah|\ge-\sqrt D\,|h|,\qquad\forall a\in\mcl A^\circ.
\end{equation}
In particular, $-f_a(t,q_a)+\log\det a\le\sqrt D\,|h|+\log\det a$, whereas its value at $a=I_D$ is the finite number $-f_{I_D}(t,q)$. Choose $\kappa\in(0,1)$ so that
\begin{equation}
\label{eq:enriched-kappa-choice}
\sqrt D\,|h|+\frac12\log\kappa<-f_{I_D}(t,q)-1.
\end{equation}
Since $I_D\in\mcl A_\kappa$, the two restricted logarithms satisfy
\begin{align*}
\lim_{N\to\infty}\frac1N\E\log Z_N^{\geq\kappa}(t,q)&\stackrel{\eqref{eq:enriched-restricted-Laplace}}{\ge}-f_{I_D}(t,q),\\
\limsup_{N\to\infty}\frac1N\E\log Z_N^{<\kappa}(t,q)&\stackrel{\eqref{eq:enriched-boundary-estimate}}{\le}\sqrt D\,|h|+\frac12\log\kappa\stackrel{\eqref{eq:enriched-kappa-choice}}{<}-f_{I_D}(t,q)-1.
\end{align*}
The first sequence converges to a number at least $-f_{I_D}(t,q)$, while the upper limit of the second is less than $-f_{I_D}(t,q)-1$. Thus the maximum in \eqref{eq:two-part-log-sum} is the term with $Z_N^{\geq\kappa}$ for large $N$. Using $Z_N=Z_N^{\geq\kappa}+Z_N^{<\kappa}$ and \eqref{eq:enriched-restricted-Laplace}, we obtain
\begin{equation}\label{e.ElogZ=supvarphi}
\lim_{N\to\infty}\frac1N\E\log Z_N(t,q)=\sup_{a\in\mcl A_\kappa}\{-f_a(t,q_a)+\log\det a\}.
\end{equation}

It remains to compare the supremum on $\mcl A_\kappa$ with that on $\mcl A^\circ$. For $a\in\mcl A^\circ$ with $\det(a^2)<\kappa$, the lower bound on $f_a$ in \eqref{eq:enriched-fiber-nonnegative} gives
\begin{equation*}
-f_a(t,q_a)+\log\det a\le\sqrt D\,|h|+\frac12\log\det(a^2)<\sqrt D\,|h|+\frac12\log\kappa.
\end{equation*}
Taking the supremum over these $a$ and using \eqref{eq:enriched-kappa-choice} and $I_D\in\mcl A_\kappa$ yields
\begin{align*}
\sup_{\substack{a\in\mcl A^\circ\\\det(a^2)<\kappa}}\{-f_a(t,q_a)+\log\det a\}\le\sqrt D\,|h|+\frac12\log\kappa<-f_{I_D}(t,q)\le\sup_{a\in\mcl A_\kappa}\{-f_a(t,q_a)+\log\det a\}.
\end{align*}
Every $a\in\mcl A^\circ$ belongs either to $\mcl A_\kappa$ or to $\{a\in\mcl A^\circ:\det(a^2)<\kappa\}$, so
\begin{equation*}
\sup_{a\in\mcl A^\circ}\{-f_a(t,q_a)+\log\det a\}=\sup_{a\in\mcl A_\kappa}\{-f_a(t,q_a)+\log\det a\}.
\end{equation*}
Using~\eqref{e.ElogZ=supvarphi} and recalling $\bar F_N(t,q)=-N^{-1}\E\log Z_N(t,q)$ proves \eqref{eq:global-enriched-limit}. Since $a\mapsto-f_a(t,q_a)+\log\det a$ is continuous on the compact set $\mcl A_\kappa\subset\mcl A^\circ$, it attains its maximum there. The same optimizer also minimizes \eqref{eq:global-enriched-limit}.

\smallskip
\noindent\emph{Step 5: Extension to $q\in\mcl Q_2$.} For finite-step $q,q'$, the argument of Lemma~\ref{lem:fixed-path-Lipschitz}, using \cite[Proposition~5.1]{chen2025free} and $|N^{-1}\sigma\sigma'^{\intercal}|\le D$ on $S_N^D$, gives
\begin{equation}
\label{eq:global-enriched-q-Lipschitz}
|\bar F_N(t,q)-\bar F_N(t,q')|\le D\int_0^1|q(s)-q'(s)|\,\d s.
\end{equation}
This estimate defines a unique $L^1$-continuous extension to $\mcl Q_1$ and remains valid for every $q,q'\in\mcl Q_1$. Moreover, $|a|_{\mathrm{op}}\le\sqrt D$ for $a\in\mcl A$, so \eqref{eq:qa-intro} yields
\begin{equation}
\label{eq:qa-global-Lipschitz}
\sup_{a\in\mcl A}\int_0^1|q_a(s)-q'_a(s)|\,\d s\le D\int_0^1|q(s)-q'(s)|\,\d s.
\end{equation}
For the given $q\in\mcl Q_2$, set $q_n(s):=q(\min\{s,1-n^{-1}\})$ for $n\in\{2,3,\ldots\}$. Then $q_n\in\mcl Q_\infty$ and $\lim_{n\to\infty}|q_n-q|_{L^1}=0$, as in the proof of Corollary~\ref{cor:fixed-a-limit}. Equations \eqref{eq:Hopf-fixed-path-Lipschitz} and \eqref{eq:qa-global-Lipschitz} give
\begin{equation}
\label{eq:enriched-Hopf-path-approximation}
\sup_{a\in\mcl A^\circ}|f_a(t,(q_n)_a)-f_a(t,q_a)|\le D^{3/2}|q_n-q|_{L^1}.
\end{equation}
Equation \eqref{eq:enriched-fiber-nonnegative} and \eqref{eq:enriched-Hopf-path-approximation} bound both objectives below uniformly in $a$, while testing $a=I_D$ gives finite upper bounds for their infima. Since the term $-\log\det a$ is unchanged, the two infima satisfy
\begin{equation}
\label{eq:enriched-infimum-path-approximation}
\left|\inf_{a\in\mcl A^\circ}\{f_a(t,(q_n)_a)-\log\det a\}-\inf_{a\in\mcl A^\circ}\{f_a(t,q_a)-\log\det a\}\right|\le D^{3/2}|q_n-q|_{L^1}.
\end{equation}
Combining \eqref{eq:enriched-infimum-path-approximation} with \eqref{eq:global-enriched-q-Lipschitz} and the formula already proved at $q_n$, we obtain
\begin{equation*}
\limsup_{N\to\infty}\left|\bar F_N(t,q)-\inf_{a\in\mcl A^\circ}\{f_a(t,q_a)-\log\det a\}\right|\le(D+D^{3/2})|q_n-q|_{L^1}.
\end{equation*}
Let $n\to\infty$ to prove \eqref{eq:global-enriched-limit}. For attainment, each function $a\mapsto f_a(t,(q_n)_a)$ is continuous on $\mcl A_\kappa$ by \eqref{eq:enriched-limit-a-modulus}, and \eqref{eq:enriched-Hopf-path-approximation} gives their uniform convergence to $a\mapsto f_a(t,q_a)$ there. The limit is therefore continuous. Passing to the limit in \eqref{eq:enriched-fiber-nonnegative} gives the same lower bound for this $q$. Choose $\kappa$ using \eqref{eq:enriched-kappa-choice}, now with $f_{I_D}(t,q)$. Step 4 then excludes $\det(a^2)<\kappa$ from the minimization, and the continuous objective attains its minimum on the compact set $\mcl A_\kappa\subset\mcl A^\circ$.
\end{proof}

\subsection{Formula for \texorpdfstring{$\mathsf F_N$}{\unichar{"1D5A5}\_N}}
\label{s:positive-outer}

The proof follows the same steps as Section~\ref{s:global-enriched}: local uniform convergence in $a$, disintegration by \eqref{eq:measure-disintegration}, and an estimate for $\det(a^2)<\kappa$. Here the partition function has no correction or cascade fields. Corollary~\ref{cor:fixed-a-limit} supplies the pointwise limit, and Lemma~\ref{lem:total-concentration-continuous} supplies the concentration estimate. Recall $f_a$ from \eqref{eq:Hopf-fixed-a-intro}. For $a\in\mcl A^\circ$, define the random free energy by
\begin{equation}
\mathsf F_{N,a}^{\omega}:=-\frac1N\log\int_{\mathbb U_N} \exp\Ll(H_N(au)+h\cdot\sum_{i=1}^Nau_i\Rr)\nu_N(\mathrm d u) \quad\text{and thus}\quad \E \mathsf F_{N,a}^{\omega}\stackrel{\eqref{eq:fixed-a-free-energy-intro}}{=}\mathsf F_{N,a}.
\label{eq:random-fiber-free-energy}
\end{equation}
Setting $t=1/2$ and $q=0$ in \eqref{eq:enriched-fixed-a-H} gives the relation
\begin{equation*}
\mathsf F_{N,a}^{\omega}=F_{N,a}\Ll(\tfrac12,0\Rr)-\tfrac12\xi(a^2).
\end{equation*}
Taking expectations recovers \eqref{eq:xia-I}. Recall $\mcl A_\kappa$ from \eqref{eq:Skappa-enriched}. We use the same interaction tensors for every $a$, as in \eqref{eq:fiber-cross-covariance}. Continuity of $\xi$ and compactness of $\mcl A$ give $\max_{a\in\mcl A}\xi(a^2)<\infty$. 

\begin{lemma}
\label{lem:a-regularity}
For every $\kappa\in(0,1)$, there exists $C_{h,\kappa}<+\infty$, independent of $N$, such that, for every $N\in\N$ with $N\ge D$ and every $\delta\in(0,1)$, we have
\begin{equation}
\label{eq:a-random-modulus}
    \E\sup_{\substack{a,a'\in \mcl A_\kappa\\|a-a'|\le\delta}}
    |\mathsf F_{N,a}^{\omega}-\mathsf F_{N,a'}^{\omega}|
    \le C_{h,\kappa}\delta.
\end{equation}
Moreover, there exists $C<\infty$, independent of $N$, such that
\begin{equation}
\label{eq:a-concentration}
\sup_{a\in\mcl A^\circ}\E|\mathsf F_{N,a}^{\omega}-\mathsf F_{N,a}|\le CN^{-1/2}.
\end{equation}
\end{lemma}

\begin{proof}
\smallskip
\noindent\emph{Step 1: Proof of \eqref{eq:a-random-modulus}.} Recall $|N^{-1}\sum_i u_i|\le1$ and the comparison of exponents used to prove \eqref{eq:matrix-factor-value-comparison}. Since $|a|=|a'|=\sqrt D$, the segment from $a$ to $a'$ lies in $\{|w|\le\sqrt D\}$. The same integration along this segment as in the proof of \eqref{eq:matrix-factor-modulus} gives
\begin{equation*}
\sup_{\substack{a,a'\in \mcl A_\kappa\\|a-a'|\le\delta}}|\mathsf F_{N,a}^{\omega}-\mathsf F_{N,a'}^{\omega}|\le\delta\bigg(\frac1N\sup_{\substack{|w|\le\sqrt D,\ |b|\le1\\u\in\mathbb U_N}}|(b\cdot\nabla_w)H_N(wu)|+|h|\bigg).
\end{equation*}
Taking expectation and applying \eqref{eq:matrix-interaction-derivative-sup} proves \eqref{eq:a-random-modulus}.

\smallskip
\noindent\emph{Step 2: Proof of \eqref{eq:a-concentration}.} Use the case $J=\varnothing$ of Lemma~\ref{lem:total-concentration-continuous}, whose partition function reduces to $\int_\Omega e^{G_0(x)+A(x)}\mu(\mathrm d x)$. Take $\Omega=\mathbb U_N$, $\mu=\nu_N$, $G_0(u)=H_N(au)$, and $A(u)=ah\cdot\sum_i u_i$. This is the partition function in \eqref{eq:random-fiber-free-energy}. For each fixed $N$, $A$ is bounded, and the variance parameter in \eqref{eq:can-concentration-variance} is
\begin{equation*}
V=\sup_{u\in\mathbb U_N}\E H_N(au)^2=N\xi(a^2)\le N\max_{a\in\mcl A}\xi(a^2).
\end{equation*}
Dividing \eqref{eq:can-total-concentration} by $N$ gives
\begin{equation}
\label{eq:a-mean-concentration}
\E|\mathsf F_{N,a}^{\omega}-\mathsf F_{N,a}|\le\frac C N\sqrt{1+N\xi(a^2)}.
\end{equation}
Taking the supremum over $a\in\mcl A^\circ$ and using $\max_{a\in\mcl A}\xi(a^2)<\infty$ proves \eqref{eq:a-concentration}. This application uses the full Gaussian field $H_N$ and requires no further truncation. The bound also holds when $\max_{a\in\mcl A}\xi(a^2)=0$, so no separate case is needed.
\end{proof}

For later use, define the limit identified in Corollary~\ref{cor:fixed-a-limit} by
\begin{equation*}
\mathsf F_{\infty,a}:=\lim_{N\to\infty}\mathsf F_{N,a}=f_a\Ll(\tfrac12,0\Rr)-\tfrac12\xi(a^2).
\end{equation*}
We next strengthen this pointwise convergence to the uniform convergence needed to integrate over $a$.

\begin{lemma}
\label{lem:local-uniform-a} 
For every $\kappa\in(0,1)$, we have 
\begin{equation}
\label{eq:deterministic-local-uniform}
\lim_{N\to\infty}\sup_{a\in \mcl A_\kappa}|\mathsf F_{N,a}-\mathsf F_{\infty,a}|=0\qquad\text{and}\qquad
\lim_{N\to\infty}\E\sup_{a\in \mcl A_\kappa}|\mathsf F_{N,a}^{\omega}-\mathsf F_{\infty,a}|=0.
\end{equation}
\end{lemma}

\begin{proof}
\smallskip
\noindent\emph{Step 1: Continuity in $a$.} For $a,a'\in\mcl A_\kappa$ with $|a-a'|<1$, the triangle inequality for the expectation and \eqref{eq:a-random-modulus} give the first inequality below. Passing to the pointwise limit gives the second:
\begin{equation*}
|\mathsf F_{N,a}-\mathsf F_{N,a'}|\le C_{h,\kappa}|a-a'|\qquad\text{and}\qquad|\mathsf F_{\infty,a}-\mathsf F_{\infty,a'}|\le C_{h,\kappa}|a-a'|.
\end{equation*}
Thus $\mathsf F_{\infty,\cdot}$ is continuous on $\mcl A_\kappa$.

\smallskip
\noindent\emph{Step 2: The two uniform limits.} Apply the finite-net argument in Step 2 of Lemma~\ref{lem:enriched-local-uniform}, using \eqref{eq:a-random-modulus}, \eqref{eq:a-concentration}, and the continuity estimate from Step 1. For $\delta\in(0,1)$ and a finite $\delta$-net $a_1,\ldots,a_m$ of $\mcl A_\kappa$, that argument gives
\begin{equation*}
\E\sup_{a\in\mcl A_\kappa}|\mathsf F_{N,a}^{\omega}-\mathsf F_{\infty,a}|\le2C_{h,\kappa}\delta+\max_{j\le m}|\mathsf F_{N,a_j}-\mathsf F_{\infty,a_j}|+CmN^{-1/2}.
\end{equation*}
For fixed $\delta$, both $m$ and the net points are independent of $N$. Pointwise convergence makes the last two terms tend to zero as $N\to\infty$. Letting $\delta\downarrow0$ proves the second identity in \eqref{eq:deterministic-local-uniform}. The first follows from $\sup_{a\in\mcl A_\kappa}|\mathsf F_{N,a}-\mathsf F_{\infty,a}|\le\E\sup_{a\in\mcl A_\kappa}|\mathsf F_{N,a}^{\omega}-\mathsf F_{\infty,a}|$.
\end{proof}

\begin{proposition}
\label{prop:outer-a-free-energy}
We have 
\begin{equation}
\label{eq:outer-a-free-energy}
\lim_{N\to\infty}\mathsf F_N =\inf_{a\in\mcl A^\circ}\left\{\mathsf F_{\infty,a}-\log\det a\right\}.
\end{equation}
The infimum is attained in $\mcl A^\circ$. 
\end{proposition}

\begin{proof}
When $D=1$, the set $\mcl A^\circ$ consists of $I_1$ and the result is the second identity in \eqref{eq:fixed-a-negative-limit}. Assume $D\ge2$.

\smallskip
\noindent\emph{Step 1: The integral over $\mcl B_\kappa$.} Write the unrestricted partition function as
\begin{equation*}
Z_N:=\int_{S_N^D}\exp\Ll(H_N(\sigma)+h\cdot\sum_{i=1}^N\sigma_i\Rr)P_N(\mathrm d\sigma).
\end{equation*}
Then $\mathsf F_N=-N^{-1}\E\log Z_N$ by \eqref{eq:physical-free-energy}. Apply \eqref{eq:measure-disintegration} with $G(\sigma)=\exp(H_N(\sigma)+h\cdot\sum_i\sigma_i)$. For $b\in\mcl B^\circ$, the inner integral is identified by \eqref{eq:random-fiber-free-energy}:
\begin{equation*}
\int_{\mathbb U_N}\exp\Ll(H_N(b^{1/2}u)+h\cdot\sum_{i=1}^Nb^{1/2}u_i\Rr)\nu_N(\mathrm d u)=\exp\Ll(-N\mathsf F_{N,b^{1/2}}^\omega\Rr).
\end{equation*}
Substituting this identity into \eqref{eq:measure-disintegration} gives almost surely
\begin{equation}
\label{eq:partition-disintegration}
Z_N=c_{N,D}\int_{\mcl B^\circ}\exp\Ll(-N\mathsf F_{N,b^{1/2}}^{\omega}\Rr)(\det b)^{(N-D-1)/2}\,\mathrm d b_{\mathrm{off}}.
\end{equation}
Fix $\kappa\in(0,1)$, and let $Z_{N}^{\geq\kappa}$ and $Z_{N}^{<\kappa}$ be the integrals over $\mcl B_\kappa$ and $\{b\in\mcl B^\circ:\det b<\kappa\}$, respectively, including the factor $c_{N,D}$. Thus $Z_N=Z_{N}^{\geq\kappa}+Z_{N}^{<\kappa}$. For $b\in\mcl B_\kappa$, we have $b^{1/2}\in\mcl A_\kappa$ by \eqref{eq:Skappa-enriched}. For this comparison, put $r_N:=\sup_{a\in\mcl A_\kappa}|\mathsf F_{N,a}^\omega-\mathsf F_{\infty,a}|$. For every $b\in\mcl B_\kappa$, the inequality $|\mathsf F_{N,b^{1/2}}^\omega-\mathsf F_{\infty,b^{1/2}}|\le r_N$ gives
\begin{equation*}
e^{-Nr_N}e^{-N\mathsf F_{\infty,b^{1/2}}}\le e^{-N\mathsf F_{N,b^{1/2}}^\omega}\le e^{Nr_N}e^{-N\mathsf F_{\infty,b^{1/2}}}.
\end{equation*}
Integrate against $c_{N,D}(\det b)^{(N-D-1)/2}\,\mathrm d b_{\mathrm{off}}$ over $\mcl B_\kappa$. The middle integral is $Z_N^{\geq\kappa}$ by \eqref{eq:partition-disintegration}, and division by the positive deterministic integral gives
\begin{equation*}
e^{-Nr_N}\le\frac{Z_N^{\geq\kappa}}{c_{N,D}\int_{\mcl B_\kappa}e^{-N\mathsf F_{\infty,b^{1/2}}}(\det b)^{(N-D-1)/2}\,\mathrm d b_{\mathrm{off}}}\le e^{Nr_N}.
\end{equation*}
Taking logarithms, dividing by $N$, and substituting the definition of $r_N$, we obtain
\begin{align}
\label{eq:positive-Laplace-error}
\left|\frac1N\log Z_N^{\geq\kappa}-\frac1N\log\left(c_{N,D}\int_{\mcl B_\kappa}e^{-N\mathsf F_{\infty,b^{1/2}}}(\det b)^{(N-D-1)/2}\,\mathrm d b_{\mathrm{off}}\right)\right|\le\sup_{a\in\mcl A_\kappa}|\mathsf F_{N,a}^\omega-\mathsf F_{\infty,a}|.
\end{align}
Taking expectation in \eqref{eq:positive-Laplace-error}, the right-hand side tends to zero by the second identity in \eqref{eq:deterministic-local-uniform}. The determinant weight satisfies
\begin{equation*}
(\det b)^{(N-D-1)/2}=\exp\Ll(\tfrac N2\log\det b\Rr)(\det b)^{-(D+1)/2},
\end{equation*}
so the resulting deterministic integral has exponent $N(-\mathsf F_{\infty,b^{1/2}}+\tfrac12\log\det b)$. Its remaining factor lies between $1$ and $\kappa^{-(D+1)/2}$. The function $a\mapsto\mathsf F_{\infty,a}$ is continuous on $\mcl A_\kappa$ by Lemma~\ref{lem:local-uniform-a}. The proof of \eqref{eq:compact-Laplace-upper}--\eqref{eq:compact-Laplace-lower} used only continuity of the objective on $\mcl A_\kappa$ and positive volume in every relative neighborhood of a maximizer in $\mcl B_\kappa$. It therefore applies with $a\mapsto-\mathsf F_{\infty,a}+\log\det a$. The remaining factor changes the normalized logarithm by at most $(D+1)|\log\kappa|/(2N)$, and $N^{-1}\log c_{N,D}$ tends to zero by \eqref{eq:normalization-subexp}. We obtain
\begin{equation}
\label{eq:restricted-Laplace}
\lim_{N\to\infty}\frac1N\E\log Z_{N}^{\geq\kappa}=\sup_{a\in \mcl A_\kappa}\left\{-\mathsf F_{\infty,a}+\log\det a\right\}.
\end{equation}

\smallskip
\noindent\emph{Step 2: Bounds for $\det b<\kappa$ and for $-\mathsf F_{\infty,a}$.} The variance of $H_N(au)$ is $N\xi(a^2)$, so we have $\E e^{H_N(au)}=\exp\Ll(\frac N2\xi(a^2)\Rr)\le \exp\Ll(\frac N2\max_{a\in\mcl A}\xi(a^2)\Rr)$.
Also, $|a|_{\mathrm{op}}\le\sqrt D$ and $|N^{-1}\sum_i u_i|\le1$ give $\left|h\cdot\sum_{i=1}^Nau_i\right|\le N\sqrt D\,|h|$.
Taking expectations in \eqref{eq:partition-disintegration} restricted to $\det b<\kappa$, we therefore obtain, for $N>D+1$,
\begin{equation*}
\E Z_{N}^{<\kappa}\le c_{N,D}\operatorname{vol}(\mcl B)\exp\Ll(N\Ll(\frac12\max_{a\in\mcl A}\xi(a^2)+\sqrt D\,|h|\Rr)\Rr)\kappa^{(N-D-1)/2}.
\end{equation*}
Jensen's inequality gives
\begin{equation*}
\frac1N\E\log Z_{N}^{<\kappa}\le\frac{\log c_{N,D}+\log\operatorname{vol}(\mcl B)}N+\frac12\max_{a\in\mcl A}\xi(a^2)+\sqrt D\,|h|+\frac{N-D-1}{2N}\log\kappa.
\end{equation*}
For fixed $\kappa$, take the upper limit as $N\to\infty$ and use \eqref{eq:normalization-subexp}. This yields
\begin{equation}
\label{eq:positive-boundary-rate}
\limsup_{N\to\infty}\frac1N\E\log Z_{N}^{<\kappa}\le\frac12\max_{a\in\mcl A}\xi(a^2)+\sqrt D\,|h|+\frac12\log\kappa.
\end{equation}

We also need a bound on the limiting objective for these values of $a$. By Jensen's inequality and Tonelli's theorem, the Gaussian expectation computed in Step 2 gives
\begin{equation*}
-\mathsf F_{N,a}\le\frac1N\log\E\int_{\mathbb U_N}e^{H_N(au)+ah\cdot\sum_i u_i}\nu_N(\mathrm d u)\le\frac12\xi(a^2)+\sqrt D\,|h|.
\end{equation*}
For a lower bound, convexity of the log partition function in the centered Gaussian field gives
\begin{equation*}
-\mathsf F_{N,a}\ge\frac1N\log\int_{\mathbb U_N}\exp\Ll(ah\cdot\sum_{i=1}^Nu_i\Rr)\nu_N(\mathrm d u)\ge0.
\end{equation*}
The last inequality is Jensen's inequality, since $\int u_i\nu_N(\mathrm d u)=0$. Passing to the limit proves
\begin{equation}
\label{eq:fiber-annealed-bounds}
0\le-\mathsf F_{\infty,a}\le\frac12\xi(a^2)+\sqrt D\,|h|\le\frac12\max_{a\in\mcl A}\xi(a^2)+\sqrt D\,|h|.
\end{equation}

\smallskip
\noindent\emph{Step 3: Choosing $\kappa$ and computing the unrestricted limit.} Choose $\kappa\in(0,1)$ so that
\begin{equation}
\label{eq:positive-kappa-choice}
\frac12\max_{a\in\mcl A}\xi(a^2)+\sqrt D\,|h|+\frac12\log\kappa<0.
\end{equation}
Since $I_D\in \mcl A_\kappa$ and $-\mathsf F_{\infty,I_D}\ge0$, the limit in \eqref{eq:restricted-Laplace} is nonnegative. By \eqref{eq:positive-boundary-rate}, the complementary integral has a strictly negative limiting upper rate.

The two restrictions of $P_N$ have positive masses, as proved for \eqref{eq:enriched-restricted-concentration}. Normalize each to a probability measure. Apply Lemma~\ref{lem:total-concentration-continuous} with either normalized measure, $G_0(\sigma)=H_N(\sigma)$, $A(\sigma)=h\cdot\sum_i\sigma_i$, and $J=\varnothing$. Since $N^{-1}\sigma\sigma^{\intercal}=a^2$ for some $a\in\mcl A$, the variance parameter satisfies $V\le N\max_{a\in\mcl A}\xi(a^2)$. Restoring the deterministic mass adds a constant to the logarithm and leaves its centered value unchanged. Consequently, we have
\begin{equation}
\label{eq:positive-restricted-concentration}
\E\left|N^{-1}\log Z_{N}^{\geq\kappa}-N^{-1}\E\log Z_{N}^{\geq\kappa}\right|+\E\left|N^{-1}\log Z_{N}^{<\kappa}-N^{-1}\E\log Z_{N}^{<\kappa}\right|\le CN^{-1/2}.
\end{equation}
By \eqref{eq:restricted-Laplace}, $N^{-1}\E\log Z_N^{\geq\kappa}$ converges to a nonnegative value. By \eqref{eq:positive-boundary-rate} and \eqref{eq:positive-kappa-choice}, $N^{-1}\E\log Z_N^{<\kappa}$ has strictly negative upper limit. Hence, for all sufficiently large $N$, we have
\begin{equation*}
\frac1N\E\log Z_N^{<\kappa}<\frac1N\E\log Z_N^{\geq\kappa}.
\end{equation*}
For these $N$, we have $\max\{\E\log Z_N^{\geq\kappa},\E\log Z_N^{<\kappa}\}=\E\log Z_N^{\geq\kappa}$. The identity $Z_N=Z_N^{\geq\kappa}+Z_N^{<\kappa}$ and the bound $x+y\le2\max\{x,y\}$ for $x,y>0$ give the first inequality below. For the second, subtract the expectation from each restricted logarithm and use this ordering of their expectations:
\begin{align*}
\log Z_N&\le\max\{\log Z_N^{\geq\kappa},\log Z_N^{<\kappa}\}+\log2\\
&\le\E\log Z_N^{\geq\kappa}+|\log Z_N^{\geq\kappa}-\E\log Z_N^{\geq\kappa}|+|\log Z_N^{<\kappa}-\E\log Z_N^{<\kappa}|+\log2.
\end{align*}
Also, $Z_N\ge Z_N^{\geq\kappa}$. Taking expectations, dividing by $N$, and recalling $\mathsf F_N=-N^{-1}\E\log Z_N$, we therefore obtain
\begin{align*}
\left|\mathsf F_N+\frac1N\E\log Z_N^{\geq\kappa}\right|
=\frac1N\E\log Z_N-\frac1N\E\log Z_N^{\geq\kappa}\stackrel{\eqref{eq:positive-restricted-concentration}}{\le}CN^{-1/2}+\frac{\log2}N.
\end{align*}
The right-hand side tends to zero as $N\to\infty$. Taking the limit of the restricted logarithm by \eqref{eq:restricted-Laplace}, we obtain
\begin{equation*}
\lim_{N\to\infty}\mathsf F_N=-\sup_{a\in \mcl A_\kappa}\{-\mathsf F_{\infty,a}+\log\det a\}=\inf_{a\in \mcl A_\kappa}\{\mathsf F_{\infty,a}-\log\det a\}.
\end{equation*}
\smallskip
\noindent\emph{Step 4: The infimum over $\mcl A^\circ$ and its attainment.} For the variational problem, \eqref{eq:fiber-annealed-bounds} and \eqref{eq:positive-kappa-choice} give
\begin{equation*}
\inf_{\substack{a\in\mcl A^\circ\\\det(a^2)<\kappa}}\{\mathsf F_{\infty,a}-\log\det a\}\ge-\frac12\max_{a\in\mcl A}\xi(a^2)-\sqrt D\,|h|-\frac12\log\kappa>0\ge\mathsf F_{\infty,I_D}.
\end{equation*}
Thus extending the infimum from $\mcl A_\kappa$ to $\mcl A^\circ$ does not change its value, proving \eqref{eq:outer-a-free-energy}. By Lemma~\ref{lem:local-uniform-a}, $a\mapsto\mathsf F_{\infty,a}$ is continuous on $\mcl A_\kappa$. The term $-\log\det a$ is continuous there as well, since $\det(a^2)\ge\kappa>0$. Thus the objective attains its minimum on the compact set $\mcl A_\kappa\subset\mcl A^\circ$. 
\end{proof}

\begin{proof}[Proof of Theorem~\ref{thm:main}]
Corollary~\ref{cor:fixed-a-limit} identifies $\mathsf F_{\infty,a}=f_a(\tfrac12,0)-\tfrac12\xi(a^2)$ for every $a\in\mcl A^\circ$. Substituting this identity into Proposition~\ref{prop:outer-a-free-energy}, we obtain
\begin{equation*}
\lim_{N\to\infty}\mathsf F_N=\inf_{a\in\mcl A^\circ}\{\mathsf F_{\infty,a}-\log\det a\}=\inf_{a\in\mcl A^\circ}\Big\{f_a\Big(\tfrac12,0\Big)-\tfrac12\xi(a^2)-\log\det a\Big\}.
\end{equation*}
Proposition~\ref{prop:outer-a-free-energy} also proves attainment of this infimum. The second identity in Theorem~\ref{thm:main}, together with attainment of its infimum, is Proposition~\ref{prop:global-enriched}.
\end{proof}

\bibliographystyle{abbrv}

\end{document}